\documentclass[12pt]{amsart}

\usepackage{lahn}
\usepackage{ALLM}
\usepackage{fullpage}

\title{Unmarked simple length spectral rigidity for covers}

\author{Tarik Aougab}
\address{Department of Mathematics \\ Haverford College \\ Haverford, PA 19041}
\email{taougab@haverford.edu}

\author{Max Lahn}
\address{Department of Mathematics \\ University of Michigan \\ Ann Arbor, MI 48109}
\email{maxlahn@umich.edu}

\author{Marissa Loving}
\address{Department of Mathematics \\ University of Wisconsin \\ Madison, WI 53706}
\email{mloving2@wisc.edu}

\author{Nicholas Miller}
\address{Department of Mathematics \\ University of Oklahoma \\ Norman, OK 73019}
\email{nickmbmiller@ou.edu}

\begin{document}

\begin{abstract}

We prove that every closed orientable surface $ S $ of negative Euler characteristic admits a pair of finite-degree covers which are length isospectral over $ S $ but generically not simple length isospectral over $ S $. To do this, we first characterize when two finite-degree covers of a connected, orientable surface of negative Euler characteristic are isomorphic in terms of which curves have simple elevations. We also construct hyperbolic surfaces $ X $ and $ Y $ with the same full unmarked length spectrum but so that for each $ k $, the sets of lengths associated to curves with at most $ k $ self-intersections differ.

\end{abstract}

\maketitle

\section{Introduction}

Otal \cite{Ota90} showed that any negatively curved metric on $ S $ is determined up to isotopy by its \emph{marked length spectrum}, the function associating each closed curve to its infimal length. More generally, one can ask this question for different sets of closed curves and/or different families of metrics (see Fricke, Wolpert \cite{Wol79}, and Hamenst\"{a}dt's \cite{Ham03} work on hyperbolic metrics and Duchin--Leininger--Rafi \cite{DLR10} and Loving's \cite{Loving20} work on flat metrics).

The roots of this paper lie in a family of interrelated questions distinguished from the above work by the absence of markings. Concretely, consider a $ \Mod \of{ S } $-invariant set $ \mathcal{ G } $ of closed curves on $ S $. For example, we will consider for each $ k \in \Z_{\ge 0} \cup \set{ \infty } $ the set $ \mathcal{ G }_{ k } $ of curves whose self-intersection number is at most $ k $. The \emph{unmarked $ \mathcal{ G } $-length spectrum} of a hyperbolic metric $ \mu $ on $ S $ is the multi-set of infimal lengths of closed curves in $ \mathcal{ G } $ in the metric $ \mu $. Particularly relevant examples are the \emph{unmarked length spectrum} and the \emph{unmarked simple length spectrum}, which arise in the cases when $ \mathcal{ G } = \mathcal{ G }_{ \infty } $ is the set of all closed curves on $ S $ and $ \mathcal{ G } = \mathcal{ G }_{ 0 } $ is the set of simple curves on $ S $, respectively. 

\begin{question+} \label{que:UnmarkedLengthSpectralRigidity}

For which $ \mathcal{ G } $ is a hyperbolic metric $ \mu $ on $ S $ determined up to isometry by its unmarked $ \mathcal{ G } $-length spectrum?

\end{question+}

Unlike in the case of marked spectra, it is well-known that many classes of metrics on $ S $ are not unmarked length spectrally rigid. Indeed, examples of unmarked length isospectral but non-isometric complete hyperbolic surface metrics were first constructed by Vign\'{e}ras \cite{Vig80} in the arithmetic setting. Since then, a celebrated construction of Sunada \cite{Sun85} yields a robust method for producing unmarked length isospectral but non-isometric surface metrics arising as pullbacks of metrics along covers of a fixed surface.

On the other hand, relatively little is known about the question of unmarked simple length spectral rigidity for hyperbolic metrics. Baik--Choi--Kim \cite{BCK21} recently showed that, up to excluding the loci of some length identities, the family of hyperbolic metrics with the property that lengths in the simple length spectrum have multiplicity $ 1 $ (that is, for each positive real number $ L $, there is at most one simple closed geodesic of length $ L $) is simple length spectrally rigid. This result motivates our examination of pullbacks of hyperbolic metrics along finite-degree covers, since high multiplicity is common in this regime.

Explicitly, fix for the remainder of this article two finite-degree covers $ p \colon X \to S $ and $ q \colon Y \to S $ of a closed, orientable surface $ S $ of negative Euler characteristic.

\begin{restatable}{question+}{OurQuestion} \label{que:genericIsospectrality}

If $ X $ and $ Y $ are generically unmarked simple length isospectral over $ S $, are the covers $ p $ and $ q $ isomorphic?

\end{restatable}

Here (and in what follows), we say that a property of hyperbolic metrics on $ X $ and/or $ Y $ holds generically over $ S $ if it holds for the pullbacks of hyperbolic metrics in a dense and full measure subset of the Teichm\"{u}ller space of $ S $. Our strategy for answering \autoref{que:genericIsospectrality} is discussed in prior work of Aougab--Lahn--Loving--Xiao \cite{ALLX20} and at the beginning of \autoref{sec:simpleLengthSpectra}.

\subsection{Main results}
\label{sub:main-results}

Our first result implies that there are infinitely many pairs of hyperbolic surfaces (i.e. topological surfaces equipped with a choice of hyperbolic metric) which are unmarked length isospectral but not unmarked simple length isospectral.

\setcounter{theorem+}{0}

\begin{restatable}{theorem+}{examples} \label{thm:examples}

Every closed surface $ S $ of negative Euler characteristic admits a pair of finite-degree covers which are unmarked length isospectral over $ S $ but generically not unmarked simple length isospectral over $ S $.

\end{restatable}

In fact, the pairs of isospectral hyperbolic surfaces we obtain in the proof of \autoref{thm:examples} have the stronger property that the respective multi-sets of lengths associated to curves with any bounded number of self intersections do not coincide.

\begin{restatable}{corollary+}{k-isospectral} \label{cor:k-isospectral}

There exist hyperbolic surfaces $ X $ and $ Y $ which are unmarked length isospectral but so that for each $ k \geq 0 $, $ X $ and $ Y $ are not unmarked $ \mathcal{ G }_{ k } $-length isospectral.

\end{restatable}

The statement of \autoref{thm:examples} combines work of Maunchang \cite{Mau13} for a surface of genus $ 2 $ with our work for surfaces of genus at least $ 3 $ (see \autoref{sec:examples}). Maungchang uses Sunada's construction to build a pair of finite-degree covers $ X \to S $ and $ Y \to S $ of a closed surface $ S $ of genus $ 2 $ so that a generic hyperbolic metric on $ S $ pulls back to metrics on $ X $ and $ Y $ which have the same length spectrum but not the same simple length spectrum. Our proof of \autoref{thm:examples} generalizes Maungchang's technique and relies crucially on the following strengthening of \cite[Theorem~1.1]{ALLX20}, which uses entirely new proof techniques in order to remove the regularity hypothesis, thus resolving \cite[Conjecture~1.4]{ALLX20}. Given a finite cover $ p \colon X \to S $, a curve $ \alpha $ on $ X $ is an \emph{elevation} of a curve $ \gamma $ on $ S $ along $ p $ if $ p \of{ \alpha } $ is a nonzero iterate of $ \gamma $. 

\begin{restatable}{theorem+}{noneffective} \label{thm:non-effective}

The covers $ p \colon X \to S $ and $ q \colon Y \to S $ are isomorphic if and only if for every curve $ \gamma $ on $ S $, $ \gamma $ has a simple elevation along $ p $ to $ X $ if and only if it has a simple elevation along $ q $ to $ Y $.

\end{restatable}

In \autoref{thm:examples} and \autoref{cor:k-isospectral}, we use Sunada's construction to guarantee that the covers have the same unmarked length spectrum. We choose the covers somewhat carefully to satisfy the desired properties, but in the spirit of \autoref{que:genericIsospectrality}, one might wonder whether \emph{any} pair of non-isomorphic covers arising from Sunada's construction are generically simple length non-isospectral over $ S $.

\subsection{Applications to simple length isospectrality of covers} \label{sub:applications-length}

Our applications of \autoref{thm:non-effective} aim towards answering \autoref{que:UnmarkedLengthSpectralRigidity} for the unmarked simple length spectra of hyperbolic metrics. More precisely, \autoref{thm:non-effective} also allows us to specify a finite-degree cover of $ S $ by the data of which curves have simple elevations. When one can find a curve that lifts simply to one cover but not the other (\autoref{thm:non-effective} guarantees that one can always do this) and which has no \emph{length twins} (distinct curves on $ S $ whose lengths are commensurable in any hyperbolic metric, see \autoref{sec:traceCommensurability}), this lifting data is sufficient to distinguish the simple length spectra of the pullbacks of generic hyperbolic metrics on $ S $ to $ X $ and $ Y $. In particular, the following proposition together with \autoref{thm:non-effective} yields a tool for verifying that such a pair of covers are not simple length isospectral.

\begin{restatable}{proposition+}{lengthSpectra} \label{prop:lengthSpectra}

Suppose that there is a curve $ \gamma $ on $ S $ and a positive integer $ m \geq 1 $ with the following properties:
\begin{enumerate}

\item[(i)]

$ \gamma $ has no length twins on $ S $; and

\item[(ii)]

$ \gamma $ has different numbers of simple elevations of degree $ m $ along each of the covers $ p \colon X \to S $ and $ q \colon Y \to S $.

\end{enumerate}
Then $ X $ and $ Y $ are generically not unmarked simple length isospectral over $ S $.

\end{restatable}

\noindent We use \autoref{prop:lengthSpectra} to provide an affirmative answer to \autoref{que:genericIsospectrality} when $ p $ and $ q $ are regular and \emph{simply generated}, meaning that the fundamental groups of $ X $ and $ Y $ are generated by the set of all elevations of simple closed curves on $ S $ (see \autoref{sec:simplyGenerated} for a formal definition and discussion).

\begin{restatable}{theorem+}{simplyGeneratedCovers} \label{thm:simplyGeneratedCovers} 

If $ p \colon X \to S $ and $ q \colon Y \to S $ are regular and simply generated, then they are either isomorphic or generically not unmarked simple length isospectral over $ S $.

\end{restatable}

It is known by work of Koberda--Santharoubane \cite{KS16} and Malestein--Putman \cite{MP19} that not all regular covers are simply generated, but their examples are necessarily of very high degree. We therefore think of \autoref{thm:simplyGeneratedCovers} as addressing \autoref{que:genericIsospectrality} in the setting of low-degree regular covers.

We also use \autoref{prop:lengthSpectra} to exhibit a method for identifying when pairs of covers produced from Sunada's construction are unmarked simple length non-isospectral. Sunada's construction begins with a surjective homomorphism of a closed surface group onto a finite group, from which one finds finite covers by pulling back a pair of almost conjugate subgroups. We exhibit a general strategy for identifying when these covers arising from Sunada's construction satisfy the hypotheses of \autoref{prop:lengthSpectra}, and are therefore generically not simple length isospectral over the original surface. We use this method to produce three distinct infinite families of such pairs. This provides a vast generalization of the work of Maungchang \cite{Mau13}, producing infinitely many new examples of pairs of isospectral surfaces with distinct simple length spectra, and proving \autoref{thm:examples}. 

Identifying and controlling the presence of pairs of curves which have the same geodesic length in any hyperbolic metric presents a major challenge when attempting to distinguish covers by their simple length spectra. Such curves were shown to exist and studied extensively by Horowitz \cite{Hor72}. Leininger \cite{Lei03} proved that any such pair has the property that the modulus of their traces must coincide under any representation $ \pi_{1} \of{ S } \to \SL_{ 2 } \of{ \C } $. Indeed, the existence of such curves is a phenomenon that arises from trace identities specific to $ \SL_{ 2 } \of{ \C } $. The existence (or non-existence) of words in the free group $ F_{ 2 } $ on two generators for which the traces coincide in any representation to some $ \SL_{ n } \of{ \C } $ with $ n \geq 3 $ remains a significant open question \cite{HT10}.

\subsection{A remark on compact surfaces with boundary}

Although we phrase our results for closed surfaces of negative Euler characteristic, they also apply to compact surfaces with boundary. Indeed, a short argument involving the double of such a surface extends \cite[Theorem~1.4]{Lei03}, the key ingredient.

\subsection{Organization of paper}
\label{sub:organization}

In \autoref{sec:preliminaries}, we provide an overview of the relevant preliminary concepts and results that we will use throughout the paper. In \autoref{sec:non-effective}, we prove \autoref{thm:non-effective} and a corollary relating to covers produced via Sunada's construction which are both length isospectral and simple length isospectral. In \autoref{sec:simpleLengthSpectra}, we discuss and prove \autoref{prop:lengthSpectra} and \autoref{thm:simplyGeneratedCovers}, and then construct the aforementioned three infinite families of examples of covers which are length isospectral but generically simple length non-isospectral, thereby proving \autoref{thm:examples}.

\subsection{Acknowledgements}

We gratefully acknowledge support from NSF grants DMS-1939936 (Aougab), DMS-1906441 (Lahn), DMS-1902729 \& DMS-2231286 (Loving), and DMS-2005438 \& 2300370 (Miller). We thank Aaron Calderon, Siddhi Krishna, Chris Leininger, Ian Runnels, Nick Salter, Lukas Scheiwiller, Robert Tang, and Samuel Taylor for helpful conversations. We would also like to thank Dongryul Kim, Hugo Parlier, and Andy Putman for their feedback on a draft of this paper.

\section{Preliminaries}\label{sec:preliminaries}

\subsection{Curves and their elevations}

By a \emph{(closed) curve} on a surface $ S $, we mean an equivalence class of continuous maps from the circle $ \S^{ 1 } $ to $ S $, where two such maps are considered equivalent if one is freely homotopic to a reparametrization of the other. Note that with this definition, curves are not oriented. A curve $ \gamma $ on a surface $ S $ is:
\begin{enumerate}

\item[(i)]

\emph{essential} if it has no freely nullhomotopic representative, and \emph{inessential} otherwise;

\item[(ii)]

\emph{primitive} if none of its representatives factors through a non-trivial finite-degree covering map of the circle $ \S^{ 1 } $, and \emph{non-primitive} otherwise; and

\item[(iii)]

\emph{simple} if one of its representatives is injective, and \emph{non-simple} otherwise.

\end{enumerate}
We will implicitly assume that all curves under consideration are essential. Under this assumption, all simple curves are primitive.

Now consider a finite-degree cover $ p \colon X \to S $. A curve $ \alpha $ on $ X $ is an \emph{elevation} of a curve $ \gamma $ on $ S $ along $ p $ of degree $ d \geq 1 $ if there are parametrized representatives such that $ p \circ \alpha $ traces $ d $ times along $ \gamma $. More formally, $ \alpha $ is an elevation of $ \gamma $ along $ p $ of degree $ d $ if there are parametrized representatives so that the following square commutes:
\[
\begin{tikzcd}
\S^{ 1 } \ar[ d , " f " swap ] \ar[ r , " \alpha " ] & X \ar[ d , " p " ] \\
\S^{ 1 } \ar[ r , " \gamma " swap ] & S
\end{tikzcd}
\]
where $ f \colon \S^{ 1 } \to \S^{ 1 } $ is a covering map of degree $ d $. Note that the degree of an elevation is bounded above by the degree of the cover $ p \colon X \to S $.

The \emph{deck group} of a finite-degree cover $ p \colon X \to S $ is the group of homeomorphisms of $ X $ which cover the identity homeomorphism of $ S $; that is,
\[
\Deck \of{ p } \coloneqq \set{ \varphi \in \Homeo \of{ X } : p \circ \varphi = p } .
\]
$ \Deck \of{ p } $ acts freely on each fiber of $ p $ and so is a finite group with $ \card{ \Deck \of{ p } } \leq \deg \of{ p } $. When $ \card{ \Deck \of{ p } } = \deg \of{ p } $, then $ \Deck \of{ p } $ also acts transitively on each fiber. In this case, the covering map $ p \colon X \to S $ is the quotient associated to this action and is called \emph{regular} (these covers are also sometimes called \emph{normal} or \emph{Galois}). Note that the action of the deck group on $ X $ induces an action on the set of degree $ d $ elevations of any curve $ \gamma $ on $ S $.

A choice of basepoint $ s \in S $ yields a bijection between the set of curves on $ S $ and the set of equivalence classes of non-identity elements of the fundamental group $ \pi_{ 1 } \of{ S , s } $, where two elements $ \alpha , \beta \in \pi_{ 1 } \of{ S , s } \setminus \set{ \text{id} } $ are considered equivalent if either $ \alpha $ and $ \beta $ are conjugate or $ \alpha $ and $ \beta^{ -1 } $ are conjugate. From this perspective, primitive curves correspond to equivalence classes of elements which cannot be written non-trivially as a power of another element.

Moreover, choices of basepoints $ x \in X $ and $ s \in S $ so that $ p \of{ x } = s $ give an alternate characterization of elevations of degree $ d $: given non-identity elements $ \alpha \in \pi_{ 1 } \of{ X , x } $ and $ \gamma \in \pi_{ 1 } \of{ S , s } $ which cannot be written non-trivially as a power of another element, the primitive curve on $ X $ represented by $ \alpha $ is an elevation of the primitive curve on $ S $ represented by $ \gamma $ along $ p $ of degree $ d $ if and only if $ p_{ * } \of{ \alpha } $ is conjugate in $ \pi_{ 1 } \of{ S , s } $ to $ \gamma^{ d } $ or $ \gamma^{ -d } $, where $ p_{ * } \colon \pi_{ 1 } \of{ X , x } \into \pi_{ 1 } \of{ S , s } $ is the injective homomorphism induced by $ p $.

\subsection{Length and the Teichm\"{u}ller space}\label{sec:introteich}

The \emph{Teichm\"{u}ller\footnote{The authors would like to draw the reader's attention to increasingly relevant historical context: in addition to his current status as the namesake of many fundamental objects in the field, Oswald Teichm\"{u}ller, may his name be erased, was an ardent Nazi. The authors support ongoing efforts to rename these objects in ways that do not glorify the perpetrators of genocide.} space} $ \Teich \of{ S } $ of a closed surface $ S $ of negative Euler characteristic is the space of isotopy classes of hyperbolic metrics on $ S $. More explicitly, two hyperbolic metrics $ \mu $ and $ \nu $ are \emph{isotopic} if there is an isometry from $ \of{ S , \mu } $ to $ \of{ S , \nu } $ which is isotopic to the identity map on $ S $. Teichm\"{u}ller space is equipped with \emph{Fenchel-Nielsen coordinates} \cite{FN02}, which parametrize $ \Teich \of{ S } $ as a complex manifold homeomorphic to $ \R^{ 3 \abs{ \chi \of{ S } } } $. $ \Teich \of{ S } $ is also equipped with a proper geodesic metric called the \emph{Teichm\"{u}ller metric} as well as an isometric and properly discontinuous action of the mapping class group $ \Mod \of{ S } $, whose associated quotient is the moduli space of hyperbolic surfaces homeomorphic to $ S $.

Pulling back hyperbolic metrics along a finite-degree cover contravariantly defines an embedding of the corresponding Teichm\"{u}ller spaces. More explicitly, given a finite-degree cover $ p \colon X \to S $, there is well-defined map $ p^{ * } \colon \Teich \of{ S } \to \Teich \of{ X } $ defined by the formula $ p^{ * } \of{ \class{ \mu } } \coloneqq \class{ p^{ * } \mu } $. When $ \Teich \of{ S } $ and $ \Teich \of{ X } $ are equipped with their respective Teichm\"{u}ller metrics, this map $ p^{ * } $ is a proper isometric embedding.

A choice of basepoint $ s \in S $ yields another perspective on $ \Teich \of{ S } $. Via the association of a metric with any of its monodromy representations, we may also view $ \Teich \of{ S } $ as the space of $ \PGL_{ 2 } \of{ \R } $-conjugacy classes of discrete and faithful representations of the fundamental group $ \pi_{ 1 } \of{ S , s } $ into $ \PSL_{ 2 } \of{ \R } $. From this vantage, pulling back a metric along a cover $ p \colon X \to S $ corresponds to the restriction of the monodromy representation to the appropriate subgroup; that is, given a choice of basepoints $ x \in X $ and $ s \in S $ such that $ p \of{ x } = s $, $ p^{ * } \of{ \class{ \rho } } = \class{ \rho \circ p_{ * } } $ for all discrete and faithful representations $ \rho \colon \pi_{ 1 } \of{ S , s } \to \PSL_{ 2 } \of{ \R } $.

The \emph{(infimal) length} $ \len_{ \mu } \of{ \gamma } $ of a curve $ \gamma $ with respect to a hyperbolic metric $ \mu $ on $ S $ is defined to be the minimal length of any of its parametrized representatives. Equivalently, $ \len_{ \mu } \of{ \gamma } $ is the length of the unique (up to reparametrization) geodesic representative of $ \gamma $. Since $ \gamma $ has the same length in any isotopic hyperbolic metric, $ \len_{ \of{ \oparg } } \of{ \gamma } $ descends to a real-valued function on the Teichm\"{u}ller space $ \Teich \of{ S } $ defined by the formula $ \class{ \mu } \mapsto \len_{ \mu } \of{ \gamma } $. In fact, these length functions are real-analytic, with
\[
\len_{ \class{ \rho } } \of{ \class{ \gamma } } = 2 \cosh^{ -1 } \of{ \frac{ \abs{ \tr \of{ \rho \of{ \gamma } } } }{ 2 } },
\]
for all elements $ \gamma \in \pi_{ 1 } \of{ S , s } $ and discrete and faithful representations $ \rho \colon \pi_{ 1 } \of{ S , s } \into \PSL_{ 2 } \of{ \R } $.

\subsection{The Thurston compactification of Teichm\"{u}ller space}

Both the Teichm\"{u}ller space $ \Teich \of{ S } $ and the curve complex $ \CC \of{ S } $ (which will be introduced in \autoref{sec:curvegraph}) of a surface $ S $ admit compactifications whose boundary points may be described in the language of geodesic laminations. A \emph{geodesic lamination} on a hyperbolic surface $ S $ is a compact subset of $ S $ which is the union of pairwise disjoint simple geodesics (some, all, or none of which may be closed) called \emph{leaves}. We note that geodesic laminations may also be defined independently of a reference hyperbolic metric on $ S $.

A \emph{transverse measure} on a geodesic lamination $ \lambda $ is a positive real-valued function on the set of arcs transverse to $ \lambda $ which is invariant under leaf-preserving isotopy and absolutely continuous with respect to the Lebesgue measure. A \emph{measured lamination} on $ S $ is a pair $ \of{ \lambda , \mu } $, where $ \lambda $ is a geodesic lamination on $ S $ and $ \mu $ is a tranverse measure on $ \lambda $. We will denote by $ \ML \of{ S } $ the space of measured laminations on $ S $ topologized by the weak* topology.

Two measured laminations on $ S $ are called \emph{projectively equivalent} if they have the same underlying geodesic lamination and their transverse measures differ by multiplication by a positive constant. A projective equivalence class of measured laminations is called a \emph{projective measured lamination}. We will denote the space of such projective measured laminations by $ \PML \of{ S } $, which is homeomorphic to a sphere of dimension $ 6 g - 7 $. Thurston \cite[Theorem~5.3]{Thu82} introduced a compactification $ \overline{ \Teich \of{ S } } $ of Teichm\"{u}ller space, now known as the \emph{Thurston compactification}, in which $ \overline{ \Teich \of{ S } } \setminus \Teich \of{ S } = \PML \of{ S } $. For this reason, $ \PML \of{ S } $ is called the \emph{Thurston boundary} of Teichm\"{u}ller space.

Many of the relevant properties of the Teichm\"{u}ller space extend naturally to the Thurston boundary and Thurston compactification. For example, the mapping class group $ \Mod \of{ S } $ acts naturally on $ \PML \of{ S } $ by homeomorphisms, and this action along with the action of $ \Mod \of{ S } $ on $ \Teich \of{ S } $ jointly comprise an action of $ \Mod \of{ S } $ on $ \overline{ \Teich \of{ S } } = \Teich \of{ S } \sqcup \PML \of{ S } $ by homeomorphisms. Moreover, the isometric embedding of Teichm\"{u}ller spaces induced by a finite-degree cover extends to continuous embeddings of their Thurston boundaries and Thurston compactifications.

\begin{theorem*}[{\cite[Theorem~1]{BMN99}}]

If $ p \colon X \to S $ is a finite-degree cover of a closed, orientable surface $ S $ of negative Euler characteristic, then the induced isometric embedding $ p^{ * } \colon \Teich \of{ S } \into \Teich \of{ X } $ extends to a continuous embedding $ p^{ * } \colon \overline{ \Teich \of{ S } } \into \overline{ \Teich \of{ X } } $.

\end{theorem*}

\subsection{Intersection numbers and the curve complex}\label{sec:curvegraph}

Any pair of curves $ \alpha $ and $ \beta $ on $ S $ has a \emph{(geometric) intersection number} $ i \of{ \alpha , \beta } $, which is the minimal number of transverse crossings between $ \alpha $ and $ \beta $ among all parametrized representatives. When $ \alpha = \beta $, this number $ i \of{ \alpha , \alpha } $ is called the \emph{self-intersection number}. A curve $ \alpha $ on $ S $ is simple if and only if $ i \of{ \alpha , \alpha } = 0 $, and a pair of curves $ \alpha $ and $ \beta $ on $ S $ are called \emph{disjoint} if $ i \of{ \alpha , \beta } = 0 $.

A choice of parametrized representatives for a set of curves on $ S $ is said to be in \emph{minimal position} if every representative achieves its self-intersection number and every pair of representatives achieves its geometric intersection number. It is well-known that any set of closed geodesics on a hyperbolic surface are in minimal position. In particular, this implies that every set of curves on $ S $ has a set of parametrized representatives in minimal position.

The \emph{curve complex} $ \CC \of{ S } $ of a closed surface $ S $ of genus $ g \geq 2 $ is a flag simplicial complex introduced by Harvey \cite{Har81}. The simplices in $ \CC \of{ S } $ correspond to \emph{multi-curves} on $ S $, which are collections of pairwise disjoint simple curves on $ S $; that is, the vertices of $ \CC \of{ S } $ are simple curves on $ S $ and a set of vertices span a simplex if they are pairwise disjoint. There is a natural simplicial action of the mapping class group $ \Mod \of{ S } $ on $ \CC \of{ S } $. Moreover, $ \CC \of{ S } $ may be equipped with a metric $ d_{ S } $ so that every simplex is a regular Euclidean simplex of side length $ 1 $. Harvey \cite{Har81} observed that the $ 1 $-skeleton of $ \CC \of{ S } $ is a connected metric graph.

Masur--Minsky \cite[Theorem~1.1]{MM99} proved that, considered with the metric $ d_{ S } $, the curve complex $ \CC \of{ S } $ is a Gromov hyperbolic geodesic space, and Klarreich \cite{Kla18} constructed an explicit parametrization of its Gromov boundary by the space $ \EL \of{ S } $ of ending laminations of $ S $, which we define in \autoref{sec:CCboundary}.

As proven by Rafi--Schleimer \cite{RS09}, finite-degree covers of surfaces give rise to quasi-isometric embeddings of their curve complexes. Specifically, given a finite-degree cover $ p \colon X \to S $, we define $ \widetilde{ p } \colon \CC \of{ S } \to \CC \of{ X } $ so that for each simple curve $ \gamma \in \CC \of{ S } $, $ \widetilde{ p } \of{ \gamma } $ is a choice of elevation of $ \gamma $ along $ p $ to $ X $. At first glance, our definition of $ \widetilde{ p } $ appears to rely significantly on a choice of elevation for every simple curve on $ S $ and thus its properties will not be well behaved when defined as described. However, any two elevations of a simple curve $ \gamma $ along $ p $ are disjoint on $ X $, since any intersection between distinct elevations on $ X $ would descend to a self-intersection of $ \gamma $ on $ S $. Thus different choices of elevations for $ \widetilde{ p } $ yield maps which have uniformly bounded distance at most $ 1 $, and $ \widetilde{ p } \colon \CC \of{ S } \to \CC \of{ X } $ is said to be \emph{coarsely well-defined}.

\begin{theorem*}[{\cite[Theorem~7.1]{RS09}}] \label{thm:CCembedding}

If $ p \colon X \to S $ is a finite-degree cover of a closed, orientable surface $ S $ of negative Euler characteristic, then $ \widetilde{ p } \colon \CC \of{ S } \to \CC \of{ X } $ is a quasi-isometric embedding.

\end{theorem*}

Note that in \cite{RS09}, the above theorem is formulated in terms of the many-to-one relation which associates a simple curve on the base surface to its full pre-image on the cover, which is a multi-curve. Since a multi-curve has diameter at most $ 1 $ in the curve complex, this distinction will not matter for our purposes.

\subsection{The Gromov boundary of the curve complex}\label{sec:CCboundary}

Recall that a geodesic lamination of a closed hyperbolic surface $ S $ is a closed subset of $ S $ foliated by simple geodesics, each of which is called a \emph{leaf}. We will call a geodesic lamination on a surface $ S $:
\begin{enumerate}

\item[(i)]

\emph{minimal} if each leaf is dense in $ S $;

\item[(ii)]

\emph{filling} if every component of the complement of $ \lambda $ in $ S $ is simply connected; and

\item[(iii)]

\emph{ending} if $ \lambda $ is both minimal and filling.

\end{enumerate}
Note that these properties are well-defined irrespective of the choice of hyperbolic metric and that, necessarily, a minimal geodesic lamination cannot contain a closed leaf.

The space $ \EL \of{ S } $ of ending laminations on $ S $ is topologized as follows: Denote by $ \PMEL \of{ S } $ the space of projective measured laminations whose underlying geodesic lamination is both minimal and filling. There is a natural surjective map $ \of{ \oparg }_{ \topol } \colon \PMEL \of{ S } \to \EL \of{ S } $ defined by the formula $ \class{ \of{ \lambda , \mu } }_{ \topol } = \lambda $; that is, this map sends a projective measured ending lamination to its underlying geodesic lamination. We equip $ \EL \of{ S } $ with the quotient topology associated to this map. An equivalent characterization of this topology is described by Hamenst\"adt \cite{Ham06}, in which it is called the \emph{coarse Hausdorff topology}.

As previously mentioned, Klarreich \cite{Kla18} parametrized the Gromov boundary of the curve complex $ \CC \of{ S } $ of a closed surface $ S $ of negative Euler characteristic by the space $ \EL \of{ S } $ of ending laminations.

\begin{theorem*}[{\cite[Theorem~1.3]{Kla18}}]\label{thm:boundaryCC}

There is a natural homeomorphism $ \Lambda_{ S } \colon \CC \of{ S } \sqcup \partial \CC \of{ S } \to \CC \of{ S } \sqcup \EL \of{ S } $ which is the identity on $ \CC \of{ S } $ and restricts to a homeomorphism between $ \partial \CC \of{ S } $ and $ \EL \of{ S } $.

\end{theorem*}

Recall that a finite-degree cover $ p \colon X \to S $ induces a coarsely well-defined quasi-isometric embedding $ \widetilde{ p } \colon \CC \of{ S } \to \CC \of{ X } $. This in turn induces a continuous embedding $\widetilde{ p } \colon \partial\CC \of{ S } \into \partial\CC \of{ X } $, which is both well-defined and canonical by the discussion in the previous subsection. By conjugating by the homeomorphisms provided in \autoref{thm:boundaryCC}, we obtain a continuous embedding $ \EL \of{ S } \into \EL \of{ X } $ so that the following diagram commutes
\[
\begin{tikzcd}
\partial \CC \of{ S } \ar[ r , " \widetilde{ p } " , hook ] \ar[ d , " \Lambda_{ S } " swap ] & \partial \CC \of{ X } \ar[ d , " \Lambda_{ X } " ] \\
\EL \of{ S } \ar[ r , hook ] & \EL \of{ X }
\end{tikzcd}
\]
One can verify that in the presence of a fixed hyperbolic metric on $ S $, this map takes an ending lamination $ \lambda \in \EL \of{ S } $ to its full pre-image
\[
\of{ \Lambda_{ X } \circ \widetilde{ p } \circ \Lambda_{ S }^{ -1 } } \of{ \lambda } = \Lambda_{ X } \of{ \widetilde{ p } \of{ \Lambda_{ S }^{ -1 } \of{ \lambda } } } = p^{ -1 } \of{ \lambda } = \set{ x \in X : p \of{ x } \in \lambda } ,
\]
which is an ending lamination on $ X $ when $ X $ is equipped with the pullback metric.

In the rest of this article, we will implicitly use the above theorem to identify the space of ending laminations with the Gromov boundary of the curve complex, and we will abuse notation by writing $ \widetilde{ p } \colon \EL \of{ S } \into \EL \of{ X } $.

\subsection{The simple length spectrum and Sunada's construction} \label{sec:Sunada}

The \emph{(unmarked) length spectrum} (resp. \emph{simple length spectrum}) of a closed hyperbolic surface $ X $ is the multi-set of lengths of closed geodesics (resp. simple closed geodesics), counted with multiplicity. We will call two such surfaces \emph{length isospectral} (resp. \emph{simple length isospectral}) if they have the same length spectrum (resp. simple length spectrum).

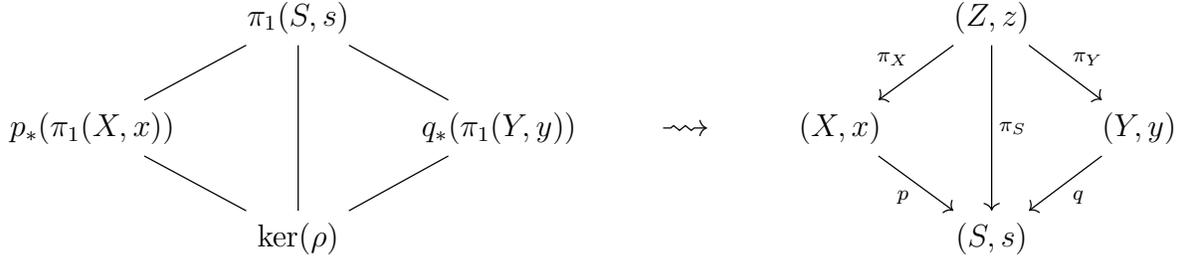
\begin{figure}[ht]
\centering
\begin{tikzcd}[ column sep = 0.6cm ]
& \pi_{ 1 } \of{ S , s } \ar[ dl , no head ] \ar[ dr , no head ] \ar[ dd , no head ] &&&&& \of{ Z , z } \ar[ dl , " \pi_{ X } " swap ] \ar[ dr , " \pi_{ Y } " ] \ar[ dd , " \pi_{ S } " ] \\
p_{ * } \of{ \pi_{ 1 } \of{ X , x } } \ar[ dr , no head ] && q_{ * } \of{ \pi_{ 1 } \of{ Y , y } } \ar[ dl , no head ] & \ar[ r , squiggly ] & \, & \of{ X , x } \ar[ dr , " p " swap ] && \of{ Y , y } \ar[ dl , " q " ] \\
& \ker \of{ \rho } &&&&& \of{ S , s }
\end{tikzcd}
\caption{A Sunada covering diamond from the Galois correspondence} \label{fig:quotientGraphsCoveringDiamonds}
\end{figure}

Most known examples of length isospectral surfaces arise from a construction due to Sunada \cite{Sun85}, which we now describe. Sunada's construction relies on a surjective homomorphism $ \rho \colon \pi_{ 1 } \of{ S , s } \onto G $, called the \emph{constructing representation}, of a closed hyperbolic surface group $ \pi_{ 1 } \of{ S , s } $ onto a finite group $ G $ and a pair of subgroups $ A , B \leq G $, called the \emph{constructing subgroups}. The pre-images $ \rho^{ -1 } \of{ A } $ and $ \rho^{ -1 } \of{ B } $ are finite-index subgroups of $ \pi_{ 1 } \of{ S , s } $ and so the Galois correspondence between subgroups of $ \pi_{ 1 } \of{ S , s } $ and based covers of $ \of{ S , s } $ yields the commutative diagram of finite-degree based covers drawn in \autoref{fig:quotientGraphsCoveringDiamonds}.

We will call a commutative diagram of this form a \emph{(based) Sunada covering diamond}. Note that since $ \ker \of{ \rho } $ is normal in $ \pi_{ 1 } \of{ S , s } $, $ \pi_{ S } \colon Z \to S $ is regular, with deck group
\[
\Deck \of{ \pi_{ S } } \cong \frac{ N_{ \pi_{ 1 } \of{ S , s } } \of{ \ker \of{ \rho } } }{ \ker \of{ \rho } } = \frac{ \pi_{ 1 } \of{ S , s } }{ \ker \of{ \rho } } \cong \im \of{ \rho } = G ,
\]
where $ N_{ \pi_{ 1 } \of{ S , s } } \of{ \ker \of{ \rho } } $ is the normalizer of $ \ker \of{ \rho } $ in $ \pi_{ 1 } \of{ S , s } $. Similarly, $ \pi_{ X } $ and $ \pi_{ Y } $ are regular, with deck groups $ \Deck \of{ \pi_{ X } } \cong A $ and $ \Deck \of{ \pi_{ Y } } \cong B $.

The constructing subgroups $ A , B \leq G $ of the finite group $ G $ are called \emph{almost conjugate} (or \emph{Gassmann equivalent}) if they contain the same number of elements of every conjugacy class in $ G $; that is, $ A $ and $ B $ are almost conjugate in $ G $ if $ \card{ A \cap C } = \card{ B \cap C } $ for every conjugacy class $ C \subseteq G $. Note that two conjugate subgroups of $ G $ are almost conjugate, but almost conjugate subgroups need not be conjugate in $ G $.

\begin{theorem*}[\cite{Sun85}] \label{thm:Sunada}

Let $ p \colon X \to S $ and $ q \colon Y \to S $ be the finite-degree covers of a closed surface $ S $ of negative Euler characteristic associated to a surjective homomorphism $ \rho \colon \pi_{ 1 } \of{ S , s } \onto G $ and a pair of subgroups $ A , B \leq G $ as described above.
\begin{enumerate}

\item[1.]

If $ A $ and $ B $ are almost conjugate in $ G $, then for every $ \mu \in \Teich \of{ S } $, $ p^{ * } \of{ \mu } $ and $ q^{ * } \of{ \mu } $ are length isospectral; and

\item[2.]

If $ A $ and $ B $ are almost conjugate but not conjugate in $ G $, then for a generic $ \mu \in \Teich \of{ S } $, $ p^{ * } \of{ \mu } $ and $ q^{ * } \of{ \mu } $ are not isometric.

\end{enumerate}

\end{theorem*}

The first part of Sunada's theorem can be stated in more general contexts involving actions of finite groups on complete and connected Riemannian manifolds by isometries with at most finitely many fixed points. A proof of this more general form (which implies the first part of \autoref{thm:Sunada}) can be found in \cite[291-294]{Bus10}.

Fix two finite-degree covers $ p \colon X \to S $ and $ q \colon Y \to S $ of a closed, orientable surface $ S $ of negative Euler characteristic. For a given $ \mu \in \Teich \of{ S } $ on $ S $, we will say that $ X $ and $ Y $ are \emph{simple length isospectral over $ S $} if $ \of{ X , p^{ * } \of{ \mu } } $ and $ \of{ Y , q^{ * } \of{ \mu } } $ are simple length isospectral. The manifold structure on $ \Teich \of{ S } $ (or equivalently, the Weil--Petersson measure on $ \Teich \of{ S } $) gives a notion of genericness for these kinds of conditions; we say that $ X $ and $ Y $ generically have a certain property over $ S $ if the set of all $ \mu \in \Teich \of{ S } $ for which $ \of{ X , p^{ * } \of{ \mu } } $ and $ \of{ Y , q^{ * } \of{ \mu } } $ do not have that property has null measure in $ \Teich \of{ S } $. In this language, we can rephrase \autoref{thm:Sunada} as saying that when the constructing subgroups are almost conjugate but not conjugate, the resulting covers are isospectral but generically not isometric.

\subsection{Trace-commensurability in free and surface groups} \label{sec:traceCommensurability}

The \emph{character} of an element $ \alpha $ in a group $ \Gamma $ is the complex-valued function $ \tr \of{ \alpha } $ on the space $ \hom \of{ \Gamma , \SL_{ 2 } \of{ \C } } $ of linear representations of $ \Gamma $ in $ \SL_{ 2 } \of{ \C } $ defined by the formula $ \tr \of{ \alpha } \of{ \rho } = \tr \of{ \rho \of{ \alpha } } $. We will call two elements $ \alpha , \beta \in \Gamma $ \emph{trace-commensurable} if there are non-zero integers $ m , n$ so that $ \tr \of{ \alpha^{ m } }^{ 2 } = \tr \of{ \beta^{ n } }^{ 2 } $ as functions on $ \hom \of{ \Gamma , \SL_{ 2 } \of{ \C } } $. More explicitly, $ \alpha $ and $ \beta $ are trace-commensurable if $ \tr \of{ \rho \of{ \alpha^{ m } } }^{ 2 } = \tr \of{ \rho \of{ \beta^{ n } } }^{ 2 } $ for all linear representations $ \rho \colon \Gamma \to \SL_{ 2 } \of{ \C } $.

This equivalence relation on $ \Gamma $ is a natural generalization of the notion of \emph{trace equivalence} defined in \cite{Lei03}, in which $ \Gamma $ is a closed hyperbolic surface group and $ m = n = 1 $.

\begin{observation*} \label{obs:traceCommensurability}

Consider two elements $ \alpha $ and $ \beta $ of a group $ \Gamma $.

\begin{enumerate}

\item[(i)]

If $ \alpha $ and $ \beta $ are trace-commensurable, then there are non-zero integers $ m $ and $ n $ such that $ \tr \of{ \alpha^{ m } }^{ 2 } = \tr \of{ \beta^{ n } }^{ 2 } $. However, the trace relation
\[
\tr \of{ \vect{ A } } \cdot \tr \of{ \vect{ B } } = \tr \of{ \vect{ A } \vect{ B } } + \tr \of{ \vect{ A } \vect{ B }^{ -1 } }
\]
for matrices $ \vect{ A } , \vect{ B } \in \SL_{ 2 } \of{ \C } $ implies that
\[
\tr \of{ \alpha^{ 2 m } } = \tr \of{ \alpha^{ m } }^{ 2 } - 2 = \tr \of{ \beta^{ n } }^{ 2 } - 2 = \tr \of{ \beta^{ 2 n } } .
\]

\item[(ii)]

Suppose that $ \alpha $ and $ \beta $ are trace-commensurable, so that there are non-zero integers $ m $ and $ n $ such that $ \tr \of{ \alpha^{ m } }^{ 2 } = \tr \of{ \beta^{ n } }^{ 2 } $. Let $ \varphi \colon \Gamma \to \Delta $ be a group homomorphism and note that for any linear representation $ \rho \colon \Delta \to \SL_{ 2 } \of{ \C } $ of $ \Delta $, the composition $ \rho \circ \varphi \colon \Gamma \to \SL_{ 2 } \of{ \C } $ is a linear representation of $ \Gamma $, which implies that
\begin{align*}
\tr \of{ \rho \of{ \varphi \of{ \alpha }^{ m } } }^{ 2 } & = \tr \of{ \rho \of{ \varphi \of{ \alpha^{ m } } } }^{ 2 } = \tr \of{ \of{ \rho \circ \varphi } \of{ \alpha^{ m } } }^{ 2 } = \tr \of{ \of{ \rho \circ \varphi } \of{ \alpha^{ m } }^{ 2 } } \\
& = \tr \of{ \of{ \rho \circ \varphi } \of{ \beta^{ n } }^{ 2 } } = \tr \of{ \of{ \rho \circ \varphi } \of{ \beta^{ n } } }^{ 2 } = \tr \of{ \rho \of{ \varphi \of{ \beta^{ n } } } }^{ 2 } = \tr \of{ \rho \of{ \varphi \of{ \beta }^{ n } } }^{ 2 } .
\end{align*}
Thus $ \varphi \of{ \alpha } $ and $ \varphi \of{ \beta } $ are trace-commensurable as elements of $ \Delta $, with the same commensurability constants.

\end{enumerate}

\end{observation*}

We will call two curves $ \alpha $ and $ \beta $ on $ S $ \emph{(hyperbolic) length-commensurable} if there are positive integers $ m , n \geq 1 $ such that
\[
m \cdot \len_{ \mu } \of{ \alpha } = n \cdot \len_{ \mu } \of{ \beta },
\]
for all points $ \mu \in \Teich \of{ S } $. This equivalence relation on the space of currents is a natural generalization of the notion of \emph{hyperbolic equivalence} of curves as in \cite{Lei03}, in which $ \alpha $ and $ \beta $ are primitive curves and $ m = n = 1 $. Leininger gave an equivalent characterization in terms of certain characters of representations of the fundamental group, which we now present.

Fix a basepoint $ s \in S $ on a closed, orientable surface $ S $ of negative Euler characteristic, and suppose that two elements $ \alpha , \beta \in \pi_{ 1 } \of{ S , s } $ are trace-commensurable, so that there are non-zero integers $ m , n $ such that $ \tr \of{ \rho \of{ \alpha }^{ m } }^{ 2 } = \tr \of{ \rho \of{ \beta }^{ n } }^{ 2 } $ for all linear representations $ \rho \colon \pi_{ 1 } \of{ S , s } \to \SL_{ 2 } \of{ \C } $. Note that any point $ \mu \in \Teich \of{ S } $ is represented by a discrete and faithful representation $ \pi_{ 1 } \of{ S , s } \to \PSL_{ 2 } \of{ \R } $, which lifts to a linear representation $ \pi_{ 1 } \of{ S , s } \to \SL_{ 2 } \of{ \R } \subsetneq \SL_{ 2 } \of{ \C } $. Thus
\begin{align*}
\abs{ m } \cdot \len_{ \mu } \of{ \alpha } & = \len_{ \mu } \of{ \alpha^{ m } } =  2 \cosh^{ -1 } \of{ \frac{ \abs{ \tr \of{ \rho \of{ \alpha^{ m } } } } }{ 2 } } \\
& = 2 \cosh^{ -1 } \of{ \frac{ \abs{ \tr \of{ \rho \of{ \beta^{ n } } } } }{ 2 } } = \len_{ \mu } \of{ \beta^{ n } } = \abs{ n } \cdot \len_{ \mu } \of{ \beta } .
\end{align*}
In summary, the above argument shows that trace-commensurable elements of $ \pi_{ 1 } \of{ S , s } $ represent length-commensurable curves on $ S $. The following theorem of Leininger provides a converse. Note that Leininger considers the slightly more restrictive relationships of trace equivalence and hyperbolic equivalence mentioned above.

\begin{theorem*}[{\cite[Theorem~1.4]{Lei03}}] \label{thm:equivalence}

Fix a basepoint $ s \in S $ on a closed, orientable surface $ S $ of negative Euler characteristic. Two elements $ \alpha , \beta \in \pi_{ 1 } \of{ S , s } $ are trace-commensurable in $ \pi_{ 1 } \of{ S , s } $ if and only if the curves $ \class{ \alpha } $ and $ \class{ \beta } $ they represent are length-commensurable.

Moreover, for any two trace-commensurable elements $ \alpha , \beta \in \pi_{ 1 } \of{ S , s } $ representing primitive curves and any simple curve $ \gamma $ on $ S $, $ i \of{ \class{ \alpha } , \gamma } = 0 $ if and only if $ i \of{ \class{ \beta } , \gamma } = 0 $.

\end{theorem*}

Consider elements $ \alpha $ and $ \beta $ of a group $ \Gamma $. If some non-zero powers $ \alpha^{ m } $ and $ \beta^{ n } $ are conjugate in $ \Gamma $, then the conjugation-invariance of the trace implies that $ \alpha $ and $ \beta $ are trace-commensurable. However, not all pairs of trace-commensurable elements need be of this form. We will call such a pair of elements, which are trace-commensurable but not for the aforementioned reason, \emph{trace twins}.

It may be surprising to learn that examples of trace twins exist in many groups, including free groups and surface groups. For example, consider the free group $ F_{ 2 } = \gen{ a , b \mid } $ on two generators $ a $ and $ b $. Horowitz computes \cite[Example~5.2]{Hor72} that
\begin{align*}
\tr \of{ \rho \of{ a^{ 2 } b^{ -1 } a b } }^{ 2 } & = \tr \of{ \rho \of{ a^{ 2 } b a b^{ -1 } } }^{ 2 } & \tr \of{ \rho \of{ a b^{ 2 } a^{ 2 } b } }^{ 2 } & = \tr \of{ \rho \of{ a b a^{ 2 } b^{ 2 } } }^{ 2 } ,
\end{align*}
under any representation $ \rho \colon F_{ 2 } \to \SL_{ 2 } \of{ \C } $, but no powers of these elements are conjugate to each other. Thus $ a^{ 2 } b^{ -1 } a b $ and $ a^{ 2 } b a b^{ -1 } $ are trace twins in $ F_{ 2 } $ and so are $ a b^{ 2 } a^{ 2 } b $ and $ a b a^{ 2 } b^{ 2 } $.

By \autoref{thm:equivalence} above, trace twins in a surface group yield homotopically distinct curves on the surface whose lengths are commensurable in any hyperbolic metric. As described in greater detail in \autoref{sec:simpleLengthSpectra}, we will call such a pair of curves \emph{length twins}. We will find it very useful to find elements in a group which have no trace twins. Horowitz gives useful characterizations of such elements in free groups and uses them to build a wealth of examples.

\begin{theorem*}[{\cite[Corollary~7.2]{Hor72}}] \label{thm:horowitz}

Let $ a , b \in F $ be primitive elements of a free group $ F $. For any integers $ j , k \in \Z $ and element $ u \in F $, if $ a^{ j } b^{ k } $ and $ u $ are trace equivalent in $ F $, then $ u $ is conjugate in $ F $ to either $ a^{ j } b^{ k } $ or its inverse.

\end{theorem*}

In \autoref{sec:noLengthTwins}, we extend Horowitz's original argument to show that elements of the form $ a^{ j } b^{ k } $ have no trace twins in $ F $ (this involves dealing with commensurability constants other than $ 1 $), and then derive a consequence for the lengths of certain curves on surfaces represented by such elements.

\section{Covers and simple elevations} \label{sec:non-effective}

In this section, we prove \autoref{thm:non-effective} and then deduce \autoref{cor:Sunada}, a corollary related to Sunada's theorem. In \autoref{sec:uniqueErgodicity} and \autoref{sec:deckGroups}, we prove some preliminary technical lemmas concerning the space of projective measured laminations and the deck group of a finite-degree cover, respectively. We use these results in the proof of \autoref{thm:non-effective} in \autoref{sec:mainthm}.

\subsection{Uniquely ergodic ending laminations} \label{sec:uniqueErgodicity}

A geodesic lamination is \emph{uniquely ergodic} if it admits a single projective equivalence class of transverse measures. Masur's criterion provides a sufficient condition for the endpoint of a Teichm\"{u}ller geodesic to be uniquely ergodic based on recurrence to the thick part of Teichm\"{u}ller space, which we now define: given a positive real number $ \epsilon > 0 $, the \emph{$ \epsilon $-thick part} of $ \Teich \of{ X } $ is the subset corresponding to metrics with injectivity radius strictly greater than $ \epsilon $. In other words, a point $ \class{ \mu } \in \Teich \of{ X } $ is in the $ \epsilon $-thick part if $ \of{ X , \mu } $ has no closed geodesic with length at most $ \epsilon $.

\begin{theorem*}[\cite{Mas92}] \label{thm:MasurCriterion}

Let $ R $ be a geodesic ray in $ \Teich \of{ X } $ with endpoint $ \lambda \in \PML \of{ X } $. If there is a positive real number $ \epsilon > 0 $ such that $ R $ recurs infinitely often to the $ \epsilon $-thick part of $ \Teich \of{ X } $, then the underlying lamination $ \lambda_{ \topol } $ is uniquely ergodic.

\end{theorem*}

The following result of Lindenstrauss--Mirzakhani implies that the mapping class group orbit of any uniquely ergodic measured lamination is dense in the space of measured laminations.

\begin{theorem*}[{\cite[Theorem~8.9]{LM08}}] \label{thm:LM08}

A measured lamination $ \lambda \in \ML \of{ X } $ on a closed, orientable surface $ X $ of negative Euler characteristic has dense $ \Mod \of{ X } $-orbit if and only if its support contains no simple closed curves.

\end{theorem*}

We now use the above two results to show that the space $ \PMEL_{ \epsilon } \of{ X } $ of projective measured ending laminations whose underlying laminations are the endpoints of Teichm\"{u}ller geodesic rays in $ \Teich \of{ X } $ which recur infinitely often to the $ \epsilon $-thick part is dense in $ \PML \of{ X } $. We will use this in the proof of \autoref{thm:non-effective}.

\begin{lemma*} \label{lem:densePMEL}

For any closed, orientable surface $ X $ of negative Euler characteristic, there is a positive real number $ \epsilon_{ X } > 0 $ such that $ \PMEL_{ \epsilon } \of{ X } $ is dense in $ \PML \of{ X } $ for any $ 0 < \epsilon \leq \epsilon_{ X } $.

\begin{proof}

For any positive real number $ \epsilon > 0 $, let $ \MEL_{ \epsilon } \of{ X } \coloneqq F^{ -1 } \of{ \PMEL_{ \epsilon } \of{ X } } $, where $ F \colon \ML \of{ X } \onto \PML \of{ X } $ is the quotient map. We will first show that there exists $ \epsilon_{ X } > 0 $ so that $ \MEL_{ \epsilon_X } \of{ X } $ is dense in $ \ML \of{ X } $.

Let $ \phi \in \Mod \of{ X } $ be a pseudo-Anosov mapping class with stable lamination $ \lambda_{ \topol } \in \EL \of{ X } $. There is a positive real number $ \epsilon_X > 0 $, which we fix for the remainder of this proof, so that the axis $ R $ of $ \phi $ in $ \Teich \of{ X } $ is contained in the $ \epsilon_X $-thick part of $ \Teich \of{ X } $. This implies that any measured lamination $ \lambda $ with underlying lamination $ \lambda_{ \topol } $ has the property that $F(\lambda)$ lies in $ \PMEL_{ \epsilon_X } \of{ X } $. Fix such a measured lamination $ \lambda $.

Since $ \Mod \of{ X } $ acts by isometries of the Teichm\"{u}ller metric which preserve the $ \epsilon_{ X } $-thick part of $ \Teich \of{ X } $, $ \MEL_{ \epsilon_{ X } } \of{ X } $ is invariant under this action. Indeed, for any mapping class $ g \in \Mod \of{ X } $, $ F \of{ g \cdot \lambda } = g \cdot F \of{ \lambda } $ is an endpoint of the Teichm\"{u}ller geodesic $ g \cdot R $, which is contained in the $ \epsilon_{ X } $-thick part of $ \Teich \of{ X } $. Thus, $ \MEL_{ \epsilon_{ X } } \of{ X } $ contains the entire orbit $ \Mod \of{ X } \cdot \lambda $.

\autoref{thm:MasurCriterion} implies that the underlying lamination $ \lambda_{ \topol } $ of $ \lambda $ is uniquely ergodic, and so \autoref{thm:LM08} implies that its orbit under the mapping class group $ \Mod \of{ X } $ is dense in $ \ML \of{ X } $. Since $ \MEL_{ \epsilon_X } \of{ X } $ contains this orbit by the above argument, it is also dense in $ \ML \of{ X } $. Therefore, since $ \PMEL_{ \epsilon_{ X } } \of{ X } = F \of{ \MEL_{ \epsilon_{ X } } \of{ X } } $ is the image of a dense subset under a continuous surjection, it is dense in $ \PML \of{ X } $. Finally, note that for any $ 0 < \epsilon \leq \epsilon_{ X } $, $ \PMEL_{ \epsilon_{ X } } \of{ X } \subseteq \PMEL_{ \epsilon } \of{ X } $ and so $ \PMEL_{ \epsilon } \of{ X } $ is also dense in $ \PML \of{ X } $, as desired. \qedhere

\end{proof}

\end{lemma*}

\subsection{The deck group of a finite-degree cover} \label{sec:deckGroups}

The pullback along a finite-degree regular cover of a hyperbolic metric on the base is a hyperbolic metric on the total space for which all deck transformations are isometries. If the base is a closed surface of genus at least three, these deck transformations generically comprise the entire isometry group.

\begin{lemma*} \label{lem:noisom}

For any finite-degree regular cover $ \pi \colon Z \to X $ of a closed, orientable surface $ X $ of genus at least $ 3 $, there is a hyperbolic metric $ \mu $ on $ X $ such that $ \Isom \of{ Z , \pi^{ * } \mu } = \Deck \of{ \pi } $.

\begin{proof}

For any hyperbolic metric $ \mu $ on $ X $, it follows directly from the definitions of the deck group and the pullback metric that $ \Deck \of{ \pi } \leq \Isom \of{ Z , \pi^{ * } \mu } $. It therefore suffices to show there is a hyperbolic metric $ \mu $ on $ X $ so that the reverse containment holds.

Fix basepoints $ x \in X $ and $ z \in Z $ such that $ \pi \of{ z } = x $, and let $ \rho_{ \mu } \colon \pi_{ 1 } \of{ X , x } \into \PSL_{ 2 } \of{ \R } $ be an associated holonomy representation. Such a choice is well-defined up to conjugation by an element of $ \PGL_{ 2 } \of{ \R } $. We will denote
\begin{align*}
\Gamma & \coloneqq \rho_{ \mu } \of{ \pi_{ 1 } \of{ X , x } } , & \Delta & \coloneqq \rho_{ \pi^{ * } \mu } \of{ \pi_{ 1 } \of{ Z , z } } ,
\end{align*}
where $ \rho_{ \pi^{ * } \mu } \coloneqq \rho_{ \mu } \circ \pi_{ * } \colon \pi_{ 1 } \of{ Z , z } \into \PSL_{ 2 } \of{ \R } $ is a holonomy representation of the cover $ \of{ Z , \pi^{ * } \mu } $. Recall that the \emph{commensurator} of $ \Gamma $ in $ \PSL_{ 2 } \of{ \R } $ is defined by
\[
\Comm \of{ \Gamma } \coloneqq \set{ A \in \PSL_{ 2 } \of{ \R } : \Gamma \cap A \Gamma A^{ -1 } \textrm{ has finite index in both } \Gamma \textrm{ and } A \Gamma A^{ -1 } } ,
\]
and that when $ \Gamma $ is non-arithmetic, $ \Comm \of{ \Gamma } $ is a lattice, and the quotient $ \H^{ 2 } / \Comm \of{ \Gamma } $ is the unique minimal area representative in the commensurability class of $ \Gamma $. Moreover, since $ \Gamma \leq \Comm \of{ \Gamma } $ and $ \Comm \of{ \Delta } = \Comm \of{ \Gamma } $, there are covering maps $ \H^{ 2 } / \Gamma \to \H^{ 2 } / \Comm \of{ \Gamma } $ and $ \H^{ 2 } / \Delta \to \H^{ 2 } / \Comm \of{ \Gamma } $. 

For the rest of the proof, we fix a metric $ \mu $ on $ X $ for which $ \Comm \of{ \Gamma } = \Gamma $. By \cite[Theorems~2,~3A]{Gre63}, the set of points in Teichm\"{u}ller space represented by such metrics is a dense open subset of $ \Teich \of{ X } $ when $ X $ has genus at least $ 3 $. Note that with this choice of $ \mu $, $ \Gamma $ is necessarily non-arithmetic by the Margulis dichotomy \cite{Mar75}, as otherwise $ \Comm \of{ \Gamma } $ would be dense in $ \PSL_{ 2 } \of{ \R } $.

Recall that $ \Isom \of{ Z , \pi^{ * } \mu } $ is canonically isomorphic to $ N_{ \PSL_{ 2 } \of{ \R } } \of{ \Delta } / \Delta $, where
\[
N_{ \PSL_{ 2 } \of{ \R } } \of{ \Delta } \coloneqq \set{ A \in \PSL_{ 2 } \of{ \R } : A \Delta A^{ -1 } = \Delta },
\]
is the normalizer of $ \Delta $ in $ \PSL_{ 2 } \of{ \R } $. Note that
\[
\Delta \leq N_{ \PSL_{ 2 } \of{ \R } } \of{ \Delta } \leq \Comm \of{ \Delta } = \Comm \of{ \Gamma } = \Gamma .
\]
Since $ \pi $ is a regular cover, $ \Delta $ is normal in each of these groups, and so
\[
\Isom \of{ Z , \pi^{ * } \mu } \cong \frac{ N_{ \PSL_{ 2 } \of{ \R } } \of{ \Delta } }{ \Delta } \leq \frac{ \Comm \of{ \Gamma } }{ \Delta } = \frac{ \Gamma }{ \Delta } \cong \Deck \of{ \pi } .
\]
By the remark at the beginning of the proof, we see that for any such choice of $ \mu $ as above, $ \Isom \of{ Z , \pi^{ * } \mu } = \Deck \of{ \pi } $ as required. \qedhere

\end{proof}

\end{lemma*}

We remark that the key reason for the hypothesis in \autoref{lem:noisom} that $ X $ has genus at least $ 3 $ is due to the hyper-elliptic involution of surfaces of genus $ 2 $. Specifically, due to the presence of this isometry, we cannot conclude that $ \Comm \of{ \Gamma } = \Gamma $ for a generic metric; in this case, we can only show that $ \ind{ \Comm \of{ \Gamma } }{ \Gamma } \leq 2 $ and therefore that $ \Isom \of{ Z , \pi^{ * } \of{ \mu } } $ is at worst a $ \Z / 2 \Z $ extension of $ \Deck \of{ \pi } $ for generic $ \mu \in \Teich \of{ X } $.

\autoref{lem:noisom} has the following consequences for the collection of mapping classes corresponding to deck transformations of a finite-degree cover.

\begin{proposition*} \label{prop:deckGroup}

Let $ \pi \colon Z \to X $ be a finite-degree cover of a closed, orientable surface $ X $ of negative Euler characteristic and consider the quotient map $ Q \colon \Homeo^{ + } \of{ Z } \onto \Mod \of{ Z } $. The restriction of $ Q $ to the deck group $ \Deck \of{ \pi } $ is injective. Moreover, if $ \pi $ is regular and $ X $ has genus at least $ 3 $, then $ Q $ has image
\[
Q \of{ \Deck \of{ \pi } } = \set{ \varphi \in \Mod \of{ Z } : \varphi \cdot \pi^{ * } \of{ \mu } = \pi^{ * } \of{ \mu } \textrm{ for all } \mu \in \Teich \of{ X } } .
\]

\begin{proof}

Let $ \varphi \in \Deck \of{ \pi } $, fix an orientation and a hyperbolic metric on $ X $, which pullback along $ \pi $ to an orientation and a hyperbolic metric on $ Z $. Then any deck transformation of $ \pi $ is an orientation-preserving isometry of $ X $.

Fix a non-separating simple closed geodesic $ \gamma $ on $ X $ and a component $ \alpha $ of its full pre-image on $ Z $, which is also a non-separating simple closed geodesic. Since $ \alpha $ is simple and non-separating, there is a simple closed geodesic $ \beta $ on $ Z $ which intersects it exactly once. Since $ \varphi $ is an isometry, $ \varphi \of{ \alpha } $ and $ \varphi \of{ \beta } $ are simple closed geodesics on $ Z $ which intersect exactly once.

Suppose that $ Q \of{ \varphi } $ is trivial in $ \Mod \of{ Z } $. Then $ Q \of{ \varphi } $ acts trivially on the homotopy classes of $ \alpha $ and $ \beta $, and so $ \varphi \of{ \alpha } $ and $ \varphi \of{ \beta } $ are homotopic to $ \alpha $ and $ \beta $, respectively. Since the metric on $ Z $ is negatively curved, closed geodesics are unique within their homotopy class, and hence $ \varphi $ preserves $ \alpha $ and $ \beta $. This implies that $ \varphi $ fixes the unique intersection point between $ \alpha $ and $ \beta $. Since any non-trivial element of the deck group $ \Deck \of{ \pi } $ acts freely on the total space $ Z $, this implies that $ \varphi = \id_{ Z } $ is trivial in $ \Homeo^{ + } \of{ Z } $.

Now suppose that $ \pi $ is regular and that $ X $ has genus at least $ 3 $. The containment
\[
Q \of{ \Deck \of{ \pi } } \subseteq \set{ \varphi \in \Mod \of{ Z } : \varphi \cdot \pi^{ * } \of{ \mu } = \pi^{ * } \of{ \mu } \textrm{ for all } \mu \in \Teich \of{ X } },
\]
is clear, since pullbacks of metrics along $ \pi $ are invariant under the deck group $ \Deck \of{ \pi } $. To show the other containment, fix a hyperbolic metric $ \mu $ on $ X $ as in \autoref{lem:noisom}; that is, we require $ \Isom \of{ Z , \pi^{ * } \mu } = \Deck \of{ \pi } $. Let $ \varphi \in \Homeo^{ + } \of{ Z } $, and suppose that $ Q \of{ \varphi } \cdot \pi^{ * } \of{ \class{ \mu } } = \pi^{ * } \of{ \class{ \mu } } $. Then 

\begin{align*}
\class{ \of{ \pi \circ \varphi }^{ * } \mu } & = \class{ \varphi^{ * } \of{ \pi^{ * } \mu } } = Q \of{ \varphi^{ -1 } } \cdot \class{ \pi^{ * } \mu } = Q \of{ \varphi^{ -1 } } \cdot \of{ Q \of{ \varphi } \cdot \class{ \pi^{ * } \mu } } \\
& = Q \of{ \varphi^{ -1 } \circ \varphi } \cdot \class{ \pi^{ * } \mu } = Q \of{ \id_{ Z } } \cdot \class{ \pi^{ * } \mu } = \class{ \pi^{ * } \mu } ,
\end{align*}
and so there is a homeomorphism $ \psi \in \Homeo^{ + }_{ 0 } \of{ Z } $ such that
\[
\of{ \varphi \circ \psi }^{ * } \of{ \pi^{ * } \mu } = \of{ \pi \circ \varphi \circ \psi }^{ * } \mu = \pi^{ * } \mu .
\]
Thus $ \varphi \circ \psi \in \Isom \of{ Z , \pi^{ * } \mu } = \Deck \of{ \pi } $ is an isometry of $ \pi^{ * } \mu $ and therefore a deck transformation of $ \pi $. In particular, this implies that $ Q \of{ \varphi } = Q \of{ \varphi \circ \psi } \in Q \of{ \Deck \of{ \pi } } $. \qedhere

\end{proof}

\end{proposition*}

\subsection{The proof of \autoref{thm:non-effective} and a corollary} \label{sec:mainthm}

We are now ready to prove \autoref{thm:non-effective}, which we restate below for the reader's convenience. Recall that $ p \colon X \to S $ and $ q \colon Y \to S $ are finite-degree covers of a closed surface $ S $ of negative Euler characteristic.

\noneffective*

Since the proof is fairly involved, we begin with a brief outline of our strategy. Given covers $ p \colon X \to S $ and $ q \colon Y \to S $ which satisfy the simple elevation criteria, we first fix a regular cover $ \pi_{ S } \colon Z \to S $ which factors through both $ X $ and $ Y $ (see \autoref{fig:quotientGraphsCoveringDiamonds}). We then show in \autoref{clm:CC} that the simple elevation criteria implies that the orbits of the curve complexes of $ X $ and $ Y $, considered as subsets of the curve complex of $ Z $, under the action of $ \Deck \of{ \pi_{ S } } $ coincide. \autoref{clm:EL}, \autoref{clm:PML}, and \autoref{clm:Teich} show that similar statements hold for $ \EL \of{ X } $ and $ \EL \of{ Y } $; for $ \PML \of{ X } $ and $ \PML \of{ Y } $; and then finally for $ \Teich \of{ X } $ and $ \Teich \of{ Y } $, respectively.

As $ \Teich \of{ X } $ and $ \Teich \of{ Y } $ are properly embedded submanifolds of $ \Teich \of{ Z } $, the dimensionality argument in \autoref{clm:singledeck} shows that if the orbits of the isometrically embedded Teichm\"{u}ller spaces coincide, these Teichm\"{u}ller spaces must be related by a single deck transformation. From this, we conclude in \autoref{clm:conj} that the deck groups of $ p \colon X \to S $ and $ q \colon Y \to S $ are conjugate in $ \Deck \of{ \pi_{ S } } $ (by this very deck transformation) and therefore that the corresponding covers are equivalent.

\begin{proof}[Proof of \autoref{thm:non-effective}]

First, suppose that for all curves $ \gamma $ on $ S $, $ \gamma $ has a simple elevation along $ p $ to $ X $ if and only if it has a simple elevation along $ q $ to $ Y $. Fix basepoints $ x \in X $, $ y \in Y $, and $ s \in S $ such that $ p \of{ x } = q \of{ y } = s $, and let $ K $ be a finite-index normal subgroup of $ \pi_{ 1 } \of{ S , s } $ contained in both $ p_{ * } \of{ \pi_{ 1 } \of{ X , x } } $ and $ q_{ * } \of{ \pi_{ 1 } \of{ Y , y } } $. For example, we may take $ K $ to be the normal core of $ p_{ * } \of{ \pi_{ 1 } \of{ X , x } } \cap q_{ * } \of{ \pi_{ 1 } \of{ Y , y } } $ in $ \pi_{ 1 } \of{ S , s } $; that is,
\[
K = \bigcap_{ \mathclap{ \gamma \in \pi_{ 1 } \of{ S , s } } }{ \gamma \of!\big!{ p_{ * } \of{ \pi_{ 1 } \of{ X , x } } \cap q_{ * } \of{ \pi_{ 1 } \of{ Y , y } } } \gamma^{ -1 } } .
\]
The Galois correspondence between finite-index subgroups of $ \pi_{ 1 } \of{ S , s } $ and finite-degree based covers then yields the based Sunada covering diamond drawn in \autoref{fig:quotientGraphsCoveringDiamonds}.


As discussed in \autoref{sec:Sunada}, $ \pi_{ S } $, $ \pi_{ X } $, and $ \pi_{ Y } $ are regular covers with deck groups $ \Deck \of{ \pi_{ S } } \cong G $, $ \Deck \of{ \pi_{ X } } \cong A $, and $ \Deck \of{ \pi_{ Y } } \cong B $, respectively. By \autoref{prop:deckGroup} above, we may identify the deck groups of these covers by their isomorphic images in the mapping class group $ \Mod \of{ Z } $.

For the remainder of the proof, it may be useful for the reader to refer back to \autoref{sec:preliminaries} for detailed definitions and relevant notation.

\begin{claim} \label{clm:CC}

$ \Deck \of{ \pi_{ S } } \cdot \widetilde{ \pi_{ X } } \of{ \CC \of{ X } } = \Deck \of{ \pi_{ S } } \cdot \widetilde{ \pi_{ Y } } \of{ \CC \of{ Y } } $ as subsets of $ \CC \of{ Z } $.

\begin{proof}[Proof of \autoref{clm:CC}]

Let $ \alpha \in \CC \of{ X } $ be a simple curve on $ X $. Then $ \alpha $ is a simple elevation of $ p \of{ \alpha } $ along $ p $ to $ X $ and so by assumption, $ p \of{ \alpha } $ has a simple elevation $ \beta \in \CC \of{ Y } $ along $ q $ to $ Y $. Note that
\[
\pi_{ S } \of{ g \cdot \widetilde{ \pi_{ X } } \of{ \alpha } } = \pi_{ S } \of{ \widetilde{ \pi_{ X } } \of{ \alpha } } = p \of{ \pi_{ X } \of{ \widetilde{ \pi_{ X } } \of{ \alpha } } } = p \of{ \alpha } = q \of{ \beta } = q \of{ \pi_{ Y } \of{ \widetilde{ \pi_{ Y } } \of{ \beta } } } = \pi_{ S } \of{ \widetilde{ \pi_{ Y } } \of{ \beta } },
\]
for any $ g \in \Deck \of{ \pi_{ S } } $. Since $ \pi_{ S } \colon Z \to S $ is a regular cover, $ g \cdot \widetilde{ \pi_{ X } } \of{ \alpha } = h \cdot \widetilde{ \pi_{ Y } } \of{ \beta } $ for some deck transformation $ h \in \Deck \of{ \pi_{ S } } $, and thus $ \Deck \of{ \pi_{ S } } \cdot \widetilde{ \pi_{ X } } \of{ \CC \of{ X } } \subseteq \Deck \of{ \pi_{ S } } \cdot \widetilde{ \pi_{ Y } } \of{ \CC \of{ Y } } $. By symmetry, the other containment also holds. \qedhere

\end{proof}

\end{claim}

\begin{claim} \label{clm:EL}

$ \Deck \of{ \pi_{ S } } \cdot \widetilde{ \pi_{ X } } \of{ \EL \of{ X } } = \Deck \of{ \pi_{ S } } \cdot \widetilde{ \pi_{ Y } } \of{ \EL \of{ Y } } $ as subsets of $ \EL \of{ Z } $.

\begin{proof}[Proof of \autoref{clm:EL}]

Let $ g \in \Deck \of{ \pi_{ S } } $ and $ \alpha \in \EL \of{ X } $. Since $ \EL \of{ X } $ is the Gromov boundary of $ \CC \of{ X } $ (as discussed in \autoref{sec:CCboundary}), some sequence $ \of{ \alpha_{ n } }_{ n = 0 }^{ \infty } $ in $ \CC \of{ X } $ converges to $ \alpha $. Thus $ \of{ g \cdot \widetilde{ \pi_{ X } } \of{ \alpha_{ n } } }_{ n = 0 }^{ \infty } $ converges to $ g \cdot \widetilde{ \pi_{ X } } \of{ \alpha } $ in $ \CC \of{ Z } $.

By \autoref{clm:CC}, for each $ n  $ there is a deck transformation $ h_{ n } \in \Deck \of{ \pi_{ S } } $ and a simple curve $ \beta_{ n } \in \CC \of{ Y } $ so that $ g \cdot \widetilde{ \pi_{ X } } \of{ \alpha_{ n } } = h_{ n } \cdot \widetilde{ \pi_{ Y } } \of{ \beta_{ n } } $. Since $ \Deck \of{ \pi_{ S } } $ is finite, some element $ h $ appears in the sequence $ \of{ h_{ n } }_{ n = 0 }^{ \infty } $ infinitely often and so, by passing to a subsequence, we may assume without loss of generality that $ h_{ n } = h $ is a fixed deck transformation for all indices $ n \in \N $.

Note that the sequence $ \of{ \beta_{ n } }_{ n = 0 }^{ \infty } $ is unbounded in $ \CC \of{ Y } $. Indeed, since $ \widetilde{ \pi_{ X } } \colon \CC \of{ X } \to \CC \of{ Z } $ is a quasi-isometric embedding, and $ g \colon \CC \of{ Z } \to \CC \of{ Z } $ is an isometry, it follows that the sequence $ \of{ g \cdot \widetilde{ \pi_{ X } } \of{ \alpha_{ n } } }_{ n = 0 }^{ \infty } = \of{ h \cdot \widetilde{ \pi_{ Y } } \of{ \beta_{ n } } }_{ n = 0 }^{ \infty } $ is unbounded in $ \CC \of{ Z } $. Similarly, since $ \widetilde{ \pi_{ Y } } \colon \CC \of{ Y } \to \CC \of{ Z } $ is a quasi-isometric embedding and $ h \colon \CC \of{ Z } \to \CC \of{ Z } $ is an isometry, it follows that $ \of{ \beta_{ n } }_{ n = 0 }^{ \infty } $ is unbounded in $ \CC \of{ Y } $.

Thus by passing to yet another subsequence, we may assume that $ \of{ \beta_{ n } }_{ n = 0 }^{ \infty } $ converges to an ending lamination $ \beta \in \EL \of{ Y } $ on $ Y $. In particular, by continuity of $ \widetilde{ \pi_{ X } } \colon \EL \of{ X } \to \EL \of{ Z } $ and $ \widetilde{ \pi_{ Y } } \colon \EL \of{ Y } \to \EL \of{ Z } $,
\begin{align*}
g \cdot \widetilde{ \pi_{ X } } \of{ \alpha } & = g \cdot \widetilde{ \pi_{ X } } \of{ \lim_{ n \to \infty }{ \alpha_{ n } } } = \lim_{ n \to \infty }{ g \cdot \widetilde{ \pi_{ X } } \of{ \alpha_{ n } } } \\
& = \lim_{ n \to \infty }{ h \cdot \widetilde{ \pi_{ Y } } \of{ \beta_{ n } } } = h \cdot \widetilde{ \pi_{ Y } } \of{ \lim_{ n \to \infty }{ \beta_{ n } } } = h \cdot \widetilde{ \pi_{ Y } } \of{ \beta } ,
\end{align*}
and so $ \Deck \of{ \pi_{ S } } \cdot \widetilde{ \pi_{ X } } \of{ \EL \of{ X } } \subseteq \Deck \of{ \pi_{ S } } \cdot \widetilde{ \pi_{ Y } } \of{ \EL \of{ Y } } $. By symmetry, the other containment also holds. \qedhere

\end{proof}

\end{claim}

\begin{claim} \label{clm:PML}

$ \Deck \of{ \pi_{ S } } \cdot \pi_{ X }^{ * } \of{ \PML \of{ X } } = \Deck \of{ \pi_{ S } } \cdot \pi_{ Y }^{ * } \of{ \PML \of{ Y } } $ as subsets of $ \PML \of{ Z } $.

\begin{proof}[Proof of \autoref{clm:PML}]

Fix $ \epsilon \coloneqq \min \of{ \epsilon_{ X } , \epsilon_{ Y } } > 0 $ as in \autoref{lem:densePMEL}, let $ g \in \Deck \of{ \pi_{ S } } $ and consider a projective measured ending lamination $ \alpha \in \PMEL_{ \epsilon } \of{ X } $ on $ X $ whose underlying geodesic lamination $ \alpha_{ \topol } \in \EL \of{ X } $ is the endpoint of a geodesic ray $ R $ in $ \Teich \of{ X } $ which recurs infinitely often to the $ \epsilon $-thick part. Since lower bounds on injectivity radius are preserved by local isometries, $ \pi_{ X }^{ * } \of{ R } $ is a geodesic ray in $ \Teich \of{ Z } $ which recurs infinitely often to the $ \epsilon $-thick part, with endpoint
\[
\of{ \pi_{ X }^{ * } \of{ \alpha } }_{ \topol } = \pi_{ X }^{ -1 } \of{ \alpha_{ \topol } } = \widetilde{ \pi_{ X } } \of{ \alpha_{ \topol } } \in \EL \of{ Z } .
\]
In particular, Masur's criterion \autoref{thm:MasurCriterion} implies that this ending lamination on $ Z $ is uniquely ergodic. By \autoref{clm:EL}, $ g \cdot \widetilde{ \pi_{ X } } \of{ \alpha_{ \topol } } = h \cdot \widetilde{ \pi_{ Y } } \of{ \beta_{ \topol } } $ for some deck transformation $ h \in \Deck \of{ \pi_{ S } } $ and some ending lamination $ \beta_{ \topol } \in \EL \of{ Y } $ on $ Y $. Since $ \widetilde{ \pi_{ X } } \of{ \alpha_{ \topol } } $ is uniquely ergodic, so too is its image $ \widetilde{ \pi_{ Y } } \of{ \beta_{ \topol } } = h^{ -1 } g \cdot \widetilde{ \pi_{ X } } \of{ \alpha_{ \topol } } $ under the deck transformation $ h^{ -1 } g $.

As $ \widetilde{ \pi_{ Y } } \of{ \beta_{ \topol } } $ comes from pullback of an ending lamination on $ Y $, it is invariant under the action of $ \Deck \of{ \pi_{ Y } } $, and so any representative of the unique projective class of transverse measures on $ \widetilde{ \pi_{ Y } } \of{ \beta_{ \topol } } $ descends to a transverse measure on $ \beta_{ \topol } $. Let $ \beta \in \PML \of{ Y } $ be the corresponding projective measured lamination, and note that $ h \cdot \pi_{ Y }^{ * } \of{ \beta } $ is a projective measured lamination on $ Z $ with underlying geodesic lamination
\[
\of{ h \cdot \pi_{ Y }^{ * } \beta }_{ \topol } = h \cdot \of{ \pi_{ Y }^{ * } \of{ \beta } }_{ \topol } = h \cdot \widetilde{ \pi_{ Y } } \of{ \beta_{ \topol } } = g \cdot \widetilde{ \pi_{ X } } \of{ \alpha_{ \topol } } = g \cdot \of{ \pi_{ X }^{ * } \of{ \alpha } }_{ \topol } = \of{ g \cdot \pi_{ X }^{ * } \of{ \alpha } }_{ \topol } .
\]

\noindent Since $ g \cdot \pi_{ X }^{ * } \of{ \alpha } $ is uniquely ergodic, this implies that $ g \cdot \pi_{ X }^{ * } \of{ \alpha } = h \cdot \pi_{ Y }^{ * } \of{ \beta } $, and so
\[
\Deck \of{ \pi_{ S } } \cdot \pi_{ X }^{ * } \of{ \PMEL_\epsilon \of{ X } } \subseteq \Deck \of{ \pi_{ S } } \cdot \pi_{ Y }^{ * } \of{ \PML \of{ Y } } .
\]
Note that for each deck transformation $ h \in \Deck \of{ \pi_{ S } } $, $ h \cdot \pi_{ Y }^{ * } \colon \PML \of{ Y } \to \PML \of{ Z } $ is a topological embedding, and so $ h \cdot \pi_{ Y }^{ * } \of{ \PML \of{ Y } } $ is compact and therefore closed in the Hausdorff space $ \PML \of{ Z } $. Thus
\[
\Deck \of{ \pi_{ S } } \cdot \pi_{ Y }^{ * } \of{ \PML \of{ Y } } = \bigcup_{ \mathclap{ h \in \Deck \of{ \pi_{ S } } } }{ h \cdot \pi_{ Y }^{ * } \of{ \PML \of{ Y } } },
\]
is closed in $ \PML \of{ Z } $, since $ \Deck \of{ \pi_{ S } } $ is finite. Finally, since $ \PMEL_{ \epsilon } \of{ X } $ is dense in $ \PML \of{ X } $ by \autoref{lem:densePMEL}, $ \pi_{ X }^{ * } $ is a topological embedding, and $ \Deck \of{ \pi_{ S } } $ is finite
\begin{align*}
\Deck \of{ \pi_{ S } } \cdot \pi_{ X }^{ * } \of{ \PML \of{ X } } & = \Deck \of{ \pi_{ S } } \cdot \pi_{ X }^{ * } \of{ \overline{ \PMEL_{ \epsilon } \of{ X } } } \\
& \subseteq \Deck \of{ \pi_{ S } } \cdot \overline{ \pi_{ X }^{ * } \of{ \PMEL_{ \epsilon } \of{ X } } }  = \overline{ \Deck \of{ \pi_{ S } } \cdot \pi_{ X }^{ * } \of{ \PMEL_{ \epsilon } \of{ X } } } \\
& \subseteq \overline{ \Deck \of{ \pi_{ S } } \cdot \pi_{ Y }^{ * } \of{ \PML \of{ Y } } } = \Deck \of{ \pi_{ S } } \cdot \pi_{ Y }^{ * } \of{ \PML \of{ Y } } .
\end{align*}
By symmetry, the other containment also holds. \qedhere

\end{proof}

\end{claim}

\begin{claim} \label{clm:Teich}

$ \Deck \of{ \pi_{ S } } \cdot \pi_{ X }^{ * } \of{ \Teich \of{ X } } = \Deck \of{ \pi_{ S } } \cdot \pi_{ Y }^{ * } \of{ \Teich \of{ Y } } $ as a subset of $\Teich(Z)$.

\begin{proof}[Proof of \autoref{clm:Teich}]

Let $ g \in \Deck \of{ \pi_{ S } } $ and $ \mu \in \Teich \of{ X } $. Fix a projective measured lamination $ \kappa \in \PML \of{ S } $, and let $ \alpha^{ + } \coloneqq p^{ * } \of{ \kappa } \in \PML \of{ X } $ and $ \beta^{ + } \coloneqq q^{ * } \of{ \kappa } \in \PML \of{ Y } $. By \cite{HM79}, $ \alpha^{ + } $ is one endpoint of a geodesic $ L_{ X } \subseteq \Teich \of{ X } $ through $ \mu $. Let $ \alpha^{ - } $ be the other endpoint of $ L_{ X } $. Since $ \pi_{ X }^{ * } \colon \Teich \of{ X } \into \Teich \of{ Z } $ is an isometric embedding and $ g \colon \Teich \of{ Z } \to \Teich \of{ Z } $ is an isometry, $ g \cdot \pi_{ X }^{ * } \of{ L_{ X } } $ is a geodesic in $ \Teich \of{ Z } $ containing $ g \cdot \pi_{ X }^{ * } \of{ \mu } $ with endpoints $ g \cdot \pi_{ X }^{ * } \of{ \alpha^{ + } } $ and $ g \cdot \pi_{ X }^{ * } \of{ \alpha^{ - } } $. 

By \autoref{clm:PML}, $ g \cdot \pi_{ X }^{ * } \of{ \alpha^{ - } } = h \cdot \pi_{ Y }^{ * } \of{ \beta^{ - } } $ for some deck transformation $ h \in \Deck \of{ \pi_{ S } } $ and projective measured lamination $ \beta^{ - } \in \PML \of{ Y } $. On the other hand, $ \pi_{ X }^{ * } \of{ \alpha^{ + } } = \pi_{ S }^{ * } \of{ \kappa } $ is $ \Deck \of{ \pi_{ S } } $-invariant, and so 
\[
g \cdot \pi_{ X }^* \of{ \alpha^{ + } } = g \cdot \pi_{ S }^{ * } \of{ \kappa } = \pi_{ S }^{ * } \of{ \kappa } = h \cdot \pi_{ S }^{ * } \of{ \kappa } = h \cdot \pi_{ Y }^{ * } \of{ \beta^{ + } } .
\]
In summary, $ g \cdot \pi_{ X }^{ * } \of{ L_{ X } } $ is a geodesic in $ \Teich \of{ X } $ containing $ g \cdot \pi_{ X }^{ * } \of{ \mu } $ with endpoints $ g \cdot \pi_{ X }^{ * } \of{ \alpha^{ + } } = h \cdot \pi_{ Y }^{ * } \of{ \beta^{ + } } $ and $ g \cdot \pi_{ Y }^{ * } \of{ \alpha^{ - } } = h \cdot \pi_{ Y }^{ * } \of{ \beta^{ - } } $.

We now show that $ \beta^{ + } $ and $ \beta^{ - } $ jointly fill $ Y $ and therefore are endpoints of a geodesic $ L_{ Y } \subseteq \Teich \of{ Y } $. Indeed, if $ \lambda \in \PML \of{ Y } $ is disjoint from $ \beta^{ + } $ and $ \beta^{ - } $, then $ \pi_{ Y }^{ * } \of{ \lambda } $ is a projective measured lamination on $ Z $ which is disjoint from $ \pi_{ Y }^{ * } \of{ \beta^{ + } } $ and $ \pi_{ Y }^{ * } \of{ \beta^{ - } } $. However, the latter projective measured laminations are the endpoints of the geodesic $ h^{ -1 } g \cdot \pi_{ X }^{ * } \of{ L_{ X } } $ in $ \Teich \of{ Z } $, hence they jointly fill $ Z $. Therefore, no such $ \lambda $ exists, and so $ \beta^{ + } $ and $ \beta^{ - } $ jointly fill $ Y $ as required.

In particular, $ g \cdot \pi_{ X }^{ * } \of{ L_{ X } } $ and $ h \cdot \pi_{ Y }^* \of{ L_{ Y } } $ are geodesics in $ \Teich \of{ Z } $ with the same endpoints, and hence are equal. Therefore
\[
g \cdot \pi_{ X }^{ * } \of{ \mu } \in g \cdot \pi_{ X }^{ * } \of{ L_{ X } } = h \cdot \pi_{ Y }^{ * } \of{ L_{ Y } } \subseteq \Deck \of{ \pi_{ S } } \cdot \pi_{ Y }^{ * } \of{ \Teich \of{ Y } } ,
\]
and so $ \Deck \of{ \pi_{ S } } \cdot \pi_{ X }^{ * } \of{ \Teich \of{ X } } \subseteq \Deck \of{ \pi_{ S } } \cdot \pi_{ Y }^{ * } \of{ \Teich \of{ Y } } $. By symmetry, the other containment also holds. \qedhere

\end{proof}

\end{claim}

\begin{claim} \label{clm:singledeck}

$ g \cdot \pi_{ X }^{ * } \of{ \Teich \of{ X } } = \pi_{ Y }^{ * } \of{ \Teich \of{ Y } } $ for some deck transformation $ g \in \Deck \of{ \pi_{ S } } $.

\begin{proof}[Proof of \autoref{clm:singledeck}]

For each deck transformation $ g \in \Deck \of{ \pi_{ S } } $ of $ \pi_{ S } \colon Z \to S $, let
\[
I_{ g } \coloneqq \of{ g \cdot \pi_{ X }^{ * } \of{ \Teich \of{ X } } } \cap \pi_{ Y }^{ * } \of{ \Teich \of{ Y } }\subset\Teich(Z) .
\]
Since $ \pi_{ X }^{ * } $ and $ \pi_{ Y }^{ * } $ are proper isometric embeddings and $ g \colon \Teich \of{ Z } \to \Teich \of{ Z } $ is an isometry, it follows that $ g \cdot \pi_{ X }^{ * } \of{ \Teich \of{ X } } $ and $ \pi_{ Y }^{ * } \of{ \Teich \of{ Y } } $ are properly embedded submanifolds of $ \Teich \of{ Z } $. Therefore, $ I_{ g } $ is a properly embedded submanifold of both $ g \cdot \pi_{ X }^{ * } \of{ \Teich \of{ X } } $ and of $ \pi_{ Y }^{ * } \of{ \Teich \of{ Y } } $. Note that
\begin{align*}
\pi_{ Y }^{ * } \of{ \Teich \of{ Y } } & = \of{ \Deck \of{ \pi_{ S } } \cdot \pi_{ Y }^{ * } \of{ \Teich \of{ Y } } } \cap \pi_{ Y }^{ * } \of{ \Teich \of{ Y } } = \of{ \Deck \of{ \pi_{ S } } \cdot \pi_{ X }^{ * } \of{ \Teich \of{ X } } } \cap \pi_{ Y }^{ * } \of{ \Teich \of{ Y } } \\
& = \of{ \bigcup_{ g }{ \of{ g \cdot \pi_{ X }^{ * } \of{ \Teich \of{ X } } } } } \cap \pi_{ Y }^{ * } \of{ \Teich \of{ Y } } = \bigcup_{ g }{ I_{ g } } ,
\end{align*}
where in the above unions $ g $ ranges over elements of $\Deck(\pi_S)$. So $ \pi_{ Y }^{ * } \of{ \Teich \of{ Y } } $ is covered by finitely many properly embedded submanifolds.

For each deck transformation $ g \in \Deck \of{ \pi_{ S } } $, if the containment $ I_{ g } \subsetneq \pi_{ Y }^{ * } \of{ \Teich \of{ Y } } $ is strict, then $ I_{ g } $ has positive codimension in $ \pi_{ Y }^{ * } \of{ \Teich \of{ Y } } $, since it is a properly embedded submanifold of a connected manifold. Since $ \pi_{ Y }^{ * } \of{ \Teich \of{ Y } } $ is covered by finitely many of these submanifolds, this cannot be the case for all such deck transformations and so $ \pi_{ Y }^{ * } \of{ \Teich \of{ Y } } = I_{ g } $ for some deck transformation $ g \in \Deck \of{ \pi_{ S } } $. In particular,
\begin{align*}
\dim \of{ \pi_{ Y }^{ * } \of{ \Teich \of{ Y } } } & = \dim \of{ I_{ g } } = \dim \of{ \of{ g \cdot \pi_{ X }^{ * } \of{ \Teich \of{ X } } } \cap \pi_{ Y }^{ * } \of{ \Teich \of{ Y } } } \\
& \leq \dim \of{ g \cdot \pi_{ X }^{ * } \of{ \Teich \of{ X } } } = \dim \of{ \pi_{ X }^{ * } \of{ \Teich \of{ X } } } .
\end{align*}
By symmetry, the reverse inequality also holds, so that the above inequality is an equality. Thus $ \pi_{ Y }^{ * } \of{ \Teich \of{ Y } } = I_{ g } $ is a properly embedded open submanifold of $ g \cdot \pi_{ X }^{ * } \of{ \Teich \of{ X } } $. Since $ g \cdot \pi_{ X }^{ * } \of{ \Teich \of{ X } } $ is connected, this implies that $ g \cdot \pi_{ X }^{ * } \of{ \Teich \of{ X } } = \pi_{ Y }^{ * } \of{ \Teich \of{ Y } } $. \qedhere

\end{proof}

\end{claim}

For the remainder of the proof of \autoref{thm:non-effective}, fix a deck transformation $ g \in \Deck \of{ \pi_{ S } } $ so that $ g \cdot \pi_{ X }^{ * } \of{ \Teich \of{ X } } = \pi_{ Y }^{ * } \of{ \Teich \of{ Y } } $, the existence of which is guaranteed by \autoref{clm:singledeck}.

\begin{claim} \label{clm:conj}

$ g \Deck \of{ \pi_{ X } } = \Deck \of{ \pi_{ Y } } g $. 

\begin{proof}[Proof of \autoref{clm:conj}]

Consider a deck transformation $ h \in \Deck \of{ \pi_{ X } } $, and note that for each point $ \mu \in \Teich \of{ Y } $, $ g^{ -1 } \cdot \pi_{ Y }^{ * } \of{ \mu } \in \pi_{ X }^{ * } \of{ \Teich \of{ X } } $, so that
\[
g h g^{ -1 } \cdot \pi_{ Y }^{ * } (\mu) = g \cdot \of{ h \cdot \of{ g^{ -1 } \cdot \pi_{ Y }^{ * } (\mu) } } = g \cdot \of{ g^{ -1 } \cdot \pi_{ Y }^{ * } (\mu) } = \pi_{ Y }^{ * } (\mu) .
\]
By \autoref{prop:deckGroup} above, this implies that $ g h g^{ -1 } \in \Deck \of{ \pi_{ Y } } $, so that $ g \Deck \of{ \pi_{ X } } \subseteq \Deck \of{ \pi_{ Y } } g $. Conversely, consider a deck transformation $ k \in \Deck \of{ \pi_{ Y } } $, and note that for each point $ \mu \in \Teich \of{ X } $, $ g \cdot \pi_{ X }^{ * } (\mu) \in \pi_{ Y }^{ * } \Teich \of{ Y } $, so that
\[
g^{ -1 } k g \cdot \pi_{ X }^{ * } (\mu) = g^{ -1 } \cdot \of{ k \cdot \of{ g \cdot \pi_{ X }^{ * } (\mu) } } = g^{ -1 } \cdot \of{ g \cdot \pi_{ X }^{ * } (\mu) } = \pi_{ X }^{ * } (\mu) .
\]
Again by \autoref{prop:deckGroup}, this implies that $ g^{ -1 } k g \in \Deck \of{ \pi_{ X } } $, so that the reverse containment $ \Deck \of{ \pi_{ Y } } g \subseteq g \Deck \of{ \pi_{ X } } $ also holds. \qedhere

\end{proof}

\end{claim}

We are now ready to prove that $ p \colon X \to S $ and $ q \colon Y \to S $ are isomorphic as covers of $ S $. To that end, consider the map $ \widetilde{ \varphi } \coloneqq \pi_{ Y } \circ g \colon Z \to Y $. \autoref{clm:conj} implies that for all deck transformations $ h \in \Deck \of{ \pi_{ X } } $, there is a deck transformation $ k \in \Deck \of{ \pi_{ Y } } $ with $ g h = k g $, so that
\begin{align*}
\widetilde{ \varphi } \of{ h \of{ z } } = \pi_{ Y } \of{ g \of{ h \of{ z } } } = \pi_{ Y } \of{ k \of{ g \of{ z } } } = \pi_{ Y } \of{ g \of{ z } } = \widetilde{ \varphi } \of{ z }
\end{align*}
for all points $ z \in Z $. In other words, $ \widetilde{ \varphi } $ is invariant under the action of the deck group $ \Deck \of{ \pi_{ X } } $, and so descends to a continuous map $ \varphi \colon X \to Y $ such that $ \varphi \circ \pi_{ X } = \widetilde{ \varphi } = \pi_{ Y } \circ g $. Similarly, the map $ \widetilde{ \psi } \coloneqq \pi_{ X } \circ g^{ -1 } \colon Z \to X $ is invariant under the action of the deck group $ \Deck \of{ \pi_{ Y } } $, and so descends to a continuous map $ \psi \colon Y \to X $ such that $ \psi \circ \pi_{ Y } = \widetilde{ \psi } = \pi_{ X } \circ g^{ -1 } $. Note that
\begin{align*}
 \pi_{ X } & = \pi_{ X } \circ g^{ -1 } \circ g = \psi \circ \pi_{ Y } \circ g = \psi \circ \varphi \circ \pi_{ X }
\end{align*}
and
\begin{align*}
\pi_{ Y } & = \pi_{ Y } \circ g \circ g^{ -1 } = \varphi \circ \pi_{ X } \circ g^{ -1 } = \varphi \circ \psi \circ \pi_{ Y } .
\end{align*}
Since $ \pi_{ X } $ and $ \pi_{ Y } $ are surjective, this implies that $ \psi \circ \varphi = \id_{ X } $ and $ \varphi \circ \psi = \id_{ Y } $. Thus $ \varphi $ and $ \psi $ are mutually inverse homeomorphisms. Note that
\begin{align*}
p \circ \pi_{ X } = \pi_{ S } = \pi_{ S } \circ g = q \circ \pi_{ Y } \circ g = q \circ \varphi \circ \pi_{ X } .
\end{align*}
Since $ \pi_{ X } $ is surjective, this implies that $ p = q \circ \varphi $ and so $ \varphi $ is an isomorphism of covers.

Conversely, now suppose that $ p $ and $ q $ are isomorphic as covers of $ S $, so that $ p = q \circ \varphi $ for some homeomorphism $ \varphi \colon X \to Y $ and let $ \gamma $ be a primitive curve on $ S $. If $ \alpha $ is a simple elevation of $ \gamma $ along $ p $ to $ X $, then $ \varphi \of{ \alpha } $ is a simple elevation of $ \gamma $ along $ q $ to $ Y $. On the other hand, if $ \beta $ is a simple elevation of $ \gamma $ along $ q $ to $ Y $, then $ \varphi^{ -1 } \of{ \beta } $ is a simple elevation of $ \gamma $ along $ p $ to $ X $. \qedhere

\setcounter{claim}{0}

\end{proof}

The following corollary can be derived fairly directly from \autoref{thm:non-effective} and the partial proof of \autoref{thm:Sunada} found in \cite[291-294]{Bus10}.

\begin{corollary*} \label{cor:Sunada}

Suppose that two closed, orientable surfaces $ X $ and $ Y $ arise from Sunada's construction. If the resulting length-preserving bijection between curves on $ X $ and curves on $ Y $ restricts to a bijection $ \CC \of{ X } \to \CC \of{ Y } $, then the constructing subgroups are conjugate and hence $ X $ and $ Y $ are isometric over $ S $.

\begin{proof}

Let $ p \colon \of{ X , x } \to \of{ S , s } $ and $ q \colon \of{ Y , y } \to \of{ S , s } $ be finite-degree based covers of a closed surface $ S $ of negative Euler characteristic arising from Sunada's construction on a surjective homomorphism $ \rho \colon \pi_{ 1 } \of{ S , s } \onto G $ and almost conjugate subgroups $ A , B \leq G $. \autoref{thm:Sunada} provides a bijection between the set of curves on $ X $ and the set of curves on $ Y $. Moreover, as can be seen from the proof in \cite[291-294]{Bus10}, for each curve $ \gamma $ on $ X $ and each positive integer $ d \geq 1 $, this bijection restricts to a bijection between the set of degree $ d $ elevations of $ \gamma $ along $ p $ to $ X $ and the set of degree $ d $ elevations of $ \gamma $ along $ q $ to $ Y $.

Suppose that the above bijection restricts to a bijection $ \psi \colon \CC \of{ X } \to \CC \of{ Y } $, and consider a primitive curve $ \gamma $ on $ S $. If $ \alpha \in \CC \of{ X } $ is a simple elevation of $ \gamma $ along $ p $ to $ X $, then $ \psi \of{ \alpha } \in \CC \of{ Y } $ is a simple elevation of $ \gamma $ along $ q $ to $ Y $. Conversely, if $ \beta \in \CC \of{ Y } $ is a simple elevation of $ \gamma $ along $ q $ to $ Y $, then $ \psi^{ -1 } \of{ \beta } \in \CC \of{ X } $ is a simple elevation of $ \gamma $ along $ p $ to $ X $.

\autoref{thm:non-effective} now implies that $ p \colon X \to S $ and $ q \colon Y \to S $ are isomorphic as covers of $ S $. The Galois correspondence yields an element $ \gamma \in \pi_{ 1 } \of{ S , s } $ such that $ \gamma \of{ p_{ * } \of{ \pi_{ 1 } \of{ X , x } } } \gamma^{ -1 } = q_{ * } \of{  \pi_{ 1 } \of{ Y , y } } $. Thus $ \rho \of{ \gamma } A \rho \of{ \gamma }^{ -1 } = B $. Moreover, any isomorphism of covers is an isometry with respect to the pullbacks of any Riemannian metric on $ S $. \qedhere

\end{proof}

\end{corollary*}

\section{Applications to simple length isospectrality} \label{sec:simpleLengthSpectra}

Recall that we call two elements of a group \emph{trace twins} if they are trace-commensurable but no non-zero powers of these elements are conjugate. When the group in question is a closed surface group, we will call the curves on the underlying surface represented by such a pair of elements \emph{(hyperbolic) length twins}. \autoref{thm:equivalence} gives an equivalent characterization of these curves; a pair of curves are length twins if they are length-commensurable but not powers of the same primitive curve.

Our current methods of distinguishing covers by their simple length spectra rely on finding curves with no length twins and comparing their elevations to each cover. \autoref{prop:lengthSpectra}, which we restate below for the reader's convenience, makes this notion precise.

\lengthSpectra*

\begin{proof}[Proof of \autoref{prop:lengthSpectra}]

Since properties (i) and (ii) also hold (possibly with a different positive integer $ m $) if we replace $ \gamma $ with a primitive element which has $ \gamma $ as a non-zero power, we may assume without loss of generality that $ \gamma $ is primitive. Suppose that $ \gamma $ has $ r \geq 1 $ simple elevations of degree $ m $ along $ p $ to $ X $, which is strictly more than the number of simple elevations of degree $ m $ along $ q $ to $ Y $.

Let $ \mathscr{ U } $ be the subset of $ \Teich \of{ S } $ consisting of all hyperbolic metrics on $ S $ for which $ X $ and $ Y $ are simple length isospectral. We will show that $ \mathscr{ U } $ has null Lebesgue measure in $ \Teich \of{ S } $, and that its complement is dense. To that end, for each primitive curve $ \delta $ on $ S $ and each positive integer $ n \geq 1 $, let $ Z_{ n } \of{ \delta } $ be the zero set of the function $ m \cdot \len_{ \of{ \oparg } } \of{ \gamma } - n \cdot \len_{ \of{ \oparg } } \of{ \delta } $; that is,
\[
Z_{ n } \of{ \delta } = \set{ \mu \in \Teich \of{ S } : m \cdot \len_{ \mu } \of{ \gamma } = n \cdot \len_{ \mu } \of{ \delta } } .
\]
This function $ m \cdot \len_{ \of{ \oparg } } \of{ \gamma } - n \cdot \len_{ \of{ \oparg } } \of{ \delta } $ is real-analytic on $ \Teich \of{ S } $, and since $ \gamma $ is primitive and has no length twins, it can only be identically zero if $ n $ divides $ m $ and $ \delta $ is the $ \of{ \frac{ m }{ n } } $th power of $ \gamma $. For all other curves $ \delta $, $ Z_{ n } \of{ \delta } $ is a set of null Lebesgue measure whose complement is open and dense. In particular, as a countable union of null measure sets,
\[
Z \coloneqq \bigcup_{ \delta \neq \gamma } \bigcup_{ n = 1 }^{ \infty }{ Z_{ n } \of{ \delta } } ,
\]
has null Lebesgue measure, where in the above union, $ \delta $ varies over all primitive curves except $ \gamma $. Moreover, the complement of $ Z $ is dense, since $ \Teich \of{ S } $ is a Baire space. It therefore suffices to show that $ \mathscr{ U } \subseteq Z $.

To that end, let $ \mu \in \mathscr{ U } $ and $ L = m \cdot \len_{ \mu } \of{ \gamma } $. Note that each elevation $ \alpha $ of $ \gamma $ of degree $ m $ along $ p $ to $ X $ has length
\[
\len_{ p^{ * } \of{ \mu } } \of{ \alpha } = m \cdot \len_{ \mu } \of{ \gamma } = L .
\]
In particular, $ \of{ X , p^{ * } \of{ \mu } } $ has at least $ r $ simple curves of length $ L $. Since $ \of{ X , p^{ * } \of{ \mu } } $ and $ \of{ Y , q^{ * } \of{ \mu } } $ are simple length isospectral, this implies that $ \of{ Y , q^{ * } \of{ \mu } } $ also has at least $ r $ simple curves of length $ L $. However, since $ \gamma $ does not have $ r $ simple elevations of degree $ m $ along $ q $ to $ Y $, some of these curves are not elevations of $ \gamma $ along $ q $.

Fix such a simple curve $ \beta $ on $ Y $ which has length $ \len_{ q^{ * } \of{ \mu } } \of{ \beta } = L $ but is not an elevation of $ \gamma $ along $ q $. Let $ \delta $ be the primitive curve corresponding to $ q \of{ \beta } $, so that $ \beta $ is a simple elevation of $ \delta $ along $ q $ to $ Y $ of degree some positive integer $ 1 \leq n \leq \deg \of{ q } $. Then
\begin{align*}
m \cdot \len_{ \mu } \of{ \gamma } = \len_{ p^{ * } \of{ \mu } } \of{ \alpha } = L = \len_{ q^{ * } \of{ \mu } } \of{ \beta } = n \cdot \len_{ \mu } \of{ \delta } .
\end{align*}
$ \gamma \neq \delta $ by assumption, and so $ \mu \in Z_{ n } \of{ \delta } \subseteq Z $. Thus $ \mathscr{ U } \subseteq Z $, and so $ \mathscr{ U } $ has null Lebesgue measure in $ \Teich \of{ S } $ and its complement is dense. \qedhere

\end{proof}

In light of \autoref{prop:lengthSpectra}, we would like to find examples of curves with no length twins. \autoref{sec:noLengthTwins} is dedicated to constructing many examples of such curves.

Given a finite-degree cover $ p \colon X \to S $ of $ S $, we will denote by $ \mathscr{ L } \of{ p } $ the set of simple elevations of curves on $ S $ with no length twins. More knowledge of the distribution of $ \mathscr{ L } \of{ p } $ in $ \CC \of{ X } $ has the potential to strengthen \autoref{thm:non-effective}. Towards this end, we pose the following question.

\begin{question*} \label{que:twinlessAccumulation}

Given a finite-degree cover $ p \colon X \to S $, which points in $ \EL \of{ X } = \partial \CC \of{ X } $ are accumulation points of $ \mathscr{ L } \of{ p } $?

\end{question*}

If the answer to \autoref{que:twinlessAccumulation} is that $ \mathscr{ L } \of{ p } $ accumulates at every point of $ \EL \of{ X } $, then we may replace $ \CC \of{ X } $ and $ \CC \of{ Y } $ with $ \mathscr{ L } \of{ p } $ and $ \mathscr{ L } \of{ q } $, respectively, in \autoref{clm:CC} of the proof of \autoref{thm:non-effective} in \autoref{sec:mainthm}. Combined with \autoref{prop:lengthSpectra} above, we would thus obtain an affirmative answer to \autoref{que:genericIsospectrality}, which we restate below for the reader's convenience:

\OurQuestion*

Even without an answer to \autoref{que:twinlessAccumulation}, we can still provide partial affirmative answers to \autoref{que:genericIsospectrality} for simply generated covers in \autoref{sec:simplyGenerated}. Further evidence for an affirmative answer is given in \autoref{sec:examples}, where we show that every closed surface of negative Euler characteristic has covers that are length isospectral but generically not simple length isospectral.

\subsection{Curves with no length twins} \label{sec:noLengthTwins}

As outlined above, a main component of our program of distinguishing covers involves finding curves which have no length twins. The goal of this section is to prove the following proposition, which provides many examples of such curves.

\begin{proposition*} \label{prop:extendedHorowitz}

Let $ S $ be a closed surface of genus $ g \geq 2 $, and suppose that elements $ \alpha , \beta \in \pi_{ 1 } \of{ S , s } $ representing distinct, disjoint, simple, non-separating curves have injective parametrized representatives $ \overline{ \alpha } $ and $ \overline{ \beta } $ meeting exactly once (at the basepoint $ s $).

\begin{enumerate}

\item[1.]

If $ g = 2 $, then for any integers $ j , k \in \Z $, the only length twin of the curve represented by $ \alpha^{ j } \beta^{ k } $ is its image under the hyperelliptic involution; and

\item[2.]

If $ g \geq 3 $, then for any integers $ j , k \in \Z $, the curve represented by $ \alpha^{ j } \beta^{ k } $ has no length twins.

\end{enumerate}

\end{proposition*}

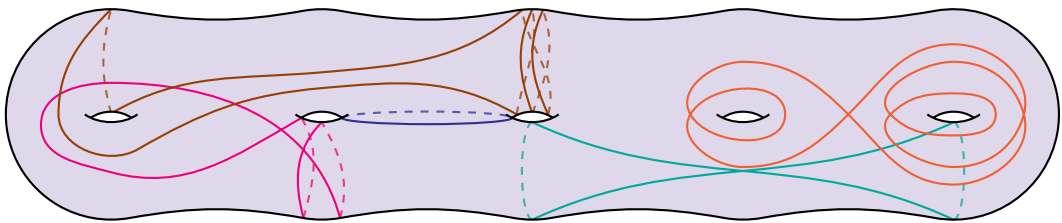
\begin{figure}[ht]
\centering
\begin{tikzpicture}[ scale = 1.4 ]

\fill[ RoyalPurple , opacity = 0.2 ]
    ( 0 , 1 ) arc[ radius = 1 , start angle = 90 , end angle = { 270 + 10 } ]
    to[ out = 10 , in = { 180 - 10 } ] ( 2 - 0.1736 , -0.9848 )
    arc[ radius = 1 , start angle = { -90 - 10 } , end angle = { -90 + 10 } ]
    to[ out = 10 , in = { 180 - 10 } ] ( 4 - 0.1736 , -0.9848 )
    arc[ radius = 1 , start angle = { -90 - 10 } , end angle = { -90 + 10 } ]
    to[ out = 10 , in = { 180 - 10 } ] ( 6 - 0.1736 , -0.9848 )
    arc[ radius = 1 , start angle = { -90 - 10 } , end angle = { -90 + 10 } ]
    to[ out = 10 , in = { 180 - 10 } ] ( 8 - 0.1736 , -0.9848 )
    arc[ radius = 1 , start angle = { -90 - 10 } , end angle = { 90 + 10 } ]
    to[ out = { 180 + 10 } , in = -10 ] ( 6 + 0.1736 , 0.9848 )
    arc[ radius = 1 , start angle = { 90 - 10 } , end angle = { 90 + 10 } ]
    to[ out = { 180 + 10 } , in = -10 ] ( 4 + 0.1736 , 0.9848 )
    arc[ radius = 1 , start angle = { 90 - 10 } , end angle = { 90 + 10 } ]
    to[ out = { 180 + 10 } , in = -10 ] ( 2 + 0.1736 , 0.9848 )
    arc[ radius = 1 , start angle = { 90 - 10 } , end angle = { 90 + 10 } ]
    to[ out = { 180 + 10 } , in = -10 ] ( 0 + 0.1736 , 0.9848 )
    arc[ radius = 1 , start angle = { 90 - 10 } , end angle = 90 ] ;

\draw[ thick , RubineRed ]
    ( 2 - 0.1926 , -0.0329 ) to[ out = 210 , in = 345 ] ( 0 , -0.55 )
    to[ out = 165 , in = 270 ] ( -2 / 3 , -1 / 10 )
    to[ out = 90 , in = 180 ] ( 0 , 3 / 10 )
    to[ out = 0 , in = 105 ] ( 2 + 0.1736 , -0.9848 )
    ( 2 , -0.075 ) to[ out = 225 , in = 105 ] ( 2 - 0.1736 , -0.9848 ) ;
\draw[ thick , dashed , RubineRed , opacity = 0.75 ]
    ( 2 + 0.1736 , -0.9848 ) to[ out = 75 , in = 300 ] ( 2 , -0.075 )
    ( 2 - 0.1736 , -0.9848 ) to[ out = 75 , in = 300 ] ( 2 - 0.1926 , -0.0329 ) ;

\draw[ thick , Blue ]
    ( 2 + 0.1926 , -0.0329 ) arc[ x radius = { 1 - 0.1926 } , y radius = 0.06 , start angle = 180 , end angle = 360 ] ;
\draw[ thick , dashed , Blue , opacity = 0.75 ]
    ( 2 + 0.1926 , -0.0329 ) arc[ x radius = { 1 - 0.1926 } , y radius = 0.06 , start angle = 180 , end angle = 0 ] ;

\draw[ thick , Emerald ]
    ( 4 , -1 ) to[ out = 25 , in = 205 ] ( 8 , 0 - 0.075 )
    ( 8 , -1 ) to[ out = 155 , in = 335 ] ( 4 , 0 - 0.075 ) ;
\draw[ thick , dashed , Emerald , opacity = 0.75 ]
    ( 4 , -1 ) arc[ x radius = 0.1 , y radius = { ( 1 - 0.075 ) / 2 } , start angle = 270 , end angle = 90 ]
    ( 8 , -1 ) arc[ x radius = 0.1 , y radius = { ( 1 - 0.075 ) / 2 } , start angle = -90 , end angle = 90 ] ;

\draw[ thick , RawSienna ]
    ( 0 , 1 ) to[ out = 225 , in = 90 ] ( -1 / 2 , 0 )
        to[ out = 270 , in = 210 ] ( 1 / 4 , -1 / 3 )
        to[ out = 30 , in = 185 ] ( 2 , 0.25 )
        to[ out = 5 , in = 150 ] ( 4 - 0.15 , 0 )
    ( 4 , 1 ) to[ out = 245 , in = 110 ] ( 4 , 0.01 )
    ( 4 + 0.1 , 1 - 0.01 ) to[ out = 245 , in = 115 ] ( 4 + 0.15 , 0 )
    ( 4 - 0.1 , 1 - 0.01 ) to[ out = 220 , in = 5 ] ( 2 , 0.4 )
        to[ out = 185 , in = 30 ] ( 0 , 0.01 ) ;
\draw[ thick , dashed , RawSienna , opacity = 0.75 ]
    ( 4 - 0.15 , 0 ) to[ out = 80 , in = 280 ] ( 4 , 1 )
    ( 4 , 0.01 ) to[ out = 70 , in = 290 ] ( 4 + 0.1 , 1 - 0.01 )
    ( 4 + 0.15 , 0 ) to[ out = 75 , in = 280 ] ( 4 - 0.1 , 1 - 0.01 )
    ( 0 , 0.01 ) to[ out = 105 , in = 255 ] ( 0 , 1 ) ;

\draw[ thick , RedOrange ]
    ( 7 , 0 ) to[ out = 135 , in = 0 ] ( 6 , 1 / 2 )
    to[ out = 180 , in = 120 ] ( 6 - 1 / 2 , 0 )
    to[ out = 300 , in = 270 ] ( 6 + 2 / 5 , 0 )
    to[ out = 90 , in = 60 ] ( 6 - 1 / 2 , 0 )
    to[ out = 240 , in = 180 ] ( 6 , -1 / 2 )
    to[ out = 0 , in = 225 ] ( 7 , 0 )
    to[ out = 45 , in = 180 ] ( 8 , 2 / 3 )
    to[ out = 0 , in = 75 ] ( 8 + 2 / 3 , 0 )
    to[ out = 255 , in = 0 ] ( 8 , -1 / 2 )
    to[ out = 180 , in = 225 ] ( 8 - 3 / 5 , 0 )
    to[ out = 45 , in = 180 ] ( 8 , 1 / 5 )
    to[ out = 0 , in = 90 ] ( 8 + 2 / 5 , 0 )
    to[ out = 270 , in = 0 ] ( 8 , -1 / 5 )
    to[ out = 180 , in = 315 ] ( 8 - 3 / 5 , 0 )
    to[ out = 135 , in = 180 ] ( 8 , 1 / 2 )
    to[ out = 0 , in = 105 ] ( 8 + 2 / 3 , 0 )
    to[ out = 285 , in = 0 ] ( 8 , -2 / 3 )
    to[ out = 180 , in = 315 ] ( 7 , 0 ) ;

\fill[ White ]
    ( 0 - 0.1926 , -0.0329 ) arc[ x radius = 0.2723 , y radius = 0.2 , start angle = 135 , end angle = 45 ]
    arc[ x radius = 0.35355 , y radius = 0.25 , start angle = 303 , end angle = 237 ] ;
\draw[ thick ]
    ( 0 - 0.25 , 0 ) arc[ x radius = 0.35355 , y radius = 0.25 , start angle = 225 , end angle = 315 ]
    ( 0 - 0.1926 , -0.0329 ) arc[ x radius = 0.2723 , y radius = 0.2 , start angle = 135 , end angle = 45 ] ;

\fill[ White ]
    ( 2 - 0.1926 , -0.0329 ) arc[ x radius = 0.2723 , y radius = 0.2 , start angle = 135 , end angle = 45 ]
    arc[ x radius = 0.35355 , y radius = 0.25 , start angle = 303 , end angle = 237 ] ;
\draw[ thick ]
    ( 2 - 0.25 , 0 ) arc[ x radius = 0.35355 , y radius = 0.25 , start angle = 225 , end angle = 315 ]
    ( 2 - 0.1926 , -0.0329 ) arc[ x radius = 0.2723 , y radius = 0.2 , start angle = 135 , end angle = 45 ] ;

\fill[ White ]
    ( 4 - 0.1926 , -0.0329 ) arc[ x radius = 0.2723 , y radius = 0.2 , start angle = 135 , end angle = 45 ]
    arc[ x radius = 0.35355 , y radius = 0.25 , start angle = 303 , end angle = 237 ] ;
\draw[ thick ]
    ( 4 - 0.25 , 0 ) arc[ x radius = 0.35355 , y radius = 0.25 , start angle = 225 , end angle = 315 ]
    ( 4 - 0.1926 , -0.0329 ) arc[ x radius = 0.2723 , y radius = 0.2 , start angle = 135 , end angle = 45 ] ;

\fill[ White ]
    ( 6 - 0.1926 , -0.0329 ) arc[ x radius = 0.2723 , y radius = 0.2 , start angle = 135 , end angle = 45 ]
    arc[ x radius = 0.35355 , y radius = 0.25 , start angle = 303 , end angle = 237 ] ;
\draw[ thick ]
    ( 6 - 0.25 , 0 ) arc[ x radius = 0.35355 , y radius = 0.25 , start angle = 225 , end angle = 315 ]
    ( 6 - 0.1926 , -0.0329 ) arc[ x radius = 0.2723 , y radius = 0.2 , start angle = 135 , end angle = 45 ] ;

\fill[ White ]
    ( 8 - 0.1926 , -0.0329 ) arc[ x radius = 0.2723 , y radius = 0.2 , start angle = 135 , end angle = 45 ]
    arc[ x radius = 0.35355 , y radius = 0.25 , start angle = 303 , end angle = 237 ] ;
\draw[ thick ]
    ( 8 - 0.25 , 0 ) arc[ x radius = 0.35355 , y radius = 0.25 , start angle = 225 , end angle = 315 ]
    ( 8 - 0.1926 , -0.0329 ) arc[ x radius = 0.2723 , y radius = 0.2 , start angle = 135 , end angle = 45 ] ;

\draw[ thick ]
    ( 0 , 1 ) arc[ radius = 1 , start angle = 90 , end angle = { 270 + 10 } ]
    to[ out = 10 , in = { 180 - 10 } ] ( 2 - 0.1736 , -0.9848 )
    arc[ radius = 1 , start angle = { -90 - 10 } , end angle = { -90 + 10 } ]
    to[ out = 10 , in = { 180 - 10 } ] ( 4 - 0.1736 , -0.9848 )
    arc[ radius = 1 , start angle = { -90 - 10 } , end angle = { -90 + 10 } ]
    to[ out = 10 , in = { 180 - 10 } ] ( 6 - 0.1736 , -0.9848 )
    arc[ radius = 1 , start angle = { -90 - 10 } , end angle = { -90 + 10 } ]
    to[ out = 10 , in = { 180 - 10 } ] ( 8 - 0.1736 , -0.9848 )
    arc[ radius = 1 , start angle = { -90 - 10 } , end angle = { 90 + 10 } ]
    to[ out = { 180 + 10 } , in = -10 ] ( 6 + 0.1736 , 0.9848 )
    arc[ radius = 1 , start angle = { 90 - 10 } , end angle = { 90 + 10 } ]
    to[ out = { 180 + 10 } , in = -10 ] ( 4 + 0.1736 , 0.9848 )
    arc[ radius = 1 , start angle = { 90 - 10 } , end angle = { 90 + 10 } ]
    to[ out = { 180 + 10 } , in = -10 ] ( 2 + 0.1736 , 0.9848 )
    arc[ radius = 1 , start angle = { 90 - 10 } , end angle = { 90 + 10 } ]
    to[ out = { 180 + 10 } , in = -10 ] ( 0 + 0.1736 , 0.9848 )
    arc[ radius = 1 , start angle = { 90 - 10 } , end angle = 90 ] ;

\end{tikzpicture}
\caption{Curves with no length twins} \label{fig:noLengthTwins}
\end{figure}

For example, each of the curves drawn in \autoref{fig:noLengthTwins} can be realized by elements of the corresponding surface group which satisfy the hypotheses of \autoref{prop:extendedHorowitz}, and therefore have no length twins.

Before beginning the proof of \autoref{prop:extendedHorowitz}, we will establish two preliminary lemmas, the first of which states that powers of simple curves have no length twins.

\begin{lemma*} \label{lem:simpleLengthTwins}

Powers of simple curves have no length twins.

\end{lemma*}

\begin{proof}

Let $ \alpha $ and $ \beta $ be a pair of length-commensurable curves on a closed surface $ S $ of negative Euler characteristic. Suppose that $ \alpha $ is a power of a simple curve $ \gamma $, $ \beta $ is a power of a primitive curve $ \delta $, and that $ \gamma $ and $ \delta $ are length-commensurable.

Since $ \gamma $ and $ \delta $ are primitive and $ i \of{ \gamma , \gamma } = 0 $, \autoref{thm:equivalence} implies that $ i \of{ \gamma , \delta } = 0 $. If $ \gamma \neq \delta $, take a simple curve $ \eta \neq \gamma $ obtained from $ \delta $ via surgery. Note that $ i \of{ \gamma , \eta } = 0 $, and so Fenchel-Nielsen coordinates yield, for any two positive real numbers $ c $ and $ e $, a hyperbolic metric $ \mu $ on $ S $ for which $ \len_{ \mu } \of{ \gamma } = c $ and $ \len_{ \mu } \of{ \eta } = e $. Note that $ \len_{ \mu } \of{ \delta } \geq \len_{ \mu } \of{ \eta } = e $. This contradicts our assumption that $\gamma$ and $\delta$ are length-commensurable, and thus $ \delta = \gamma $. We have shown that any curve which is length-commensurable to $ \alpha $ is also a power of $\gamma$, and so $ \alpha $ has no length twins. \qedhere

\end{proof}

Recall \autoref{thm:equivalence}, which states that length-commensurable curves on a surface are represented by trace-commensurable elements in the corresponding surface group. Finding curves without length twins is thus equivalent to finding elements with no trace twins. The following argument augments Horowitz's proof of \autoref{thm:horowitz} in \cite[Corollary~7.2]{Hor72} to find examples of such elements in free groups.

\begin{lemma*} \label{lem:extendedHorowitz}

For any integers $ j , k \in \Z $ and primitive elements $ a, b \in F$ of a free group $ F $, the element $ w = a^{ j } b^{ k } $ has no trace twins in $ F $.

\end{lemma*}

\begin{proof}

Note that \cite[Theorem~1.1]{Hor72} implies that we may assume that $ F $ has finite rank. If $ j = 0 $, $ k = 0 $, or $ a = b $, then the desired result follows immediately from \autoref{thm:horowitz}. Thus we may assume that $ j $ and $ k $ are both non-zero and that $ a \neq b $. Suppose that $ w = a^{ j } b^{ k } $ and some element $ x \in F $ are trace-commensurable, so that
\[
\tr \of!\big!{ \underbrace{ a^{ j } b^{ k } \dotsm a^{ j } b^{ k } }_{ m } }^{ 2 } = \tr \of{ \of!\big!{ a^{ j } b^{ k } }^{ m } }^{ 2 } = \tr \of{ w^{ m } }^{ 2 } = \tr \of{ x^{ n } }^{ 2 }
\]
for some non-zero integers $ m , n \in \Z \setminus \set{ 0 } $. Extend $ \set{ a , b } $ to a basis $ B $ for $ F $, and let $ y = x^{ n } $.

\cite[Theorem~3.1]{Hor72} gives a way to express $ \tr \of{ w^{ m } } $ as a polynomial in $ \tr \of{ a } $, $ \tr \of{ b } $, and $ \tr \of{ a b } $. Similarly, this theorem gives a way to express $ \tr \of{ y } $ as a polynomial in products of the elements of $ B $ which appear in a freely reduced word representing $ y $. Since these traces must be equal for all linear representations of $ F $, this implies that $ y \in \gen{ a , b } $, so that $ y $ is conjugate to a word of the form
\[
z = a^{ e_{ 1 } } b^{ f_{ 1 } } a^{ e_{ 2 } } b^{ f_{ 2 } } \dotsm a^{ e_{ p } } b^{ f_{ p } }
\]
for some integers $ e_{ 1 } , \dotsc , e_{ p } , f_{ 1 } , \dotsc , f_{ p } \in \Z $, at most one of which is zero. We will show that $ w $ is either $ z $ or $ z^{ -1 } $. Since the trace is invariant under conjugation,
\[
\tr \of!\big!{ \underbrace{ a^{ j } b^{ k } \dotsm a^{ j } b^{ k } }_{ m } }^{ 2 } = \tr \of{ x^{ n } }^{ 2 } = \tr \of{ y }^{ 2 } = \tr \of{ a^{ e_{ 1 } } b^{ f_{ 1 } } \dotsm a^{ e_{ p } } b^{ f_{ p } } }^{ 2 } .
\]
\cite[Lemma~6.1]{Hor72} implies that $ m = p $, $ \abs{ e_{ 1 } } = \dotsb = \abs{ e_{ p } } = \abs{ j } $, and $ \abs{ f_{ 1 } } = \dotsb = \abs{ f_{ p } } = \abs{ k } $. Let $ \varphi \colon F_{ 2 } \to \Z $ be the homomorphism defined by the conditions $ \varphi \of{ a } = \sgn \of{ j } $ and $ \varphi \of{ b } = \sgn \of{ k } $, where
\[
\sgn \of{ t } =
\begin{dcases*}
-1 , & if $ x < 0 $; \\
0 , & if $ x = 0 $; and \\
1 , & otherwise if $ x > 0 $.
\end{dcases*}
\]

Since $ \varphi $ is a homomorphism, $ \tr \of{ \varphi \of{ w^{ m } } }^{ 2 } = \tr \of{ \varphi \of{ z } }^{ 2 } $ as functions on $ \hom \of{ \Z , \SL_{ 2 } \of{ \C } } $, and so $ \abs{ \varphi \of{ w^{ m } } } = \abs{ \varphi \of{ z } } $ by \autoref{thm:horowitz}. Thus
\begin{align*}
\abs{ m } \cdot \abs!\big!{ \abs{ j } + \abs{ k } } & = \abs{ m } \cdot \abs{ j \cdot \sgn \of{ j } + k \cdot \sgn \of{ k } } = \abs!\big!{ m \cdot \varphi \of!\big!{ a^{ j } b^{ k } } } = \abs{ m \cdot \varphi \of{ w } } = \abs{ \varphi \of{ w^{ m } } } \\
& = \abs{ \varphi \of{ z } } = \abs{ \varphi \of{ a^{ e_{ 1 } } b^{ f_{ 1 } } \dotsm a^{ e_{ m } } b^{ f_{ m } } } } = \abs{ \sgn \of{ j } \sum_{ i = 1 }^{ m }{ e_{ i } } + \sgn \of{ k } \sum_{ i = 1 }^{ m }{ f_{ i } } } .
\end{align*}
Since $ \abs{ e_{ 1 } } = \dotsb = \abs{ e_{ m } } = \abs{ j } $ and $ \abs{ f_{ 1 } } = \dotsb = \abs{ f_{ m } } = \abs{ k } $, the only way for this to happen is for there to be no cancellation in the sum appearing in the bottom right of the above equations. Thus there are two cases; either $ e_{ 1 } = \dotsb = e_{ m } = j $ and $ f_{ 1 } = \dotsb = f_{ m } = k $, in which case $ w = z $, or alternatively $ e_{ 1 } = \dotsb = e_{ m } = -j $ and $ f_{ 1 } = \dotsb = f_{ m } = -k $, in which case $ w $ is conjugate to $ z^{ -1 } $.

In summary, we have shown that for any element $ x \in F $, if $ w $ and $ x $ are trace-commensurable, then some non-zero powers of $ w $ and $ x $ are conjugate in $ F $ (and in fact conjugate in $ \gen{ a , b } $). Therefore, $ w = a^{ j } b^{ k } $ has no trace twins in $ F $. \qedhere

\end{proof}

We are now ready to prove \autoref{prop:extendedHorowitz}. Note that we will make use of the second part of \autoref{obs:traceCommensurability}, that the images of trace-commensurable elements of a group under a homomorphism are themselves trace-commensurable.

\begin{proof}[Proof of \autoref{prop:extendedHorowitz}]

For ease of readability, we will suppress the distinction between an element of $ \pi_{ 1 } \of{ S , s } $ and the curve on $ S $ it represents. Under the above hypotheses, there are elements $ \gamma_{ 1 } , \dotsc , \gamma_{ 2 g - 2 } \in \pi_{ 1 } \of{ S , s } $ representing simple curves so that
\[
\pi_{ 1 } \of{ S , s } = \gen{ \alpha , \beta , \gamma_{ 1 } , \dotsc , \gamma_{ 2 g - 2 } \mid \comm{ \alpha }{ \gamma_{ 1 } } \comm{ \beta }{ \gamma_{ 2 } } \comm{ \gamma_{ 3 } }{ \gamma_{ 4 } } \dotsm \comm{ \gamma_{ 2 g - 3 } }{ \gamma_{ 2 g - 2 } } }
\]
is the usual presentation for $ \pi_{ 1 } \of{ S , s } $, the generators of which are drawn in \autoref{fig:generatorsWithPants}.

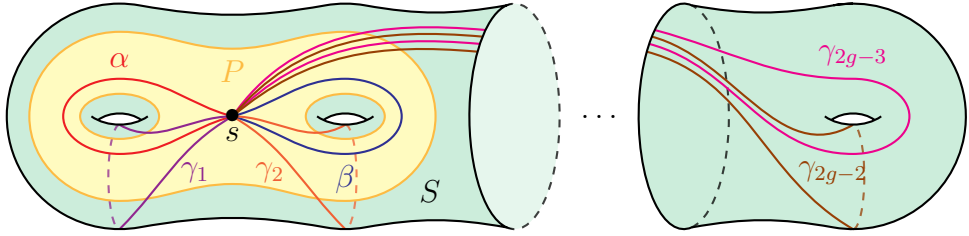
\begin{figure}[ht]
\centering
\begin{tikzpicture}[ scale = 1.5 ]

\fill[ Green , opacity = 0.2 ]
    ( 0 , 1 ) arc[ radius = 1 , start angle = 90 , end angle = { 270 + 10 } ]
    to[ out = 10 , in = { 180 - 10 } ] ( 2 - 0.1736 , -0.9848 )
    arc[ radius = 1 , start angle = { -90 - 10 } , end angle = { -90 + 10 } ]
    to[ out = 10 , in = { 180 - 10 } ] ( 3.5 - 0.4 * 0.1736 , -0.9848 )
    arc[ x radius = { 0.4 * 1 } , y radius = 1 , start angle = { 270 - 10 } , end angle = { 90 + 10 } ]
    to[ out = { 180 + 10 } , in = -10 ] ( 2 + 0.1736 , 0.9848 )
    arc[ radius = 1 , start angle = { 90 - 10 } , end angle = { 90 + 10 } ]
    to[ out = { 180 + 10 } , in = -10 ] ( 0 + 0.1736 , 0.9848 )
    arc[ radius = 1 , start angle = { 90 - 10 } , end angle = 90 ]
    
    ( 6.5 , 1 ) arc[ radius = 1 , start angle = 90 , end angle = { 90 + 10 } ]
    to[ out = { 180 + 10 } , in = -10 ] ( 5 + 0.4 * 0.1736 , 0.9848 )
    arc[ x radius = 0.4 , y radius = 1 , start angle = { 90 - 10 } , end angle = { 270 + 10 } ]
    to[ out = 10 , in = { 180 - 10 } ] ( 6.5 - 0.1736 , -0.9848 )
    arc[ radius = 1 , start angle = { -90 - 10 } , end angle = { 90 + 10 } ] ;

\fill[ Green , opacity = 0.1 ]
    ( 3.5 - 0.4 * 0.1736 , -0.9848 ) arc[ x radius = 0.4 , y radius = 1 , start angle = { -90 - 10 } , end angle = { 270 - 10 } ] ;

\fill[ White ]
    ( 0 , 3 / 4 ) arc[ x radius = 0.8 , y radius = 3 / 4 , start angle = 90 , end angle = 270 ]
    to[ out = 0 , in = 180 ] ( 1 , -3 / 5 )
    to [ out = 0 , in = 180 ] ( 2 , -3 / 4 ) 
    arc[ x radius = 0.8 , y radius = 3 / 4 , start angle = -90 , end angle = 90 ]
    to[ out = 180 , in = 0 ] ( 1 , 3 / 5 )
    to[ out = 180 , in = 0 ] ( 0 , 3 / 4 )
    ( 0 , 0 ) circle[ x radius = 0.35 , y radius = 0.2 ]
    ( 2 , 0 ) circle[ x radius = 0.35 , y radius = 0.2 ] ;

\fill[ Yellow , opacity = 0.25 ]
    ( 0 , 3 / 4 ) arc[ x radius = 0.8 , y radius = 3 / 4 , start angle = 90 , end angle = 270 ]
    to[ out = 0 , in = 180 ] ( 1 , -3 / 5 )
    to [ out = 0 , in = 180 ] ( 2 , -3 / 4 ) 
    arc[ x radius = 0.8 , y radius = 3 / 4 , start angle = -90 , end angle = 90 ]
    to[ out = 180 , in = 0 ] ( 1 , 3 / 5 )
    to[ out = 180 , in = 0 ] ( 0 , 3 / 4 )
    ( 0 , 0 ) circle[ x radius = 0.35 , y radius = 0.2 ]
    ( 2 , 0 ) circle[ x radius = 0.35 , y radius = 0.2 ] ;

\fill[ White ]
    ( 0 , 0 ) circle[ x radius = 0.35 , y radius = 0.2 ]
    ( 2 , 0 ) circle[ x radius = 0.35 , y radius = 0.2 ] ;

\fill[ Green , opacity = 0.2 ]
    ( 0 , 0 ) circle[ x radius = 0.35 , y radius = 0.2 ]
    ( 2 , 0 ) circle[ x radius = 0.35 , y radius = 0.2 ] ;

\draw[ thick , dashed , opacity = 0.75 ]
    ( 3.5 - 0.4 * 0.1736 , -0.9848 ) arc[ x radius = 0.4 , y radius = 1 , start angle = { -90 - 10 } , end angle = { 90 + 10 } ]
    ( 5 + 0.4 * 0.1736 , -0.9848 ) arc[ x radius = 0.4 , y radius = 1 , start angle = { -90 + 10 } , end angle = { 90 - 10 } ] ;

\draw[ thick , dashed , Plum , opacity = 0.75 ]
    ( 0 , -0.075 ) arc[ x radius = 0.1 , y radius = { ( 1 - 0.075 ) / 2 } , start angle = 90 , end angle = 270 ] ;

\draw[ thick , dashed , RedOrange , opacity = 0.75 ]
    ( 2 , -0.075 ) arc[ x radius = 0.1 , y radius = { ( 1 - 0.075 ) / 2 } , start angle = 90 , end angle = -90 ] ;

\draw[ thick , dashed , RawSienna , opacity = 0.75 ]
    ( 6.5 , -1 ) arc[ x radius = 0.1 , y radius = { ( 1 - 0.075 ) / 2 } , start angle = -90 , end angle = 90 ] ;

\draw ( 4.25 , 0 ) node {$ \hdots $} ;

\draw[ thick , Dandelion ]
    ( 0 , 3 / 4 ) arc[ x radius = 0.8 , y radius = 3 / 4 , start angle = 90 , end angle = 270 ]
    to[ out = 0 , in = 180 ] ( 1 , -3 / 5 )
    to [ out = 0 , in = 180 ] ( 2 , -3 / 4 ) 
    arc[ x radius = 0.8 , y radius = 3 / 4 , start angle = -90 , end angle = 90 ]
    to[ out = 180 , in = 0 ] ( 1 , 3 / 5 ) node[ anchor = north ] {$ P $}
    to[ out = 180 , in = 0 ] ( 0 , 3 / 4 )
    ( 0 , 0 ) circle[ x radius = 0.35 , y radius = 0.2 ]
    ( 2 , 0 ) circle[ x radius = 0.35 , y radius = 0.2 ] ;

\draw[ thick , Red ]
    ( 1 , 0 )
    to[ out = 165 , in = 0 ] ( 0 , 1 / 3 ) node[ anchor = south ] {$ \alpha $}
    arc[ x radius = 0.5 , y radius = 1 / 3 , start angle = 90 , end angle = 180 ]
    arc[ x radius = 0.5 , y radius = 1 / 3 , start angle = 180 , end angle = 270 ]
    to[ out = 0 , in = 195 ] ( 1 , 0 ) ;

\draw[ thick , Blue ]
    ( 1 , 0 )
    to[ out = -15 , in = 180 ] ( 2 , -1 / 3 ) node[ anchor = north ] {$ \beta $}
    arc[ x radius = 1 / 2 , y radius = 1 / 3 , start angle = -90 , end angle = 0 ]
    arc[ x radius = 1 / 2 , y radius = 1 / 3 , start angle = 0 , end angle = 90 ]
    to[ out = 180 , in = 15 ] ( 1 , 0 ) ;

\draw[ thick , Plum ]
    ( 1 - 0.56 , -1 / 2 ) node[ anchor = west ] {$ \gamma_{ 1 } $}
    ( 1 , 0 ) to[ out = 180 , in = -30 ] ( 0 , -0.075 )
    ( 0 , -1 ) to[ out = 45 , in = 210 ] ( 1 , 0 ) ;

\draw[ thick , RedOrange ]
    ( 1 + 0.56 , -1 / 2 ) node[ anchor = east ] {$ \gamma_{ 2 } $}
    ( 1 , 0 ) to[ out = 0 , in = 210 ] ( 2 , -0.075 )
    ( 2 , -1 ) to[ out = 135 , in = -30 ] ( 1 , 0 ) ;

\draw[ thick , Magenta ]
    ( 1 , 0 ) to[ out = 45 , in = { 180 - 5 } ] ( 3.5 - 0.4 * 0.7660 , 0.6428 )
    ( 1 , 0 ) to[ out = 55 , in = { 180 - 5 } ] ( 3.5 - 0.4 * 0.6428 , 0.7660 )
    ( 5 - 0.4 * 0.7660 , 0.6428 ) to[ out = -15 , in = 180 ] ( 6.5 , -1 / 3 )
    arc[ x radius = 1 / 2 , y radius = 1 / 3 , start angle = -90 , end angle = 0 ]
    arc[ x radius = 1 / 2 , y radius = 1 / 3 , start angle = 0 , end angle = 90 ] node[ anchor = south ] {$ \gamma_{ 2 g - 3 } $}
    to[ out = 180 , in = -15 ] ( 5 - 0.4 * 0.6428 , 0.7660 ) ;

\draw[ thick , RawSienna ]
    ( 6.5 - 0.585 , -1 / 2 ) node[ anchor = west ] {$ \gamma_{ 2 g - 2 } $}
    ( 1 , 0 ) to[ out = 50 , in = { 180 - 5 } ] ( 3.5 - 0.4 * 0.7071 , 0.7071 )
    ( 1 , 0 ) to[ out = 40 , in = { 180 - 5 } ] ( 3.5 - 0.4 * 0.8192 , 0.5736 )
    ( 5 - 0.4 * 0.7071 , 0.7071 ) to[ out = -15 , in = { 180 + 30 } ] ( 6.5 , -0.075 )
    ( 5 - 0.4 * 0.8192 , 0.5736 ) to[ out = -15 , in = { 180 - 30 } ] ( 6.5 , -1 ) ;

\draw[ thick ]
    ( 0 , 1 ) arc[ radius = 1 , start angle = 90 , end angle = { 270 + 10 } ]
    to[ out = 10 , in = { 180 - 10 } ] ( 2 - 0.1736 , -0.9848 )
    arc[ radius = 1 , start angle = { -90 - 10 } , end angle = { -90 + 10 } ]
    to[ out = 10 , in = { 180 - 10 } ] ( 3.5 - 0.4 * 0.1736 , -0.9848 )
    arc[ x radius = { 0.4 * 1 } , y radius = 1 , start angle = { 270 - 10 } , end angle = { 90 + 10 } ]
    to[ out = { 180 + 10 } , in = -10 ] ( 2 + 0.1736 , 0.9848 )
    arc[ radius = 1 , start angle = { 90 - 10 } , end angle = { 90 + 10 } ]
    to[ out = { 180 + 10 } , in = -10 ] ( 0 + 0.1736 , 0.9848 )
    arc[ radius = 1 , start angle = { 90 - 10 } , end angle = 90 ]
    
    ( 6.5 , 1 ) arc[ radius = 1 , start angle = 90 , end angle = { 90 + 10 } ]
    to[ out = { 180 + 10 } , in = -10 ] ( 5 + 0.4 * 0.1736 , 0.9848 )
    arc[ x radius = 0.4 , y radius = 1 , start angle = { 90 - 10 } , end angle = { 270 + 10 } ]
    to[ out = 10 , in = { 180 - 10 } ] ( 6.5 - 0.1736 , -0.9848 )
    arc[ radius = 1 , start angle = { -90 - 10 } , end angle = { 90 + 10 } ] ;

\fill[ White ]
    ( -0.1926 , -0.0329 ) arc[ x radius = 0.2723 , y radius = 0.2 , start angle = 135 , end angle = 45 ]
    arc[ x radius = 0.35355 , y radius = 0.25 , start angle = 303 , end angle = 237 ]
    ;

\draw[ thick ]
    ( -0.25 , 0 ) arc[ x radius = 0.35355 , y radius = 0.25 , start angle = 225 , end angle = 315 ]
    ( -0.1926 , -0.0329 ) arc[ x radius = 0.2723 , y radius = 0.2 , start angle = 135 , end angle = 45 ] ;

\fill[ White ]
    ( 2 - 0.1926 , -0.0329 ) arc[ x radius = 0.2723 , y radius = 0.2 , start angle = 135 , end angle = 45 ]
    arc[ x radius = 0.35355 , y radius = 0.25 , start angle = 303 , end angle = 237 ]
    ;

\draw[ thick ]
    ( 2 - 0.25 , 0 ) arc[ x radius = 0.35355 , y radius = 0.25 , start angle = 225 , end angle = 315 ]
    ( 2 - 0.1926 , -0.0329 ) arc[ x radius = 0.2723 , y radius = 0.2 , start angle = 135 , end angle = 45 ] ;

\fill[ White ]
    ( 6.5 - 0.1926 , -0.0329 ) arc[ x radius = 0.2723 , y radius = 0.2 , start angle = 135 , end angle = 45 ]
    arc[ x radius = 0.35355 , y radius = 0.25 , start angle = 303 , end angle = 237 ]
    ;

\draw[ thick ]
    ( 6.5 - 0.25 , 0 ) arc[ x radius = 0.35355 , y radius = 0.25 , start angle = 225 , end angle = 315 ]
    ( 6.5 - 0.1926 , -0.0329 ) arc[ x radius = 0.2723 , y radius = 0.2 , start angle = 135 , end angle = 45 ] ;

\draw ( 1 , 0 ) node {$ \bullet $} node[ anchor = north ] {$ s $}
    ( 2.75 , -2 / 3 ) node {$ S $} ;

\end{tikzpicture}
\caption{Generators of the surface group $ \pi_{ 1 } \of{ S , s } $} \label{fig:generatorsWithPants}
\end{figure}

A closed regular neighborhood $ P $ of $ \overline{ \alpha } \cup \overline{ \beta } $ is homeomorphic to a pair of pants, with boundary components represented by $ \alpha $, $ \beta $, and exactly one of $ \alpha \beta $ and $ \alpha^{ -1 } \beta $ (which product represents this boundary curve depends on a choice of orientation for $ \alpha $ and $ \beta $). In particular, the boundary components of $ P $ are distinct. Fix integers $ j , k \in \Z $. If $ \alpha^{ j } \beta^{ k } $ is simple, then by \autoref{lem:simpleLengthTwins} it has no length twins. Thus we may assume that $ \alpha^{ j } \beta^{ k } $ is not simple.

Suppose that $ \alpha^{ j } \beta^{ k } $ and a curve $ \delta $ are length-commensurable. \autoref{thm:equivalence} then implies that $ \delta $ does not intersect any of the boundary curves of $ P $, so that either $ \delta \subseteq P $ or $ \delta \subseteq S \setminus P $.

If $ g = 2 $, then $ P $ and $ S \setminus P $ are pairs of pants interchanged by the hyperelliptic involution. Since the hyperelliptic involution is an isometry for any hyperbolic metric, it suffices in this case to show that some non-zero powers of $ \alpha^{ j } \beta^{ k } $ and $ \delta $ are conjugate if $ \delta \subseteq P $.

On the other hand, if $ g \geq 3 $, then $ S \setminus P $ is a surface of genus at least $ 1 $ with $ 3 $ boundary components. In particular, if $ \delta \subseteq S \setminus P $, then it intersects a simple curve contained entirely within $ S \setminus P $. By \autoref{thm:equivalence}, this is inconsistent with $ \alpha^{ j } \beta^{ k } $ and $ \delta $ being length-commensurable. Thus in this case, $ \delta \subseteq P $ so it also suffices in this case to show that some non-zero powers of $ \alpha^{ j } \beta^{ k } $ and $ \delta $ are conjugate if $ \delta \subseteq P $.

To that end, we will assume that $ \delta \subseteq P $, so that $ \delta \in \pi_{ 1 } \of{ P , s } $. Let $ \varphi \colon \pi_{ 1 } \of{ S , s } \onto F_{ 2 } $ be the surjective homomorphism onto the free group $ F_{ 2 } = \gen{ a , b \mid } $ sending $ \alpha $ to $ a $ and $ \beta $ to $ b $. $ \varphi \of!\big!{ \alpha^{ j } \beta^{ k } } = a^{ j } b^{ k } $ and $ \varphi \of{ \delta } $ are trace-commensurable, and so by \autoref{lem:extendedHorowitz}, some non-zero powers of these elements are conjugate in $ F_{ 2 } $. However, the restriction of $ \varphi $ to $ \pi_{ 1 } \of{ P , s } $ is injective, and some non-zero powers of $ \alpha^{ j } \beta^{ k } $ and $ \delta $ are conjugate in $ \pi_{ 1 } \of{ P , s } $. \qedhere

\end{proof}

\subsection{Regular covers of low degree} \label{sec:simplyGenerated}

Fix a closed, orientable surface $ S $ of negative Euler characteristic, and consider for each choice of basepoint $ s \in S $ the set
\[
\mathcal{ S }_{ \N } = \set{ \gamma^{ m } \in \pi_{ 1 } \of{ S , s } : \textrm{$ \gamma $ is simple and $ m \neq 0 $} },
\]
of powers of simple elements. We will call a finite-degree, regular cover $ p \colon X \to S $ \emph{simply generated} if $ \pi_{ 1 } \of{ X , x } = \gen{ p_{ * }^{ -1 } \of{ \mathcal{ S }_{ \N } } } $ for some (equivalently every) choice of basepoint $ x \in X $. Note that finite-degree abelian covers, which we considered in \cite[Theorem~1.2]{ALLX20}, are simply generated, but the pathological regular covers of \cite{KS16} and \cite{MP19} are not.

In this section, we aim to prove \autoref{thm:simplyGeneratedCovers}, which we restate below for the reader's convenience.

\simplyGeneratedCovers*

The proof will rely on a strengthening of \cite[Theorem~1.2]{ALLX20}. For each simply generated, finite-degree, regular cover $ p \colon X \to S $ and each simple element $ \gamma \in \pi_{ 1 } \of{ S , s } $, we will denote
\[
n_{ p } \of{ \gamma } = \ind{ \gen{ \gamma } }{ \gen{ \gamma } \cap p_{ * } \of{ \pi_{ 1 } \of{ X , x } } } .
\]
That is, $ n_{ p } \of{ \gamma } $ is the smallest positive integer $ m $ such that $ \gamma^{ m } \in p_{ * } \of{ \pi_{ 1 } \of{ X , x } } $, or equivalently the smallest degree of an elevation of the curve $ \gamma $ represents along $ p $ to $ X $. Since $ p $ is regular, $ n_{ p } \of{ \gamma } $ is also the degree of any elevation of this curve.

\begin{lemma*} \label{lem:simplyGenElevations}

Two simply generated, finite-degree, regular covers $ p \colon X \to S $ and $ q \colon Y \to S $ are isomorphic if and only if for some (equivalently, any) choice of basepoints $ x \in X $, $ y \in Y $, and $ s \in S $ with $ p \of{ x } = q \of{ y } = s $, $ n_{ p } \of{ \gamma } = n_{ q } \of{ \gamma } $ for all simple elements $ \gamma \in \pi_{ 1 } \of{ S , s } $.

\begin{proof}

Suppose that $ n_{ p } \of{ \gamma } = n_{ q } \of{ \gamma } $ for all simple elements $ \gamma \in \pi_{ 1 } \of{ S , s } $. We will first show that $ p_{ * } \of{ \pi_{ 1 } \of{ X , x } } \cap \mathcal{ S }_{ \N }  = q_{ * } \of{ \pi_{ 1 } \of{ Y , y } } \cap \mathcal{ S }_{ \N }  $. To that end, let $ \eta \in p_{ * } \of{ \pi_{ 1 } \of{ X , x } } \cap \mathcal{ S }_{ \N }  $. Then $ \eta = p_{ * } \of{ \alpha } = \gamma^{ m } $ for some element $ \alpha \in \pi_{ 1 } \of{ X , x } $, some simple element $ \gamma \in \pi_{ 1 } \of{ S , s } $, and a non-zero integer $ m \neq 0 $. Since $ \alpha $ represents an elevation of degree $ \abs{ m } $ of the curve represented by $ \gamma $ along $ p $ to $ X $, then
\[
n_{ q } \of{ \gamma } = n_{ p } \of{ \gamma } = \abs{ m }. 
\]
So there is an element $ \beta \in \pi_{ 1 } \of{ Y , y } $ such that $ q_{ * } \of{ \beta } = \gamma^{ m } = \eta $. Thus $ p_{ * } \of{ \pi_{ 1 } \of{ X , x } } \cap \mathcal{ S }_{ \N }  \subseteq q_{ * } \of{ \pi_{ 1 } \of{ Y , y } } \cap \mathcal{ S }_{ \N }  $. By symmetry, the other containment also holds. This implies that
\begin{align*}
p_{ * } \of{ \pi_{ 1 } \of{ X , x } } & = p_{ * } \of{ \gen{ p_{ * }^{ -1 } \of{\mathcal{ S }_{ \N }  } } } = \gen{ p_{ * } \of{ \pi_{ 1 } \of{ X , x } } \cap \mathcal{ S }_{ \N }  } \\
& = \gen{ q_{ * } \of{ \pi_{ 1 } \of{ Y , y } } \cap \mathcal{ P }_{ s } } = q_{ * } \of{ \gen{ q_{ * }^{ -1 } \of{ \mathcal{ S }_{ \N } } } } = q_{ * } \of{ \pi_{ 1 } \of{ Y , y } } ,
\end{align*}
and so $ p \colon X \to S $ and $ q \colon Y \to S $ are isomorphic as covers of $ S $. \qedhere

\end{proof}

\end{lemma*}

Armed with this characterization of isomorphism for simply generated regular covers, we are now ready to prove \autoref{thm:simplyGeneratedCovers}.

\begin{proof}[Proof of \autoref{thm:simplyGeneratedCovers}]

Suppose that $ p $ and $ q $ are not isomorphic. Then \autoref{lem:simplyGenElevations} implies that there are basepoints $ x \in X $, $ y \in Y $, and $ s \in S $ such that $ p \of{ x } = q \of{ y } = s $ and a simple element $ \gamma \in \pi_{ 1 } \of{ S , s } $ such that $ n_{ p } \of{ \gamma } \neq n_{ q } \of{ \gamma } $. Note that $ n_{ p } \of{ \gamma } $ and $ n_{ q } \of{ \gamma } $ denote the degree of any elevation of $ \gamma $ along $ p $ and $ q $, respectively; that all such elevations have the same degree follows from the regularity of $ p $ and $ q $.

Thus $ \gamma $ represents a simple curve on $ S $ with a simple elevation of degree $ n_{ p } \of{ \gamma } $ along $ p $ to $ X $ but no simple elevation of degree $ n_{ p } \of{ \gamma } $ along $ q $ to $ Y $. \autoref{lem:simpleLengthTwins} and \autoref{prop:lengthSpectra} imply that $ X $ and $ Y $ are generically not simple length isospectral over $ S $. \qedhere

\end{proof}

\subsection{Examples of generically simple length non-isospectral covers} \label{sec:examples}

We now provide a plethora of examples of pairs of finite-degree covers of a surface $ S $ of negative Euler characteristic which are length isospectral over $ S $ but generically not simple length isospectral over $ S $. The structure of our argument in all three examples follows the same basic structure:
\begin{enumerate}

\item[1.]

We construct a surjective homomorphism of $ \pi_{ 1 } \of{ S } $ onto a finite group $ G $ and identify a pair of finite-degree covers of $ S $ corresponding to a pair of almost conjugate subgroups of $ G $.

\item[2.]

We analyze the corresponding quotients of the Cayley graph of $ G $ to find a curve on $ S $ without length twins which has more simple elevations to one cover than it does to the other.

\item[3.]

We apply \autoref{thm:Sunada} and \autoref{prop:lengthSpectra} to conclude that the covers are length isospectral over $ S $ but generically not simple length isospectral over $ S $.

\end{enumerate}

In the following examples, we will follow the above outline, choosing the finite groups to be $ \of{ \Z / 8 \Z }^{ \times } \ltimes \of{ \Z / 8 \Z } $, $ \SL_{ 3 } \of{ \F_{ 2 } } $, and $ S_{ 6 } $, respectively.

\begin{example*} \label{example:semidirect}

Let $ G \coloneqq \of{ \Z / 8 \Z }^{ \times } \ltimes \of{ \Z / 8 \Z } $, which has the group law $ \of{ a , b } \cdot \of{ c , d } = \of{ a c , b + a d } $, and let $ S $ be a closed surface of genus $ g \geq 3 $. One can check (or refer to \cite[Example~11.2.2]{Bus10}) that the subgroups of $ G $ given by
\begin{align*}
A & \coloneqq \set{ \of{ 1 , 0 } , \of{ 3 , 0 } , \of{ 5 , 0 } , \of{ 7 , 0 } } , & B & \coloneqq \set{ \of{ 1 , 0 } , \of{ 3 , 4 } , \of{ 5 , 4 } , \of{ 7 , 0 } }
\end{align*}
are almost conjugate but not conjugate in $ G $, and moreover that $ G $ is generated by the three elements $ a_{ 1 } = \of{ 5 , 0 } $, $ a_{ 2 } = \of{ 3 , 0 } $, and $ a_{ 3 } = \of{ 1 , 1 } $.

Fix a basepoint $ s \in S $, and consider the standard presentation
\[
\pi_{ 1 } \of{ S , s } = \gen{ \alpha_{ 1 } , \dotsc , \alpha_{ g } , \beta_{ 1 } , \dotsc , \beta_{ g } \mid \comm{ \alpha_{ 1 } }{ \beta_{ 1 } } \dotsm \comm{ \alpha_{ g } }{ \beta_{ g } } = 1 },
\]
of the corresponding surface group $ \pi_{ 1 } \of{ S , s } $, the generators of which are drawn in \autoref{fig:generators}.

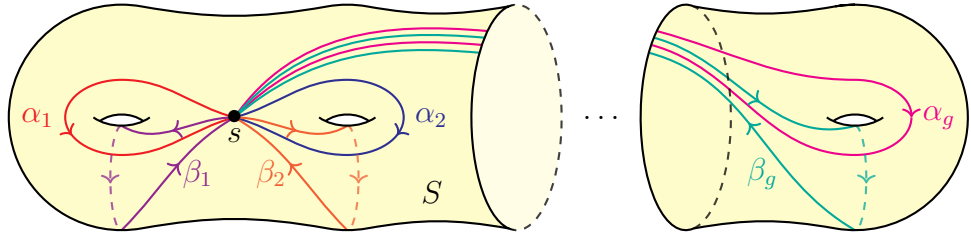
\begin{figure}[ht]
\begin{center}
\begin{tikzpicture}[ scale = 1.5 ]

\fill[ Yellow , opacity = 0.2 ]
    ( 0 , 1 ) arc[ radius = 1 , start angle = 90 , end angle = { 270 + 10 } ]
    to[ out = 10 , in = { 180 - 10 } ] ( 2 - 0.1736 , -0.9848 )
    arc[ radius = 1 , start angle = { -90 - 10 } , end angle = { -90 + 10 } ]
    to[ out = 10 , in = { 180 - 10 } ] ( 3.5 - 0.4 * 0.1736 , -0.9848 )
    arc[ x radius = { 0.4 * 1 } , y radius = 1 , start angle = { 270 - 10 } , end angle = { 90 + 10 } ]
    to[ out = { 180 + 10 } , in = -10 ] ( 2 + 0.1736 , 0.9848 )
    arc[ radius = 1 , start angle = { 90 - 10 } , end angle = { 90 + 10 } ]
    to[ out = { 180 + 10 } , in = -10 ] ( 0 + 0.1736 , 0.9848 )
    arc[ radius = 1 , start angle = { 90 - 10 } , end angle = 90 ]
    
    ( 6.5 , 1 ) arc[ radius = 1 , start angle = 90 , end angle = { 90 + 10 } ]
    to[ out = { 180 + 10 } , in = -10 ] ( 5 + 0.4 * 0.1736 , 0.9848 )
    arc[ x radius = 0.4 , y radius = 1 , start angle = { 90 - 10 } , end angle = { 270 + 10 } ]
    to[ out = 10 , in = { 180 - 10 } ] ( 6.5 - 0.1736 , -0.9848 )
    arc[ radius = 1 , start angle = { -90 - 10 } , end angle = { 90 + 10 } ] ;

\fill[ Yellow , opacity = 0.1 ]
    ( 3.5 - 0.4 * 0.1736 , -0.9848 ) arc[ x radius = 0.4 , y radius = 1 , start angle = { -90 - 10 } , end angle = { 270 - 10 } ]
    ;

\draw ( 4.25 , 0 ) node {$ \hdots $} ;

\draw[ thick , Red , ->- ]
    ( 1 , 0 )
    to[ out = 165 , in = 0 ] ( 0 , 1 / 3 )
    arc[ x radius = 0.5 , y radius = 1 / 3 , start angle = 90 , end angle = 180 ] node[ anchor = east ] {$ \alpha_{ 1 } $}
    arc[ x radius = 0.5 , y radius = 1 / 3 , start angle = 180 , end angle = 270 ]
    to[ out = 0 , in = 195 ] ( 1 , 0 ) ;

\draw[ thick , Blue , ->- ]
    ( 1 , 0 ) to[ out = 15 , in = 180 ] ( 2 , 1 / 3 )
    arc[ x radius = 1 / 2 , y radius = 1 / 3 , start angle = 90 , end angle = 0 ] node[ anchor = west ] {$ \alpha_{ 2 } $}
    arc[ x radius = 1 / 2 , y radius = 1 / 3 , start angle = 0 , end angle = -90 ]
    to[ out = 180 , in = -15 ] ( 1 , 0 ) ;

\draw[ thick , Plum , ->- ]
    ( 1 - 0.56 , -1 / 2 ) node[ anchor = west ] {$ \beta_{ 1 } $}
    ( 1 , 0 ) to[ out = 180 , in = -30 ] ( 0 , -0.075 ) ;
\draw[ thick , Plum , ->- ]
    ( 0 , -1 ) to[ out = 45 , in = 210 ] ( 1 , 0 ) ;

\draw[ thick , dashed , Plum , opacity = 0.75 , ->- ]
    ( 0 , -0.075 ) arc[ x radius = 0.1 , y radius = { ( 1 - 0.075 ) / 2 } , start angle = 90 , end angle = 270 ] ;

\draw[ thick , RedOrange , ->- ]
    ( 1 + 0.56 , -1 / 2 ) node[ anchor = east ] {$ \beta_{ 2 } $}
    ( 1 , 0 ) to[ out = 0 , in = 210 ] ( 2 , -0.075 ) ;
\draw[ thick , RedOrange , ->- ]
    ( 2 , -1 ) to[ out = 135 , in = -30 ] ( 1 , 0 ) ;

\draw[ thick , dashed , RedOrange , opacity = 0.75 , ->- ]
    ( 2 , -0.075 ) arc[ x radius = 0.1 , y radius = { ( 1 - 0.075 ) / 2 } , start angle = 90 , end angle = -90 ] ;

\draw[ thick , Magenta , -> ]
    ( 1 , 0 ) to[ out = 45 , in = { 180 - 5 } ] ( 3.5 - 0.4 * 0.7660 , 0.6428 )
    ( 1 , 0 ) to[ out = 55 , in = { 180 - 5 } ] ( 3.5 - 0.4 * 0.6428 , 0.7660 )
    ( 5 - 0.4 * 0.6428 , 0.7660 ) to[ out = -15 , in = 180 ] ( 6.5 , 1 / 3 )
    arc[ x radius = 1 / 2 , y radius = 1 / 3 , start angle = 90 , end angle = 0 ] node[ anchor = west ] {$ \alpha_{ g } $} ;
\draw[ thick , Magenta ]
    ( 7 , 0 ) arc[ x radius = 1 / 2 , y radius = 1 / 3 , start angle = 360 , end angle = 270 ]
    to[ out = 180 , in = -15 ] ( 5 - 0.4 * 0.7660 , 0.6428 ) ;

\draw[ thick , JungleGreen ]
    ( 6.5 - 0.585 , -1 / 2 ) node[ anchor = east ] {$ \beta_{ g } $}
    ( 1 , 0 ) to[ out = 50 , in = { 180 - 5 } ] ( 3.5 - 0.4 * 0.7071 , 0.7071 )
    ( 1 , 0 ) to[ out = 40 , in = { 180 - 5 } ] ( 3.5 - 0.4 * 0.8192 , 0.5736 ) ;
\draw[ thick , JungleGreen , ->- ]
    ( 5 - 0.4 * 0.7071 , 0.7071 ) to[ out = -15 , in = { 180 + 15 } ] ( 6.5 , -0.075 ) ;
\draw[ thick , JungleGreen , ->- ]
    ( 6.5 , -1 ) to[ out = { 180 - 30 } , in = -15 ] ( 5 - 0.4 * 0.8192 , 0.5736 ) ;

\draw[ thick , dashed , JungleGreen , opacity = 0.75 , ->- ]
    ( 6.5 , -0.075 ) arc[ x radius = 0.1 , y radius = { ( 1 - 0.075 ) / 2 } , start angle = 90 , end angle = -90 ] ;

\draw[ thick ]
    ( 0 , 1 ) arc[ radius = 1 , start angle = 90 , end angle = { 270 + 10 } ]
    to[ out = 10 , in = { 180 - 10 } ] ( 2 - 0.1736 , -0.9848 )
    arc[ radius = 1 , start angle = { -90 - 10 } , end angle = { -90 + 10 } ]
    to[ out = 10 , in = { 180 - 10 } ] ( 3.5 - 0.4 * 0.1736 , -0.9848 )
    arc[ x radius = { 0.4 * 1 } , y radius = 1 , start angle = { 270 - 10 } , end angle = { 90 + 10 } ]
    to[ out = { 180 + 10 } , in = -10 ] ( 2 + 0.1736 , 0.9848 )
    arc[ radius = 1 , start angle = { 90 - 10 } , end angle = { 90 + 10 } ]
    to[ out = { 180 + 10 } , in = -10 ] ( 0 + 0.1736 , 0.9848 )
    arc[ radius = 1 , start angle = { 90 - 10 } , end angle = 90 ]
    
    ( 6.5 , 1 ) arc[ radius = 1 , start angle = 90 , end angle = { 90 + 10 } ]
    to[ out = { 180 + 10 } , in = -10 ] ( 5 + 0.4 * 0.1736 , 0.9848 )
    arc[ x radius = 0.4 , y radius = 1 , start angle = { 90 - 10 } , end angle = { 270 + 10 } ]
    to[ out = 10 , in = { 180 - 10 } ] ( 6.5 - 0.1736 , -0.9848 )
    arc[ radius = 1 , start angle = { -90 - 10 } , end angle = { 90 + 10 } ] ;

\draw[ thick , dashed , opacity = 0.75 ]
    ( 3.5 - 0.4 * 0.1736 , -0.9848 ) arc[ x radius = 0.4 , y radius = 1 , start angle = { -90 - 10 } , end angle = { 90 + 10 } ]
    ( 5 + 0.4 * 0.1736 , -0.9848 ) arc[ x radius = 0.4 , y radius = 1 , start angle = { -90 + 10 } , end angle = { 90 - 10 } ] ;

\fill[ White ]
    ( -0.1926 , -0.0329 ) arc[ x radius = 0.2723 , y radius = 0.2 , start angle = 135 , end angle = 45 ]
    arc[ x radius = 0.35355 , y radius = 0.25 , start angle = 303 , end angle = 237 ]
    ;

\draw[ thick ]
    ( -0.25 , 0 ) arc[ x radius = 0.35355 , y radius = 0.25 , start angle = 225 , end angle = 315 ]
    ( -0.1926 , -0.0329 ) arc[ x radius = 0.2723 , y radius = 0.2 , start angle = 135 , end angle = 45 ] ;

\fill[ White ]
    ( 2 - 0.1926 , -0.0329 ) arc[ x radius = 0.2723 , y radius = 0.2 , start angle = 135 , end angle = 45 ]
    arc[ x radius = 0.35355 , y radius = 0.25 , start angle = 303 , end angle = 237 ]
    ;

\draw[ thick ]
    ( 2 - 0.25 , 0 ) arc[ x radius = 0.35355 , y radius = 0.25 , start angle = 225 , end angle = 315 ]
    ( 2 - 0.1926 , -0.0329 ) arc[ x radius = 0.2723 , y radius = 0.2 , start angle = 135 , end angle = 45 ] ;

\fill[ White ]
    ( 6.5 - 0.1926 , -0.0329 ) arc[ x radius = 0.2723 , y radius = 0.2 , start angle = 135 , end angle = 45 ]
    arc[ x radius = 0.35355 , y radius = 0.25 , start angle = 303 , end angle = 237 ]
    ;

\draw[ thick ]
    ( 6.5 - 0.25 , 0 ) arc[ x radius = 0.35355 , y radius = 0.25 , start angle = 225 , end angle = 315 ]
    ( 6.5 - 0.1926 , -0.0329 ) arc[ x radius = 0.2723 , y radius = 0.2 , start angle = 135 , end angle = 45 ] ;

\draw ( 1 , 0 ) node {$ \bullet $} node[ anchor = north ] {$ s $}
    ( 2.75 , -2 / 3 ) node {$ S $} ;

\end{tikzpicture}
\caption{Generators of the surface group $ \pi_{ 1 } \of{ S , s } $} \label{fig:generators}
\end{center}
\end{figure}

Consider the surjective homomorphism $ \rho \colon \pi_{ 1 } \of{ S , s } \onto G $ defined by the conditions $ \rho \of{ \alpha_{ i } } = a_{ i } $ for $ i \in \set{ 1 , 2 , 3 } $ and $ \rho \of{ \alpha_{ 4 } } = \dotsb = \rho \of{ \alpha_{ g } } = \rho \of{ \beta_{ 1 } } = \dotsb = \rho \of{ \beta_{ g } } = 1 $, which is well-defined in light of the fact that
\[
\comm{ \rho \of{ \alpha_{ 1 } } }{ \rho \of{ \beta_{ 1 } } } \dotsm \comm{ \rho \of{ \alpha_{ g } } }{ \rho \of{ \beta_{ g } } } = \comm{ a_{ 1 } }{ 1 } \comm{ a_{ 2 } }{ 1 } \comm{ a_{ 3 } }{ 1 } \cdot 1 \dotsm 1 = 1 .
\]
As described in \autoref{sec:Sunada}, the homomorphism $ \rho \colon \pi_{ 1 } \of{ S , s } \onto G $ and the subgroups $ A $ and $ B $ give rise to a pair of based covers $ p \colon \of{ X , x } \to \of{ S , s } $, $ q \colon \of{ Y , y } \to \of{ S , s } $ corresponding to $ \rho^{ -1 } \of{ A } $ and $ \rho^{ -1 } \of{ B } $, respectively. Note that $ X $ has genus
\begin{align*}
\frac{ 2 - \chi \of{ X } }{ 2 } & = \frac{ 2 - \deg \of{ p } \cdot \chi \of{ S } }{ 2 } = \frac{ 2 - \ind{ \pi_{ 1 } \of{ S , s } }{ \rho^{ -1 } \of{ A } } \cdot \of{ 2 - 2 g } }{ 2 } \\
& = 1 + \ind{ G }{ A } \of{ g - 1 } = 8 g - 7 .
\end{align*}
Similarly, $ Y $ also has genus $ 8 g - 7 $.

Since $ A $ and $ B $ are almost conjugate, it follows from \autoref{thm:Sunada} that $ X $ and $ Y $ are length isospectral over $ S $; that is, $ \of{ X , p^{ * } \of{ \mu } } $ and $ \of{ Y , q^{ * } \of{ \mu } } $ are length isospectral for any hyperbolic metric $ \mu \in \Teich \of{ S } $. To see that these covers are generically not simple length isospectral (and hence generically not isometric), it suffices by \autoref{prop:lengthSpectra} to produce a closed curve $ \gamma $ on $ S $ which has no length twins and different numbers of degree $ 1 $ simple elevations to $ X $ and $ Y $. We will find such a curve $ \gamma $ by analyzing quotients of the Cayley graph $ \mathcal{ G } $ of $ G $ with respect to the above generating set $ \set{ a_{ 1 } , a_{ 2 } , a_{ 3 } } $.

Specifically, let $ \mathcal{ G } $ be the directed graph with vertex set $ G $ and edge set
\[
\set{ \of{ g , h } \in G^{ 2 } : h = g s \textrm{ for some } s \in \set{ a_{ 1 } , a_{ 2 } , a_{ 3 } } } .
\]
Associating an edge $ \of{ g , h } $ in $ \mathcal{ G } $ with its difference $ g^{ -1 } h \in \set{ a_{ 1 } , a_{ 2 } , a_{ 3 } } $ gives an edge coloring of $ \mathcal{ G } $. Note that the action of $ G $ on itself by left multiplication extends to an action on $ \mathcal{ G } $ which preserves edges and colors, so that we may form the quotient graphs $ A \setminus \mathcal{ G } $ and $ B \setminus \mathcal{ G } $. These quotient graphs are drawn in \autoref{fig:quotientGraphs1}. Note that each vertex of $ A \setminus \mathcal{ G } $ and $ B \setminus \mathcal{ G } $ has one incoming edge and one outgoing edge of each color. These graphs are drawn in \autoref{fig:quotientGraphs1} so that self-loops are omitted and pairs of edges with the same endpoints and color but with opposite orientations are drawn as a single undirected edge.

\begin{figure}[ht]
\centering
\begin{tikzpicture}[ scale = 1.9 ]

\draw[ thick , ->- ] ( 0 , 0 ) to ( 1 , 0 ) ;
\draw[ thick , ->- ] ( 1 , 0 ) to ( 2 , 0 ) ;
\draw[ thick , ->- ] ( 2 , 0 ) to ( 3 , 0 ) ;
\draw[ thick , ->- ] ( 3 , 0 ) to ( 4 , 0 ) ;
\draw[ thick , ->- ] ( 4 , 0 ) to ( 5 , 0 ) ;
\draw[ thick , ->- ] ( 5 , 0 ) to ( 6 , 0 ) ;
\draw[ thick , ->- ] ( 6 , 0 ) to ( 7 , 0 ) ;
\draw[ thick , ->- ]
    ( 6 , -1 ) to[ out = 180 , in = 0 ] ( 1 , -1 ) ;
\draw[ thick ]
    ( 7 , 0 ) arc[ radius = 1 , start angle = 360 , end angle = 270 ]
    ( 0 , 0 ) arc[ radius = 1 , start angle = 180 , end angle = 270 ] ;

\draw[ thick , dashed ]
    ( 1 , 0 ) arc[ radius = 1 / 3 , start angle = 180 , end angle = 270 ]
        to[ out = 0 , in = 180 ] ( 3 - 1 / 3 , -1 / 3 )
        arc[ radius = 1 / 3 , start angle = 270 , end angle = 360 ]
    ( 2 , 0 ) arc[ radius = 1 / 3 , start angle = 180 , end angle = 90 ]
        to[ out = 0 , in = 180 ] ( 6 - 1 / 3 , 1 / 3 )
        arc[ radius = 1 / 3 , start angle = 90 , end angle = 0 ]
    ( 5 , 0 ) arc[ radius = 1 / 3 , start angle = 180 , end angle = 270 ]
        to[ out = 0 , in = 180 ] ( 7 - 1 / 3 , -1 / 3 )
        arc[ radius = 1 / 3 , start angle = 270 , end angle = 360 ] ;

\draw[ thick , dotted ]
    ( 1 , 0 ) arc[ radius = 2 / 3 , start angle = 180 , end angle = 270 ]
        to[ out = 0 , in = 180 ] ( 5 - 2 / 3 , -2 / 3 )
        arc[ radius = 2 / 3 , start angle = 270 , end angle = 360 ]
    ( 3 , 0 ) arc[ radius = 2 / 3 , start angle = 180 , end angle = 90 ]
        to[ out = 0 , in = 180 ] ( 7 - 2 / 3 , 2 / 3 )
        arc[ radius = 2 / 3 , start angle = 90 , end angle = 0 ] ;

\draw
    ( 0 , 0 ) node {$ \bullet $} node[ anchor = south ] {$ A $}
    ( 1 , 0 ) node {$ \bullet $} node[ anchor = south ] {$ A a_{ 3 } $}
    ( 2 , 0 ) node {$ \bullet $} node[ anchor = north ] {$ A a_{ 3 }^{ 2 } $}
    ( 3 , 0 ) node {$ \bullet $} node[ anchor = north west ] {$ A a_{ 3 }^{ 3 } $}
    ( 4 , 0 ) node {$ \bullet $} node[ anchor = north ] {$ A a_{ 3 }^{ 4 } $}
    ( 5 , 0 ) node {$ \bullet $} node[ anchor = south ] {$ A a_{ 3 }^{ 5 } $}
    ( 6 , 0 ) node {$ \bullet $} node[ anchor = north ] {$ A a_{ 3 }^{ 6 } $}
    ( 7 , 0 ) node {$ \bullet $} node[ anchor = west ] {$ A a_{ 3 }^{ 7 } $} ;

\draw[ thick , dotted ]
    ( 0 , -2 ) arc[ radius = 2 / 3 , start angle = 180 , end angle = 90 ]
        to[ out = 0 , in = 180 ] ( 4 - 2 / 3 , 2 / 3 - 2 )
        arc[ radius = 2 / 3 , start angle = 90 , end angle = 0 ]
    ( 2 , -2 ) arc[ radius = 1 / 2 , start angle = 180 , end angle = 270 ]
        to[ out = 0 , in = 180 ] ( 6 - 1 / 2 , -2 - 1 / 2 )
        arc[ radius = 1 / 2 , start angle = 270 , end angle = 360 ] ;

\draw[ thick , dashed ]
    ( 0 , -2 ) arc[ radius = 1 / 3 , start angle = 180 , end angle = 90 ]
        to[ out = 0 , in = 180 ] ( 4 - 1 / 3 , 1 / 3 - 2 )
        arc[ radius = 1 / 3 , start angle = 90 , end angle = 0 ]
    ( 1 , -2 ) arc[ radius = 3 / 4 , start angle = 180 , end angle = 270 ]
        to[ out = 0 , in = 180 ] ( 7 - 3 / 4 , -2 - 3 / 4 )
        arc[ radius = 3 / 4 , start angle = 270 , end angle = 360 ]
    ( 3 , -2 ) arc[ radius = 1 / 4 , start angle = 180 , end angle = 270 ]
        to[ out = 0 , in = 180 ] ( 5 - 1 / 4 , -2 - 1 / 4 )
        arc[ radius = 1 / 4 , start angle = 270 , end angle = 360 ] ;

\draw[ thick , ->- ] ( 0 , -2 ) to ( 1 , -2 ) ;
\draw[ thick , ->- ] ( 1 , -2 ) to ( 2 , -2 ) ;
\draw[ thick , ->- ] ( 2 , -2 ) to ( 3 , -2 ) ;
\draw[ thick , ->- ] ( 3 , -2 ) to ( 4 , -2 ) ;
\draw[ thick , ->- ] ( 4 , -2 ) to ( 5 , -2 ) ;
\draw[ thick , ->- ] ( 5 , -2 ) to ( 6 , -2 ) ;
\draw[ thick , ->- ] ( 6 , -2 ) to ( 7 , -2 ) ;
\draw[ thick , ->- ]
    ( 6 , -3 ) to[ out = 180 , in = 0 ] ( 1 , -3 ) ;
\draw[ thick ]
    ( 7 , -2 ) arc[ radius = 1 , start angle = 360 , end angle = 270 ]
    ( 0 , -2 ) arc[ radius = 1 , start angle = 180 , end angle = 270 ] ;

\draw
    ( 0 , -2 ) node {$ \bullet $} node[ anchor = east ] {$ B $}
    ( 1 , -2 ) node {$ \bullet $} node[ anchor = south ] {$ B a_{ 3 } $}
    ( 2 , -2 ) node {$ \bullet $} node[ anchor = south ] {$ B a_{ 3 }^{ 2 } $}
    ( 3 , -2 ) node {$ \bullet $} node[ anchor = south ] {$ B a_{ 3 }^{ 3 } $}
    ( 4 , -2 ) node {$ \bullet $} node[ anchor = south west ] {$ B a_{ 3 }^{ 4 } $}
    ( 5 , -2 ) node {$ \bullet $} node[ anchor = south ] {$ B a_{ 3 }^{ 5 } $}
    ( 6 , -2 ) node {$ \bullet $} node[ anchor = south ] {$ B a_{ 3 }^{ 6 } $}
    ( 7 , -2 ) node {$ \bullet $} node[ anchor = south ] {$ B a_{ 3 }^{ 7 } $} ;

\end{tikzpicture}
\caption{The quotient graphs $ A \setminus \mathcal{ G } $ and $ B \setminus \mathcal{ G } $. Edges of color $ a_{ 1 } $, $ a_{ 2 } $, and $ a_{ 3 } $ are represented by dotted, dashed, and solid lines, respectively.} \label{fig:quotientGraphs1}
\end{figure}
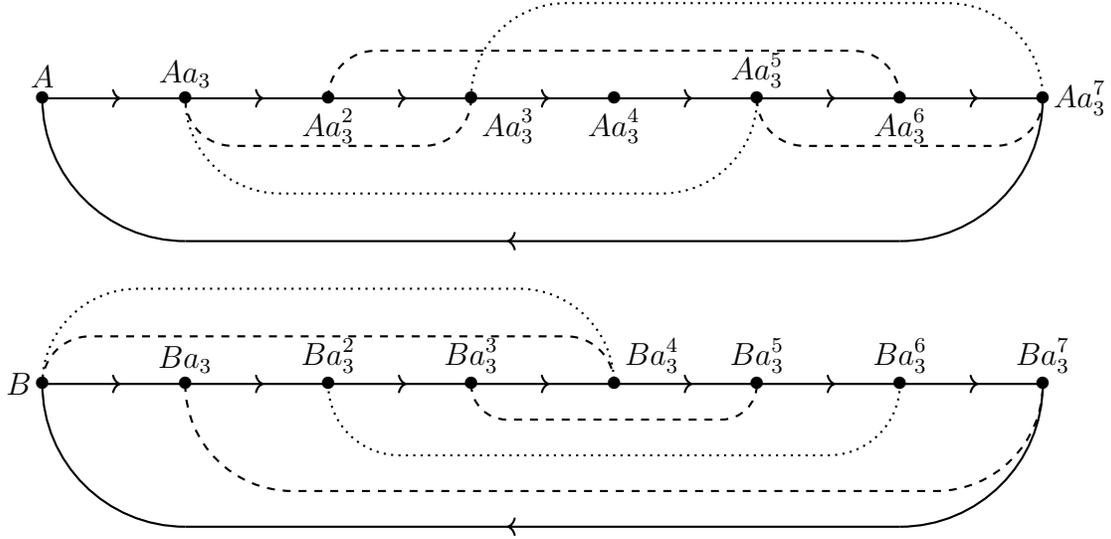

The vertex set of $ A \setminus \mathcal{ G } $ is in natural bijection with the fiber of $ p $ over the basepoint $ s \in S $, and for each index $ i \in \set{ 1 , 2 , 3 } $, the edges of color $ a_{ i } $ are in natural bijection with the path lifts of $ \alpha_{ i } $. Thus paths through $ A \setminus \mathcal{ G } $ correspond to path lifts of compositions of the curves $ \alpha_{ 1 } $, $ \alpha_{ 2 } $, and $ \alpha_{ 3 } $ and their inverses. Analogous statements hold for $ B \setminus \mathcal{ G } $.

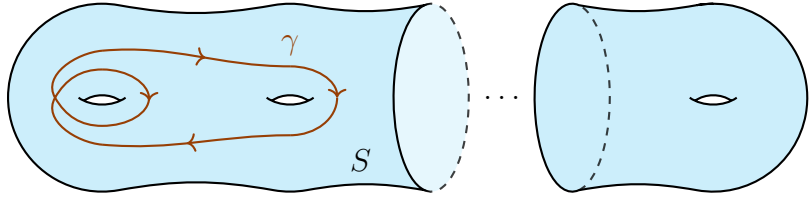
\begin{figure}[ht]
\begin{center}
\begin{tikzpicture}[ scale = 1.25 ]

\fill[ Cyan , opacity = 0.2 ]
    ( 0 , 1 ) arc[ radius = 1 , start angle = 90 , end angle = { 270 + 10 } ]
    to[ out = 10 , in = { 180 - 10 } ] ( 2 - 0.1736 , -0.9848 )
    arc[ radius = 1 , start angle = { -90 - 10 } , end angle = { -90 + 10 } ]
    to[ out = 10 , in = { 180 - 10 } ] ( 3.5 - 0.4 * 0.1736 , -0.9848 )
    arc[ x radius = { 0.4 * 1 } , y radius = 1 , start angle = { 270 - 10 } , end angle = { 90 + 10 } ]
    to[ out = { 180 + 10 } , in = -10 ] ( 2 + 0.1736 , 0.9848 )
    arc[ radius = 1 , start angle = { 90 - 10 } , end angle = { 90 + 10 } ]
    to[ out = { 180 + 10 } , in = -10 ] ( 0 + 0.1736 , 0.9848 )
    arc[ radius = 1 , start angle = { 90 - 10 } , end angle = 90 ]
    
    ( 6.5 , 1 ) arc[ radius = 1 , start angle = 90 , end angle = { 90 + 10 } ]
    to[ out = { 180 + 10 } , in = -10 ] ( 5 + 0.4 * 0.1736 , 0.9848 )
    arc[ x radius = 0.4 , y radius = 1 , start angle = { 90 - 10 } , end angle = { 270 + 10 } ]
    to[ out = 10 , in = { 180 - 10 } ] ( 6.5 - 0.1736 , -0.9848 )
    arc[ radius = 1 , start angle = { -90 - 10 } , end angle = { 90 + 10 } ] ;

\fill[ Cyan , opacity = 0.1 ]
    ( 3.5 - 0.4 * 0.1736 , -0.9848 ) arc[ x radius = 0.4 , y radius = 1 , start angle = { -90 - 10 } , end angle = { 270 - 10 } ]
    ;

\draw ( 4.25 , 0 ) node {$ \hdots $} ;

\draw[ RawSienna , thick ]
    ( -1 / 2 , 0 ) to[ out = 60 , in = 180 ] ( 0 , 0.3 )
    ( 1 / 2 , 0 ) arc[ x radius = 1 / 2 , y radius = 0.3 , start angle = 360 , end angle = 270 ]
    to[ out = 180 , in = -60 ] ( -1 / 2 , 0 )
    to[ out = 120 , in = 185 ] ( 0 , 1 / 2 )
    ( 2.5 , 0 ) arc[ x radius = 1 / 2 , y radius = 0.4 , start angle = 360 , end angle = 270 ]
    ( 0 , -1 / 2 ) to[ out = 175 , in = 240 ] ( -1 / 2 , 0 ) ;
\draw[ RawSienna , thick , -> ]
    ( 0 , 0.3 ) arc[ x radius = 1 / 2 , y radius = 1 / 3 , start angle = 90 , end angle = 0 ] ;
\draw[ thick , RawSienna , -> ]
    ( 2 , 1 / 3 ) arc[ x radius = 1 / 2 , y radius = 1 / 3 , start angle = 90 , end angle = 0 ] ;
\draw[ thick , RawSienna , ->- ]
    ( 0 , 1 / 2 ) to[ out = 5 , in = 180 ] ( 2 , 1 / 3 ) node[ anchor = south ] {$ \gamma $} ;
\draw[ thick , RawSienna , ->- ]
    ( 2 , -0.4 ) to[ out = 180 , in = -5 ] ( 0 , -1 / 2 ) ;

\draw[ thick ]
    ( 0 , 1 ) arc[ radius = 1 , start angle = 90 , end angle = { 270 + 10 } ]
    to[ out = 10 , in = { 180 - 10 } ] ( 2 - 0.1736 , -0.9848 )
    arc[ radius = 1 , start angle = { -90 - 10 } , end angle = { -90 + 10 } ]
    to[ out = 10 , in = { 180 - 10 } ] ( 3.5 - 0.4 * 0.1736 , -0.9848 )
    arc[ x radius = { 0.4 * 1 } , y radius = 1 , start angle = { 270 - 10 } , end angle = { 90 + 10 } ]
    to[ out = { 180 + 10 } , in = -10 ] ( 2 + 0.1736 , 0.9848 )
    arc[ radius = 1 , start angle = { 90 - 10 } , end angle = { 90 + 10 } ]
    to[ out = { 180 + 10 } , in = -10 ] ( 0 + 0.1736 , 0.9848 )
    arc[ radius = 1 , start angle = { 90 - 10 } , end angle = 90 ]
    
    ( 6.5 , 1 ) arc[ radius = 1 , start angle = 90 , end angle = { 90 + 10 } ]
    to[ out = { 180 + 10 } , in = -10 ] ( 5 + 0.4 * 0.1736 , 0.9848 )
    arc[ x radius = 0.4 , y radius = 1 , start angle = { 90 - 10 } , end angle = { 270 + 10 } ]
    to[ out = 10 , in = { 180 - 10 } ] ( 6.5 - 0.1736 , -0.9848 )
    arc[ radius = 1 , start angle = { -90 - 10 } , end angle = { 90 + 10 } ] ;

\draw[ thick , dashed , opacity = 0.75 ]
    ( 3.5 - 0.4 * 0.1736 , -0.9848 ) arc[ x radius = 0.4 , y radius = 1 , start angle = { -90 - 10 } , end angle = { 90 + 10 } ]
    ( 5 + 0.4 * 0.1736 , -0.9848 ) arc[ x radius = 0.4 , y radius = 1 , start angle = { -90 + 10 } , end angle = { 90 - 10 } ] ;

\fill[ White ]
    ( -0.1926 , -0.0329 ) arc[ x radius = 0.2723 , y radius = 0.2 , start angle = 135 , end angle = 45 ]
    arc[ x radius = 0.35355 , y radius = 0.25 , start angle = 303 , end angle = 237 ]
    ;

\draw[ thick ]
    ( -0.25 , 0 ) arc[ x radius = 0.35355 , y radius = 0.25 , start angle = 225 , end angle = 315 ]
    ( -0.1926 , -0.0329 ) arc[ x radius = 0.2723 , y radius = 0.2 , start angle = 135 , end angle = 45 ] ;

\fill[ White ]
    ( 2 - 0.1926 , -0.0329 ) arc[ x radius = 0.2723 , y radius = 0.2 , start angle = 135 , end angle = 45 ]
    arc[ x radius = 0.35355 , y radius = 0.25 , start angle = 303 , end angle = 237 ]
    ;

\draw[ thick ]
    ( 2 - 0.25 , 0 ) arc[ x radius = 0.35355 , y radius = 0.25 , start angle = 225 , end angle = 315 ]
    ( 2 - 0.1926 , -0.0329 ) arc[ x radius = 0.2723 , y radius = 0.2 , start angle = 135 , end angle = 45 ] ;

\fill[ White ]
    ( 6.5 - 0.1926 , -0.0329 ) arc[ x radius = 0.2723 , y radius = 0.2 , start angle = 135 , end angle = 45 ]
    arc[ x radius = 0.35355 , y radius = 0.25 , start angle = 303 , end angle = 237 ]
    ;

\draw[ thick ]
    ( 6.5 - 0.25 , 0 ) arc[ x radius = 0.35355 , y radius = 0.25 , start angle = 225 , end angle = 315 ]
    ( 6.5 - 0.1926 , -0.0329 ) arc[ x radius = 0.2723 , y radius = 0.2 , start angle = 135 , end angle = 45 ] ;

\draw ( 2.75 , -2 / 3 ) node {$ S $} ;

\end{tikzpicture}
\caption{The curve $ \gamma $ represented by $ \alpha_{ 1 }^{ -2 } \alpha_{ 2 } $.} \label{fig:curve1}
\end{center}
\end{figure}

Now let $ \gamma $ be the curve on $ S $ represented by the element $ \alpha_{ 1 }^{ -2 } \alpha_{ 2 } $. One may check by close inspection of the configuration of $ \alpha_{ 1 } $ and $ \alpha_{ 2 } $ in \autoref{fig:generators} that a closed path-lift of $ \alpha_{ 1 }^{ -2 } \alpha_{ 2 } $ based at a point $ x' $ in the fiber over $ s $ is simple precisely when the path-lift of $ \alpha_{ 1 } $ based at $ x' $ is not closed. Thus the two degree $ 1 $ elevations of $ \gamma $ along $ p $ to $ X $ (based at the points corresponding to the cosets $ A $ and $A a_{ 3 }^{ 4 } $) are not simple, but the two degree $ 1 $ elevations of $ \gamma $ along $ q $ to $ Y $ (based at the points corresponding to the cosets $ Ba_{ 3 }^{ 2 }  $ and $ Ba_{ 3 }^{ 6 }  $) are simple.

In summary, $ \gamma $ has no simple degree $ 1 $ elevations along $ p $ to $ X $, but it has two simple degree $ 1 $ elevations along $ q $ to $ Y $. $ \gamma $ has no length twins by \autoref{prop:extendedHorowitz}, so \autoref{prop:lengthSpectra} implies that $ X $ and $ Y $ are generically not simple length isospectral over $ S $.

\end{example*}

We remark that the above analysis in \autoref{example:semidirect} required the hypothesis that $ S $ has genus $ g \geq 3 $, since it used three freely related simple generators of $ \pi_{ 1 } \of{ S } $ in a configuration that isn't possible on a surface of genus $ 2 $. Moreover, our argument used this hypothesis in a more subtle way, since any non-simple curve of the type described in \autoref{prop:extendedHorowitz} has exactly one length twin on a closed surface of genus $ 2 $, which is its image under the hyperelliptic involution.

The first of these considerations does not apply to our analysis of \autoref{example:SpecialLinear} below, since $ \SL_{ 3 } \of{ \F_{ 2 } } $ is generated by two elements. However, in order to extend this example to the closed surface of genus $ 2 $, more care would be needed in dealing with any possible simple elevations of the length twin of the given curve. Although we do not perform such a computation in this article, Maungchang \cite[Theorem~1.1]{Mau13} completed this type of analysis for a pair of covers of the surface of genus $ 2 $ which are length isospectral but generically not simple length isospectral. Armed with this result and either of \autoref{example:semidirect} or \autoref{example:SpecialLinear}, we may conclude that we have proven \autoref{thm:examples}.

\begin{example*}
\label{example:SpecialLinear} 

Let $ G \coloneqq \SL_{ 3 } \of{ \F_{ 2 } } $ and let $ S $ be a closed surface of genus $ g \geq 3 $. One can check (or refer to \cite[Example~11.2.5, Lemma~11.2.6]{Bus10}) that the opposite parabolic subgroups of $ G $ given by
\begin{align*}
A & \coloneqq \set{ \vect{ a } \in \SL_{ 3 } \of{ \F_{ 2 } } : a_{ 2 , 1 } = a_{ 3 , 1 } = 0 } , & B & \coloneqq \trans{ A } = \set{ \vect{ b } \in \SL_{ 3 } \of{ \F_{ 2 } } : b_{ 1 , 2 } = b_{ 1 , 3 } = 0 }
\end{align*}
are almost conjugate but not conjugate in $ G $, and moreover that $ G $ is generated by the two elements
\begin{align*}
\vect{ a }_{ 1 } & \coloneqq \begin{pmatrix} 0 & 1 & 1 \\ 0 & 1 & 0 \\ 1 & 0 & 0 \end{pmatrix} , & \vect{ a }_{ 2 } & \coloneqq \begin{pmatrix} 1 & 0 & 0 \\ 0 & 0 & 1 \\ 0 & 1 & 1 \end{pmatrix} .
\end{align*}
As in \autoref{example:semidirect}, fix a basepoint $ s \in S $, and consider the standard presentation
\[
\pi_{ 1 } \of{ S , s } = \gen{ \alpha_{ 1 } , \dotsc , \alpha_{ g } , \beta_{ 1 } , \dotsc , \beta_{ g } \mid \comm{ \alpha_{ 1 } }{ \beta_{ 1 } } \dotsm \comm{ \alpha_{ g } }{ \beta_{ g } } = 1 },
\]
of the corresponding surface group $ \pi_{ 1 } \of{ S , s } $, the generators of which are drawn in \autoref{fig:generators}.

Consider the surjective homomorphism $ \rho \colon \pi_{ 1 } \of{ S , s } \onto G $ defined by the conditions $ \rho \of{ \alpha_{ i } } = \vect{ a }_{ i } $ for $ i \in \set{ 1 , 2 } $ and $ \rho \of{ \alpha_{ 3 } } = \dotsb = \rho \of{ \alpha_{ g } } = \rho \of{ \beta_{ 1 } } = \dotsb = \rho \of{ \beta_{ g } } = \vect{ I } $, which is well-defined in light of the fact that
\[
\comm{ \rho \of{ \alpha_{ 1 } } }{ \rho \of{ \beta_{ 1 } } } \dotsm \comm{ \rho \of{ \alpha_{ g } } }{ \rho \of{ \beta_{ g } } } = \comm{ \vect{ a }_{ 1 } }{ \vect{ I } } \comm{ \vect{ a }_{ 2 } }{ \vect{ I } } \comm{ \vect{ I } }{ \vect{ I } } \dotsm \comm{ \vect{ I } }{ \vect{ I } } = \vect{ I } .
\]
Note that the commutators of matrices in the above equation are the group-theoretic commutators of $ G = \SL_{ 3 } \of{ \F_{ 2 } } $, not the ring-theoretic commutators of the ring of $ 3 \times 3 $ matrices with entries in $ \F_{ 2 } $.

As described in \autoref{sec:Sunada}, the homomorphism $ \rho \colon \pi_{ 1 } \of{ S , s } \onto G $ and the subgroups $ A $ and $ B $ give rise to a pair of based covers $ p \colon \of{ X , x } \to \of{ S , s } $ and $ q \colon \of{ Y , y } \to \of{ S , s } $ corresponding to $ \rho^{ -1 } \of{ A } $ and $ \rho^{ -1 } \of{ B } $, respectively. Note that $ X $ has genus
\begin{align*}
\frac{ 2 - \chi \of{ X } }{ 2 } & = \frac{ 2 - \deg \of{ p } \cdot \chi \of{ S } }{ 2 } = \frac{ 2 - \ind{ \pi_{ 1 } \of{ S , s } }{ \rho^{ -1 } \of{ A } } \of{ 2 - 2 g } }{ 2 } \\
& = 1 + \ind{ G }{ A } \of{ g - 1 } = 7 g - 6 .
\end{align*}
Similarly, $ Y $ also has genus $ 7 g - 6 $.

Since $ A $ and $ B $ are almost conjugate in $ G $, it follows from \autoref{thm:Sunada} that $ X $ and $ Y $ are length isospectral over $ S $. As in \autoref{example:semidirect}, we proceed by analyzing quotients of the Cayley graph $ \mathcal{ G } $ of $ G $ with respect to the above generating set $ \set{ \vect{ a }_{ 1 } , \vect{ a }_{ 2 } } $.

Specifically, let $ \mathcal{ G } $ be the directed graph with vertex set $ G $ and edge set
\[
\set{ \of{ g , h } \in G^{ 2 } : h = g s \textrm{ for some } s \in \set{ \vect{ a }_{ 1 } , \vect{ a }_{ 2 }  } } .
\]
Associating an edge $ \of{ g , h } $ in $ \mathcal{ G } $ with its difference $ g^{ -1 } h \in \set{ \vect{ a }_{ 1 } , \vect{ a }_{ 2 } } $ gives an edge coloring of $ \mathcal{ G } $. Note that the action of $ G $ on itself by left multiplication extends to an action on $ \mathcal{ G } $ which preserves edges and colors, so that we may form the quotient graphs $ A \setminus \mathcal{ G } $ and $ B \setminus \mathcal{ G } $. Note that each vertex of $ A \setminus \mathcal{ G } $ and $ B \setminus \mathcal{ G } $ has one incoming edge and one outgoing edge of each color. These graphs are drawn in \autoref{fig:quotientGraphs2}, in which self-loops are not drawn and pairs of edges with the same endpoints and color but opposite orientations are drawn as a single un-directed edge.

\begin{figure}[ht]
\centering
\begin{tikzpicture}[ scale = 1.5 ]

\draw[ thick ] ( 0 , 0 ) to ( 1 , 0 ) ;
\draw[ thick , ->- ] ( 1.8660 , -1 / 2 ) to[ out = 60 , in = 300 ] ( 1.8660 , 1 / 2 ) ;
\draw[ thick , ->- ] ( 1.8660 , 1 / 2 ) to ( 2.8660 , 1 / 2 ) ;
\draw[ thick , ->- ] ( 2.8660 , 1 / 2 ) to[ out = 240 , in = 120 ] ( 2.8660 , -1 / 2 ) ;
\draw[ thick , ->- ] ( 2.8660 , -1 / 2 ) to ( 1.8660 , -1 / 2 ) ;

\draw[ thick , dashed , ->- ] ( 1 , 0 ) to ( 1.8660 , -1 / 2 ) ;
\draw[ thick , dashed , ->- ] ( 1.8660 , -1 / 2 ) to[ out = 120 , in = 240 ] ( 1.8660 , 1 / 2 ) ;
\draw[ thick , dashed , ->- ] ( 1.8660 , 1 / 2 ) to ( 1 , 0 ) ;
\draw[ thick , dashed , ->- ] ( 2.8660 , -1 / 2 ) to[ out = 60 , in = 300 ] ( 2.8660 , 1 / 2 ) ;
\draw[ thick , dashed , ->- ] ( 2.8660 , 1 / 2 ) to ( 3.7321 , 0 ) ;
\draw[ thick , dashed , ->- ] ( 3.7321 , 0 ) to ( 2.8660 , -1 / 2 ) ;

\draw
    ( 0 , 0 ) node {$ \bullet $} node[ anchor = south ] {$ A $}
    ( 1 , 0 ) node {$ \bullet $} node[ anchor = south ] {$ A \vect{ b }^{ 3 } $}
    ( { 1 + sqrt( 3 ) / 2 } , 1 / 2 ) node {$ \bullet $} node[ anchor = south ] {$ A \vect{ b }^{ 4 } $}
    ( { 1 + sqrt( 3 ) / 2 } , -1 / 2 ) node {$ \bullet $} node[ anchor = north ] {$ A \vect{ b }^{ 6 } $}
    ( { 2 + sqrt( 3 ) / 2 } , 1 / 2 ) node {$ \bullet $} node[ anchor = south ] {$ A \vect{ b }^{ 5 } $}
    ( { 2 + sqrt( 3 ) / 2 } , -1 / 2 ) node {$ \bullet $} node[ anchor = north ] {$ A \vect{ b }^{ 2 } $}
    ( { 2 + sqrt( 3 ) } , 0 ) node {$ \bullet $} node[ anchor = west ] {$ A \vect{ b } $} ;

\draw[ thick , ->- ]
    ( 0 , -2 ) to[ out = 45 , in = 165 ] ( 1.8660 , -2 + 1 / 2 ) ;
\draw[ thick , ->- ] ( 1.8660 , -2 + 1 / 2 ) to[ out = 0 , in = 90 ] ( 2.7321 , -2 ) ;
\draw[ thick , ->- ] ( 2.7321 , -2 ) to[ out = 270 , in = 0 ] ( 1.8660 , -2 - 1 / 2 ) ;
\draw[ thick , ->- ] ( 1.8660 , -2 - 1 / 2 ) to[ out = 195 , in = 315 ] ( 0 , -2 ) ;
\draw[ thick ] ( 3.5981 , -2 - 1 / 2 ) to[ out = 60 , in = 300 ] ( 3.5981 , -2 + 1 / 2 ) ;

\draw[ thick , dashed , ->- ] ( 1.8660 , -2 + 1 / 2 ) to ( 1 , -2 ) ;
\draw[ thick , dashed , ->- ] ( 1 , -2 ) to ( 1.8660 , -2 - 1 / 2 ) ;
\draw[ thick , dashed , ->- ] ( 1.8660 , -2 - 1 / 2 ) to ( 1.8660 , -2 + 1 / 2 ) ;
\draw[ thick , dashed , ->- ] ( 2.7321 , -2 ) to ( 3.5981 , -2 + 1 / 2 ) ;
\draw[ thick , dashed , ->- ] ( 3.5981 , -2 + 1 / 2 ) to[ out = 240 , in = 120 ] ( 3.5981 , -2 - 1 / 2 ) ;
\draw[ thick , dashed , ->- ] ( 3.5981 , -2 - 1 / 2 ) to ( 2.7321 , -2 ) ;

\draw
    ( 0 , -2 ) node {$ \bullet $} node[ anchor = east ] {$ B $}
    ( 1 , -2 ) node {$ \bullet $} node[ anchor = east ] {$ B \vect{ b }^{ 4 } $}
    ( { 1 + sqrt( 3 ) / 2 } , -2 + 1 / 2 ) node {$ \bullet $} node[ anchor = south ] {$ B \vect{ b } $}
    ( { 1 + sqrt( 3 ) / 2 } , -2 - 1 / 2 ) node {$ \bullet $} node[ anchor = north ] {$ B \vect{ b }^{ 3 } $}
    ( { 1 + sqrt( 3 ) } , -2 ) node {$ \bullet $} node[ anchor = east ] {$ B \vect{ b }^{ 6 } $}
    ( { 1 + 3 * sqrt( 3 ) / 2 } , -2 + 1 / 2 ) node {$ \bullet $} node[ anchor = south ] {$ B \vect{ b }^{ 2 } $}
    ( { 1 + 3 * sqrt( 3 ) / 2 } , -2 - 1 / 2 ) node {$ \bullet $} node[ anchor = north ] {$ B \vect{ b }^{ 5 } $} ;

\end{tikzpicture}
\caption{The quotient graphs $ A \setminus \mathcal{ G } $ and $ B \setminus \mathcal{ G } $. The cosets are enumerated by powers of the commutator $ \vect{ b } \coloneqq \comm{ \vect{ a }_{ 1 } }{ \vect{ a }_{ 2 } } $. Edges of color $ \vect{ a }_{ 1 } $ and $ \vect{ a }_{ 2 } $ are represented by solid and dashed lines, respectively.} \label{fig:quotientGraphs2}
\end{figure}
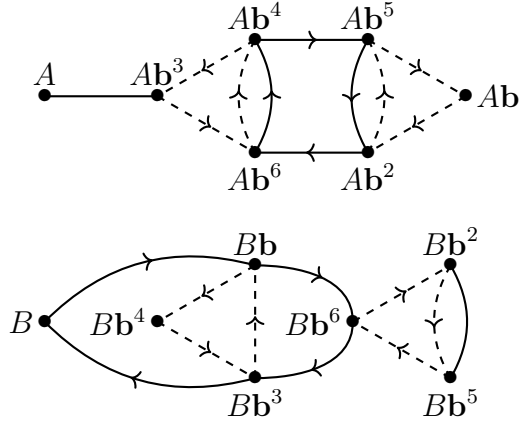

The vertex set of $ A \setminus \mathcal{ G } $ is in natural bijection with the fiber of $ p $ over the basepoint $ s \in S $, and for each index $ i \in \set{ 1 , 2 } $, the edges of color $ \vect{ a }_{ i } $ are in natural bijection with the path lifts of $ \alpha_{ i } $. Thus paths through $ A \setminus \mathcal{ G } $ correspond to path lifts of compositions of the curves $ \alpha_{ 1 } $ and $ \alpha_{ 2 } $ and their inverses. Analogous statements hold for $ B \setminus \mathcal{ G } $.

Now let $ \gamma $ be the curve on $ S $ represented by the element $ \alpha_{ 1 }^{ 2 } \alpha_{ 2 }^{ 2 } $. One may check by close inspection of the configuration of $ \alpha_{ 1 } $ and $ \alpha_{ 2 } $ in \autoref{fig:generators} that a closed path-lift of $ \alpha_{ 1 }^{ 2 } \alpha_{ 2 }^{ 2 } $ based at a point $ x' $ in the fiber over $ s $ is simple precisely when
\begin{enumerate*}

\item[(i)]

the path-lift of $ \alpha_{ 1 } $ based at $ x' $ is not closed; and

\item[(ii)]

the path-lift of $ \alpha_{ 2 } $ based at the endpoint of the path-lift of $ \alpha_{ 1 }^{ 2 } $ based at $ x' $ is not closed.

\end{enumerate*}
Thus the only degree $ 1 $ elevation of $ \gamma $ along $ p $ to $ X $ (based at the point corresponding to the coset $ A $) is not simple, but the only degree $ 1 $ elevation of $ \gamma $ along $ q $ to $ Y $ (based at the point corresponding to the coset $ B \vect{ b }^{ 3 } $) is simple.

\begin{figure}[t]
\begin{center}
\begin{tikzpicture}[ scale = 1.25 ]

\fill[ Magenta , opacity = 0.2 ]
    ( 0 , 1 ) arc[ radius = 1 , start angle = 90 , end angle = { 270 + 10 } ]
    to[ out = 10 , in = { 180 - 10 } ] ( 2 - 0.1736 , -0.9848 )
    arc[ radius = 1 , start angle = { -90 - 10 } , end angle = { -90 + 10 } ]
    to[ out = 10 , in = { 180 - 10 } ] ( 3.5 - 0.4 * 0.1736 , -0.9848 )
    arc[ x radius = { 0.4 * 1 } , y radius = 1 , start angle = { 270 - 10 } , end angle = { 90 + 10 } ]
    to[ out = { 180 + 10 } , in = -10 ] ( 2 + 0.1736 , 0.9848 )
    arc[ radius = 1 , start angle = { 90 - 10 } , end angle = { 90 + 10 } ]
    to[ out = { 180 + 10 } , in = -10 ] ( 0 + 0.1736 , 0.9848 )
    arc[ radius = 1 , start angle = { 90 - 10 } , end angle = 90 ]
    
    ( 6.5 , 1 ) arc[ radius = 1 , start angle = 90 , end angle = { 90 + 10 } ]
    to[ out = { 180 + 10 } , in = -10 ] ( 5 + 0.4 * 0.1736 , 0.9848 )
    arc[ x radius = 0.4 , y radius = 1 , start angle = { 90 - 10 } , end angle = { 270 + 10 } ]
    to[ out = 10 , in = { 180 - 10 } ] ( 6.5 - 0.1736 , -0.9848 )
    arc[ radius = 1 , start angle = { -90 - 10 } , end angle = { 90 + 10 } ] ;

\fill[ Magenta , opacity = 0.1 ]
    ( 3.5 - 0.4 * 0.1736 , -0.9848 ) arc[ x radius = 0.4 , y radius = 1 , start angle = { -90 - 10 } , end angle = { 270 - 10 } ]
    ;

\draw ( 4.25 , 0 ) node {$ \hdots $} ;

\draw[ thick , RawSienna , -> ]
    ( 1 , 0 ) to[ out = 135 , in = 0 ] ( 0 , 1 / 2 ) ;
\draw[ thick , RawSienna ]
    ( 0 , 1 / 2 ) to[ out = 180 , in = 120 ] ( -1 / 2 , 0 )
    ( 2 / 5 , 0 ) to[ out = 90 , in = 60 ] ( -1 / 2 , 0 )
    ( 0 , -1 / 2 ) to[ out = 0 , in = 225 ] ( 1 , 0 )
    ( 2 , 1 / 2 ) to[ out = 0 , in = 60 ] ( 2 + 1 / 2 , 0 )
    ( 2 - 2 / 5 , 0 ) to[ out = 90 , in = 120 ] ( 2 + 1 / 2 , 0 )
    ( 2 , -1 / 2 ) to[ out = 180 , in = 315 ] ( 1 , 0 ) ;
\draw[ thick , RawSienna , -> ]
    ( -1 / 2 , 0 ) to[ out = 240 , in = 180 ] ( 0 , -1 / 2 ) ;
\draw[ thick , RawSienna , -> ]
    ( -1 / 2 , 0 ) to[ out = -60 , in = 270 ] ( 2 / 5 , 0 ) ;
\draw[ thick , RawSienna , -> ]
    ( 1 , 0 ) to[ out = 45 , in = 180 ] ( 2 , 1 / 2 ) node[ anchor = south ] {$ \gamma $} ;
\draw[ thick , RawSienna , -> ]
    ( 2 + 1 / 2 , 0 ) to[ out = 240 , in = 270 ] ( 2 - 2 / 5 , 0 ) ;
\draw[ thick , RawSienna , -> ]
    ( 2 + 1 / 2 , 0 ) to[ out = 300 , in = 0 ] ( 2 , -1 / 2 ) ;

\draw[ thick ]
    ( 0 , 1 ) arc[ radius = 1 , start angle = 90 , end angle = { 270 + 10 } ]
    to[ out = 10 , in = { 180 - 10 } ] ( 2 - 0.1736 , -0.9848 )
    arc[ radius = 1 , start angle = { -90 - 10 } , end angle = { -90 + 10 } ]
    to[ out = 10 , in = { 180 - 10 } ] ( 3.5 - 0.4 * 0.1736 , -0.9848 )
    arc[ x radius = { 0.4 * 1 } , y radius = 1 , start angle = { 270 - 10 } , end angle = { 90 + 10 } ]
    to[ out = { 180 + 10 } , in = -10 ] ( 2 + 0.1736 , 0.9848 )
    arc[ radius = 1 , start angle = { 90 - 10 } , end angle = { 90 + 10 } ]
    to[ out = { 180 + 10 } , in = -10 ] ( 0 + 0.1736 , 0.9848 )
    arc[ radius = 1 , start angle = { 90 - 10 } , end angle = 90 ]
    
    ( 6.5 , 1 ) arc[ radius = 1 , start angle = 90 , end angle = { 90 + 10 } ]
    to[ out = { 180 + 10 } , in = -10 ] ( 5 + 0.4 * 0.1736 , 0.9848 )
    arc[ x radius = 0.4 , y radius = 1 , start angle = { 90 - 10 } , end angle = { 270 + 10 } ]
    to[ out = 10 , in = { 180 - 10 } ] ( 6.5 - 0.1736 , -0.9848 )
    arc[ radius = 1 , start angle = { -90 - 10 } , end angle = { 90 + 10 } ] ;

\draw[ thick , dashed , opacity = 0.75 ]
    ( 3.5 - 0.4 * 0.1736 , -0.9848 ) arc[ x radius = 0.4 , y radius = 1 , start angle = { -90 - 10 } , end angle = { 90 + 10 } ]
    ( 5 + 0.4 * 0.1736 , -0.9848 ) arc[ x radius = 0.4 , y radius = 1 , start angle = { -90 + 10 } , end angle = { 90 - 10 } ] ;

\fill[ White ]
    ( -0.1926 , -0.0329 ) arc[ x radius = 0.2723 , y radius = 0.2 , start angle = 135 , end angle = 45 ]
    arc[ x radius = 0.35355 , y radius = 0.25 , start angle = 303 , end angle = 237 ]
    ;

\draw[ thick ]
    ( -0.25 , 0 ) arc[ x radius = 0.35355 , y radius = 0.25 , start angle = 225 , end angle = 315 ]
    ( -0.1926 , -0.0329 ) arc[ x radius = 0.2723 , y radius = 0.2 , start angle = 135 , end angle = 45 ] ;

\fill[ White ]
    ( 2 - 0.1926 , -0.0329 ) arc[ x radius = 0.2723 , y radius = 0.2 , start angle = 135 , end angle = 45 ]
    arc[ x radius = 0.35355 , y radius = 0.25 , start angle = 303 , end angle = 237 ]
    ;

\draw[ thick ]
    ( 2 - 0.25 , 0 ) arc[ x radius = 0.35355 , y radius = 0.25 , start angle = 225 , end angle = 315 ]
    ( 2 - 0.1926 , -0.0329 ) arc[ x radius = 0.2723 , y radius = 0.2 , start angle = 135 , end angle = 45 ] ;

\fill[ White ]
    ( 6.5 - 0.1926 , -0.0329 ) arc[ x radius = 0.2723 , y radius = 0.2 , start angle = 135 , end angle = 45 ]
    arc[ x radius = 0.35355 , y radius = 0.25 , start angle = 303 , end angle = 237 ]
    ;

\draw[ thick ]
    ( 6.5 - 0.25 , 0 ) arc[ x radius = 0.35355 , y radius = 0.25 , start angle = 225 , end angle = 315 ]
    ( 6.5 - 0.1926 , -0.0329 ) arc[ x radius = 0.2723 , y radius = 0.2 , start angle = 135 , end angle = 45 ] ;

\draw ( 2.75 , -2 / 3 ) node {$ S $} ;

\end{tikzpicture}
\caption{The curve $ \gamma $ represented by $ \alpha_{ 1 }^{ 2 } \alpha_{ 2 }^{ 2 } $.}\label{fig:curve2}
\end{center}
\end{figure}
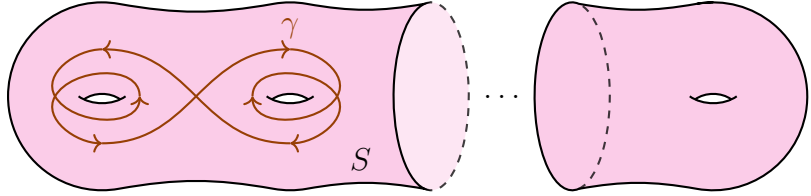

In summary, $ \gamma $ has no simple degree $ 1 $ elevations along $ p $ to $ X $, but it has one simple degree $ 1 $ elevations along $ q $ to $ Y $. $ \gamma $ has no length twins by \autoref{prop:extendedHorowitz}, so \autoref{prop:lengthSpectra} implies that $ X $ and $ Y $ are generically not simple length isospectral over $ S $.

\end{example*}

\begin{example*} \label{example:symmetric}

Let $ G \coloneqq S_{ 6 } $ be the symmetric group on $ \set{ 1 , 2 , 3 , 4 , 5 , 6 } $, and let $ S $ be a closed surface of genus $ 720 $. One can check (or refer to \cite{Gas26}) that the subgroups
\begin{align*}
A & \coloneqq \set{ \id , \of{ 1 , 2 } \of{ 3 , 4 } , \of{ 1 , 3 } \of{ 2 , 4 },\of{ 1 , 4 } \of{ 2 , 3 } } , & B & \coloneqq \set{ \id , \of{ 1 , 2 } \of{ 3 , 4 } , \of{ 1 , 2 } \of{ 5 , 6 } , \of{ 3 , 4 } \of{ 5 , 6 } },
\end{align*}
are almost conjugate. Indeed, one can check that both are subgroups isomorphic to the Klein four-group and, since the conjugacy class of an element is determined by its cycle type, these subgroups are almost conjugate. Moreover, since $ A $ acts on $ \set{ 1 , 2 , 3 , 4 , 5 , 6 } $ with two fixed points, but $ B $ does not, $ A $ and $ B $ are not conjugate.

Since there are $ 180 $ cosets for each of $ A $ and $ B $, this example differs from the previous ones in that the number of cosets makes it infeasible to give a full construction of the coset graphs.

Let $ \rho \coloneqq \of{ 1 , 5 , 3 , 6 } \of{ 2 , 4 } $ and $ \rho' \coloneqq \of{ 1 , 2 , 3 , 4 } \of{ 5 , 6 } $, and given a right coset $ A \sigma $ of $ A $, define
\begin{align*}
n_{ \rho , A \sigma } & \coloneqq \min \set{ n \in \Z : n > 0 \textrm{ and } \sigma \rho^{ n } \sigma^{ -1 } \in A } \\
n_{ \rho' , A \sigma } & \coloneqq \min \set{ n \in \Z : n > 0 \textrm{ and } \sigma \of{ \rho' }^{ n } \sigma^{ -1 } \in A } . 
\end{align*}
Analogously, given a right coset $ B \tau $ of $ B $, define
\begin{align*}
n_{ \rho , B \tau } & \coloneqq \min \set{ n \in \Z : n > 0 \textrm{ and } \tau \rho^{ n } \tau^{ -1 } \in A } \\
n_{ \rho' , B \tau } & \coloneqq \min \set{ n \in \Z : n > 0 \textrm{ and } \tau \of{ \rho' }^{ n } \tau^{ -1 } \in A } . 
\end{align*}
It is straightforward to check that these quantities are all independent of coset representative. Since $ \rho $ and $ \rho' $ have order $ 4 $, these numbers are all at most $ 4 $. Moreover, conjugacy preserves order, and since $ \rho $, $ \rho' $, $ \rho^{ -1 } = \rho^{ 3 } $, and $ \of{ \rho' }^{ -1 } = \of{ \rho' }^{ 3 } $ have order $ 4 $, they cannot be conjugate to an element of $ A $ or $ B $, all of whose elements have order $ 2 $. Thus
\[
n_{ \rho , A \sigma } , n_{ \rho' , A \sigma } , n_{ \rho , B \tau } , n_{ \rho' , B \tau } \in \set{ 2 , 4 } .
\]

\begin{claim} \label{clm:cosetsOfA}

If $ n_{ \rho , A \sigma } = 2 $, then $ n_{ \rho' , A \sigma } = 4 $, and this occurs for precisely $ 12 $ right cosets $ A \sigma $ of $ A $. Conversely, if $ n_{ \rho' , A \sigma } = 2 $, then $ n_{ \rho , A \sigma } = 4 $, and this occurs for precisely $ 12 $ right cosets $ A \sigma $ of $ A $.

\begin{proof}[Proof of \autoref{clm:cosetsOfA}]

Note that $ n_{ \rho , A \sigma } = 2 $ precisely when $ \sigma \rho^{ 2 } \sigma^{ -1 } \in A $. Similarly, $ n_{ \rho' , A \sigma } = 2 $ precisely when $ \sigma \of{ \rho' }^{ 2 } \sigma^{ -1 } \in A $. Recall that elements of $ A $ fix $ 5 $ and $ 6 $ in their action on $ \set{ 1 , \dotsc , 6 } $. Since
\[
\sigma \rho^{ 2 } \sigma^{ -1 } = \of{ \sigma \of{ 1 } , \sigma \of{ 3 } } \of{ \sigma \of{ 5 } , \sigma \of{ 6 } } , \tag{1}
\]
$ n_{ \rho , A \sigma } = 2 $ implies that $ \sigma \of{ 1 } , \sigma \of{ 3 } , \sigma \of{ 5 } , \sigma \of{ 6 } \in \set{ 1 , 2 , 3 , 4 } $. Similarly, since
\[
\sigma \of{ \rho' }^{ 2 } \sigma^{ -1 } = \of{ \sigma \of{ 1 } , \sigma \of{ 3 } } \of{ \sigma \of{ 2 } , \sigma \of{ 4 } } ,
\]
$ n_{ \rho' , A \sigma } = 2 $ implies that $ \sigma \of{ 1 } , \sigma \of{ 2 } , \sigma \of{ 3 } , \sigma \of{ 4 } \in \set{ 1 , 2 , 3 , 4 } $. Thus it cannot be that $ n_{ \rho , A \sigma } = 2 = n_{ \rho' , A \sigma } $; that is, if $ n_{ \rho , A \sigma } = 2 $, then $ n_{ \rho' , A \sigma } = 4 $, and if $ n_{ \rho' , A \sigma } = 2 $, then $ n_{ \rho , A \sigma } = 4 $.

To see that there are precisely $ 12 $ right cosets $ A \sigma $ of $ A $ for which $ n_{ \rho , A \sigma } = 2 $ and $ n_{ \rho' , A \sigma } = 4 $, note that there is a one-to-one correspondence between choices for $ \sigma \of{ 1 } $ and $ \sigma \of{ 3 } $ from $ \set{ 1 , 2 , 3 , 4 } $ and elements of $ A $. Once such a choice is made, there are two choices for $ \sigma \of{ 5 } $ from $ \set{ 1 , 2 , 3 , 4 } \setminus \set{ \sigma \of{ 1 } , \sigma \of{ 3 } } $, each of which fixes the choice of $ \sigma \of{ 6 } $. As there are two choices for $ \sigma \of{ 2 } $ from $ \set{ 5 , 6 } $, each of which fixes the choice of $ \sigma \of{ 4 } $, we find that there are $ 48 = 4 \times 3 \times 2 \times 2 $ choices of $ \sigma $, and hence $ 12 = \faktor{ 48 }{ 4 } $ right cosets $ A \sigma $.

An analogous argument shows that there are precisely $ 12 $ right cosets $ A \sigma $ of $ A $ for which $ n_{ \rho , A \sigma } = 4 $ and $ n_{ \rho' , A \sigma } = 2 $. \qedhere

\end{proof}

\end{claim}

\begin{claim} \label{clm:cosetsOfB}

$ n_{ \rho , B \tau } = n_{ \rho' , B \tau } $, and this number is $ 2 $ for precisely $ 12 $ right cosets $ B \tau $ of $ B $.

\begin{proof}[Proof of \autoref{clm:cosetsOfB}]

As argued above in \autoref{clm:cosetsOfA}, $ n_{ \rho , B \tau } = 2 $ if and only if $ \set{ \tau \of{ 1 } , \tau \of{ 3 } } $ and $ \set{ \tau \of{ 5 } , \tau \of{ 6 } } $ are a distinct pair of the sets $ \set{ 1 , 2 } $, $ \set{ 3 , 4 } $, and $ \set{ 5 , 6 } $. This occurs precisely when $ \set{ \tau \of{ 2 } , \tau \of{ 4 } } $ is one of these sets, which is exactly the condition $ \sigma \of{ \rho' }^{ 2 } \sigma^{ -1 } \in B $ and whence $ n_{ \rho' , B \tau } = 2 $.

That $ n_{ \rho , B \tau } = n_{ \rho' , B \tau } = 2 $ for precisely $ 12 $ right cosets $ B \tau $ of $ B $ now follows from the fact that $ A $ and $ B $ are almost conjugate. \qedhere

\end{proof}

\end{claim}

Given some surjective homomorphism $ \phi \colon \pi_{ 1 } \of{ S } \onto S_{ 6 } $, the above natural numbers represent the degrees of the elevations of some curve in $ \phi^{ -1 } \of{ \rho } $ (resp. $ \phi^{ -1 } \of{ \rho' } $) to each of the covers corresponding to $ \phi^{ -1 }  \of{ A } $ and $ \phi^{ -1 } \of{ B } $. There is such a homomorphism $ \phi \colon \pi_{ 1 } \of{ S } \onto S_{ 6 } $ which sends each element in a maximal non-separating multi-curve to a distinct permutation in $ S_{ 6 } $.

Concretely, fix an enumeration $ r \colon S_{ 6 } \to \set{ 1 , \dotsc , 720 } $. If we let $ \alpha_{ 1 } , \dotsc , \alpha_{ 720 } , \beta_{ 1 } , \dotsc , \beta_{ 720 } $ denote the standard generators of $ \pi_{ 1 } \of{ S } $ as drawn in \autoref{fig:generators}, and define $ \phi \colon \pi_{ 1 } \of{ S } \onto S_{ 6 } $ by the formula
\begin{align*}
\phi \of{ \alpha_{ k } } & \coloneqq r^{ -1 } \of{ k } , & \phi \of{ \beta_{ k } } & \coloneqq \id .
\end{align*}
Let $ p \colon X \to S $ and $ q \colon Y \to S $ be the covers corresponding to the subgroups $ \phi^{ -1 } \of{ A } $ and $ \phi^{ -1 } \of{ B } $, respectively. Consider the generators
\begin{align*}
\gamma & \coloneqq \alpha_{ r \of{ \rho' } } = \alpha_{ r \of{ \of{ 1 , 2 , 3 , 4 } \of{ 5 , 6 } } } & \delta & \coloneqq \alpha_{ r \of{ \rho } } = \alpha_{ r \of{ \of{ 1 , 5 , 3 , 6 } \of{ 2 , 4 } } } ,
\end{align*}
and consider the curve $ \omega $ on $ S $ represented by the element $ \gamma_{ 1 }^{ 4 } \gamma_{ 2 }^{ 2 } $, which has no length twins by \autoref{prop:extendedHorowitz}. In order to show that the covers $ p \colon X \to S $ and $ q \colon Y \to S $ are generically not simple length isospectral over $ S $, it suffices by \autoref{prop:lengthSpectra} to show that $ \omega $ has different numbers of simple elevations of degree $ 1 $ to $ X $ and $ Y $.

To that end, note that by \autoref{clm:cosetsOfA} when an elevation of $ \gamma_{ 2 } $ along $ q $ to $ Y $ has degree $ 2 $, the corresponding elevation of $ \gamma_{ 1 } $ along $ q $ to $ Y $ has degree $ 4 $, and thus the corresponding elevation of $ \omega $ will be simple and of degree $ 1 $.

On the other hand, because a given elevation of $ \gamma_{ 1 } $ along $ p $ to $ X $ will be of degree $ 2 $ if and only if the corresponding elevation of $ \gamma_{ 2 } $ is of degree $ 2 $, $ \omega $ has no simple elevations of degree $ 1 $ along $ p $ to $ X $. Thus \autoref{prop:lengthSpectra} implies that $ X $ and $ Y $ are generically not simple length isospectral over $ S $.


This construction can be easily modified to show that for each $ k \geq 1 $, the multi-sets of unmarked lengths of curves on $ X $ and $ Y $ with at most $ k $-self intersections disagree for generic hyperbolic metrics on $ S $. Indeed, consider the curve $ \omega_{ k } $ on $ S $ represented by the element $ \gamma_{ 1 }^{ 4k } \gamma_{ 2 }^{ 2 } $. This curve admits $ 12 $ elevations of degree $ 1 $ along $ q $ to $ Y $, each with self-intersection number $ k $. On the other hand, any degree $ 1 $ elevation of $ \omega_{ k } $ along $ p $ to $ X $ has at least $ k + 1 $ self intersections.  

\setcounter{claim}{0}

\subsection{Conflict of interest and data availability}

On behalf of all authors, the corresponding author states that there is no conflict of interest. Data sharing is not applicable to this article as no datasets were generated or analysed during the study.

\end{example*}

\printbibliography

\end{document}